\documentclass[11pt]{article}%
\usepackage{enumerate}
\usepackage{amsmath}
\usepackage{amssymb,latexsym}
\usepackage{amsthm}
\usepackage{mathrsfs}
\usepackage{cancel}
\usepackage{graphicx}
\usepackage{url}
\usepackage{color}
\usepackage{cite}
\usepackage{amsfonts}
\usepackage{amssymb}
\usepackage{mathrsfs}%
\providecommand{\U}[1]{\protect\rule{.1in}{.1in}}
\newcommand{\uE}{{\textup E}}
\newcommand{\uL}{{\textup L}}
\newcommand{\uS}{{\textup S}}
\newcommand{\Zp}{{\mathbb T}_p}
\newcommand{\cR}{{\mathcal R}}
\newcommand{\cT}{{\mathcal T}}
\newcommand{\cM}{{\mathcal M}}
\newcommand{\ck}{{\mathcal K}}
\newcommand{\cP}{{\mathcal P}}
\newcommand{\cL}{{\mathcal L}}
\newcommand{\cQ}{{\mathcal Q}}
\newcommand{\cX}{{\mathcal X}}
\newcommand{\F}{{\mathbb F}}
\newcommand{\K}{{\mathbb K}}
\newcommand{\GF}{\hbox{{\rm GF}}}
\newcommand{\st}{{\cal S}_t}
\newcommand{\so}{{\cal S}_1}
\newcommand{\sk}{{\cal S}_2}
\newcommand{\ad}{E}
\newcommand{\lmb}{\lambda}
\newcommand{\la}{\langle}
\newcommand{\ra}{\rangle}
\renewcommand{\proof}{\noindent{\bf Proof.}\ }
\newcommand{\Qed}{\hfill $\Box$ \medskip}
\renewcommand{\mod}{\hbox{{\rm mod}\,}}
\newtheorem{theorem}{Theorem}[section]

\newtheorem{definition}[theorem]{Definition}

\begin{document}

\title{Completeness of cubic curves in $\mathrm{PG}(2,q)$, $q\leq81$ }
\author{Daniele Bartoli
\and Stefano Marcugini
\and Fernanda Pambianco \thanks{%
Partially supported by the Italian Ministero dell'Istruzione,
dell'Universit\`a e della Ricerca (MIUR) and by the Gruppo Nazionale per le
Strutture Algebriche, Geometriche e le loro Applicazioni (GNSAGA).} \\
\\Dipartimento di Matematica e Informatica, \\Universit\`{a} degli Studi di Perugia, Italy}
\maketitle

\begin{abstract}
Theoretical results are known about the completeness of a planar algebraic
cubic curve as a $(n,3)$-arc in $\mathrm{PG}(2,q)$. They hold for $q$ big
enough and sometimes have restriction on the characteristic and on the value
of the $j$-invariant. We determine the completeness of all cubic curves for
$q\leq81$.

\end{abstract}


\section{Introduction}

In the projective plane $\mathrm{PG}(2,q)$ over the Galois field
$\mathbb{F}_{q}$ an $(n,r)$-arc is a set of $n$ points no $(r+1)$ of which are
collinear and with $r$ collinear points. An $(n,r)$-arc is called complete if
it is not contained in an $(n+1,r)$-arc of the same projective plane. An
$(n,2)$-arc is simply called an arc of size $n$ or an $n$-arc. For a detailed
description of the most important properties of these geometric structures, we
refer the reader to \cite{hirsh}. The largest size of a $(n,r)$-arc of
$\mathrm{PG}(2,q)$ is indicated by $m_{r}(2,q)$, while the minimum size of a
complete $(n,r)$-arc is indicated by $t_{r}(2,q)$. In \cite{Ball2} and
\cite{hirshStormeSurvey} bounds for $m_{r}(2,q)$ are given. In particular
$m_{3}(2,q)\leq2q+1$ for $q\geq4$; see also \cite{thas}.  The exact values of
$m_{3}(2,q)$ and $t_{3}(2,q)$ are known for $q\leq16$; see \cite{noi16} and
the references therein, whereas the complete classification of the
$(n,3)$-arcs is known for $q\leq13$; see \cite{belgi13} and the references therein.

In \cite{bier3} and \cite{hirshStormeSurvey} the relationship between the
theory of complete $(n,r)$-arcs, coding theory and mathematical statistics is presented.

A linear $[n,k,n-k]_{q}$ code $C$ such that its dual code $C^{\bot}$ has
minimum distance equal to $n$ is called NMDS. In \cite{dodu1,dodu2} it is
observed that every $[n,3,n-3]_{q}$ NMDS\ code is equivalent to an $(n,3)$-arc
in $\mathrm{PG}(2,q)$.

This means that results on the existence and completeness of $(n,3)$-arcs in
$\mathrm{PG}(2,q)$ correspond to results on the existence and
non-extendibility of $[n,3,n-3]_{q}$ NMDS\ codes, respectively.

This paper addresses the problem of the completeness of plane cubic curves as
(n,3)-arcs in $\mathrm{PG}(2,q)$. This problem is stated in \cite{hirshVoloch}%
, where the following answer is given:

\begin{theorem}
\cite[Theorem 5.1]{hirshVoloch} If $q\geq79$ and $q$ is not a power of $2$ or
$3$, then a non-singular cubic $\mathscr{C}$ with $n$ $\mathbb{F}_{q}%
$-rational points is a complete $(n,3)$-arc unless $j(\mathscr{C})=0$ in which
case the completion of $\mathscr{C}$ has at most $n+3$ points.
\end{theorem}

An extension of this result is given by \cite[Theorem 1.1, Corollary
4.2]{Giulietti2004} from which, in the case of plane cubic curves, we obtain:

\begin{theorem}
Let $q\geq121$ be an odd prime power. Let $\mathscr{C}$ be a non-singular
plane cubic curve defined over $\mathbb{F}_{q}$ whose \ $j$-invariant
$j(\mathscr{C})$ is different from $0$. Then $\mathscr{C}$ is a complete
$(n,3)$-arc.
\end{theorem}

If $j(\mathscr{C})=0$, then we have (see \cite[Theorem 4.3, Corollary
4.2]{Giulietti2004}):

\begin{theorem}
Let $q=p^{r}$, $p>3$, $q>9887$. Let $\mathscr{C}$ be a non-singular plane
cubic curve defined over $\mathbb{F}_{q}$, with $j(\mathscr{C})=0$ and having
an even number of $\mathbb{F}_{q}$-rational points. If $r$ is even or
$p\equiv1~(mod~3)$ then $\mathscr{C}$ is a complete $(n,3)$-arc.
\end{theorem}

Theoretical results on completeness of absolutely irreducible singular cubic
curves can be found in \cite{noisingolari}.

The maximum number of $\mathbb{F}_{q}$-rational points on a non-singular
algebraic curve of genus $g$, is denoted by $N_{q}(g)$. The exact value of
$N_{q}(g)$\ is not known in general, but the Hasse-Weil Bound states that
\[
q+1-2g\sqrt{q}\leq N_{q}(g)\leq q+1+2g\sqrt{q}.
\]

If a plane curve of genus $g$ has degree $n$, then $2g \leq(n-1)(n-2)$. In particular, non-singular cubic curves have $g=1$, i.e. they are elliptic curves and therefore the number of their $\mathbb{F}_{q}$-rational points  satisfies
\[
q+1-2\sqrt{q}\leq k\leq q+1+2\sqrt{q}%
\]
and can take every value in this interval other than $q+1+mp$, where $q=p^{h}$
with $p$ prime. For these results, see\cite{Hirsh24,Schoof,Hirsh2,Hirsh9}.

An absolutely irreducible singular cubic curve has $q$, $q+1$ or $q+2$,
$\mathbb{F}_{q}$-rational points depending on the fact that it contains a
node, a cusp or an isolated double point, respectively; see \cite[Table
11.7]{hirsh}.

The Hasse-Weil bound is lower than the bound on $m_{3}(2,q)$. In the cases
when the value of $m_{3}(2,q)$ is known, it is actually greater than
$q+1+2\sqrt{q}$ if $q\neq4$; see \cite{noi16}. For $q\leq16$, Table 1 gives
the values of $m_{3}(2,q)$ and of $q+1+\left\lfloor 2\sqrt{q}\right\rfloor $.
For all $q\leq81$ a non-singular cubic curve with $q+1+\left\lfloor 2\sqrt{q}\right\rfloor$
$\mathbb{F}_{q}$-rational points exists, it is complete as $(n,3)$-arc if
$q=4$ or $q\geq8$.

\bigskip\noindent\begin{minipage}{\textwidth}
\centering
\textbf{Table 1}  The values of $m_{3}(2,q)$ and $q+1+\left\lfloor 2\sqrt{q}\right\rfloor$ for  $q \leq 16$\medskip \\
\medskip
\begin{tabular}
{|c|c|c|c|c|c|c|c|c|c|c|}
\hline
$ q $ & $ 2 $ & $ 3 $ & $ 4 $ & $ 5 $ & $ 7 $ & $ 8 $ & $ 9 $ & $ 11 $ & $ 13 $ & $ 16 $\\
\hline
$m_{3}(2,q)$  & $ 7 $ & $ 9 $ & $ 9 $ & $ 11 $ & $15 $  &$ 15 $  &$ 17 $  &$ 21 $ & $ 23 $ & $ 28 $  \\
\hline
$q+1+\left\lfloor 2\sqrt{q}\right\rfloor$   & $ 5 $ & $ 7 $ & $ 9 $ & $ 10 $ & $13 $  &$ 14 $  &$ 16 $  &$ 18 $ & $ 21 $ & $ 25 $  \\
\hline
\end{tabular}
\end{minipage}

\bigskip

We follow the classification of cubic curves given in \cite{hirsh}, where two
plane curves are said to be projectively equivalent if there is a projective
transformation of the projective plane $\mathrm{PG}(2,q)$ mapping an equation
of one curve to an equation of the other.

The number of projectively distinct elliptic cubics with precisely $k$ points,
denoted by $E(k,q)$ is calculated in \cite{Schoof}.

Also a notion of isomorphism is defined for elliptic curves. We follow the
presentation given in \cite{Schoof}; see also \cite{Tate}.

\begin{definition}
Let $\mathbb{K}$ be a field; an elliptic curve $\mathscr{C}$ over $\mathbb{K}$
is a projective non-singular algebraic curve of genus 1 defined over
$\mathbb{K}$ furnished with a point 0 on $\mathscr{C}$ which is defined over
$\mathbb{K}$.
\end{definition}

Let $\overline{\mathbb{K}}$ denote an algebraic closure of $\mathbb{K}$, by
E($\overline{\mathbb{K}}$) we denote the set of points on $\mathscr{C}$
defined over $\overline{\mathbb{K}}$. This set is in a natural geometric way
an abelian group with $0$ as the zero-element. The set E($\mathbb{K}$) of
points on $\mathscr{C}$ that are defined over $\mathbb{K}$ is a subgroup of
E($\overline{\mathbb{K}}$); see \cite{Tate}.

\begin{definition}
A morphism of elliptic curves over $\mathbb{K}$, $f$: $\mathscr{C}_{1}%
\longrightarrow$\ \ $\mathscr{C}_{2}$, is an algebraic map defined over
$\mathbb{K}$ that respects the group law; in particular $f(0_{1})=0_{2}$. An
isomorphism is a morphism that has a two-sided inverse.
\end{definition}

Every elliptic curve $\mathscr{C}$ over $\mathbb{K}$ in $\mathbb{P}_{\mathbb{K}}^{2}$ is isomorphic to a curve
given by an equation:
\[
\ \ \ \ \ \ \ \ \ \ \ Y^{2}T+a_{1}XYT+a_{3}YT^{2}=X^{3}+a_{2}X^{2}%
T+a_{4}XT^{2}+a_{6}T^{3} \ \ \ \ \ \ (a_{i}\in\mathbb{K}); \ \ \ \ \ \ (1)
\]

the point $0$ is the point at infinity
$(0:1:0)$. This follows from the Riemann-Roth theorem. We have the usual formulaire:

$b_{2}=a_{1}^{2}+4a_{2},$

$b_{4}=a_{1}a_{3}+2a_{4},$

$b_{6}=a_{3}^{2}+4a_{6},$

$b_{8}=a_{1}^{2}a_{6}-a_{1}a_{3}a_{4}+4a_{2}a_{6}+a_{2}a_{3}^{2}-a_{4}^{2},$

$c_{4}=b_{2}^{2}-24b_{4}$

$c_{6}=-b_{2}^{3}+36b_{2}b_{4}-216b_{6}$

$\Delta=-b_{2}^{2}b_{8}-8b_{4}^{3}:-27b_{6}^{2}+9b_{2}b_{4}b_{6},$

$j=c_{4}^{3}/\Delta.$

A curve given by Eq. (1) is an elliptic curve if and only if the discriminant
$\Delta$ is not zero. The j-invariant of an elliptic curve $\mathscr{C}$
depends only on its isomorphism class: two elliptic curves over $\mathbb{K}$
have the same $j$-invariant if and only if they are isomorphic over
$\overline{\mathbb{K}}$. This is in general not true over $\mathbb{K}$: there
may be non-isomorphic curves over $\mathbb{K}$ that are isomorphic over
$\overline{\mathbb{K}}$. Two values of $j$ deserve special attention: they are
the values $0$ and $1728$. The elliptic curves whose $j$-invariants assume
these values correspond to the harmonic curves if $j=1728$ and to the
equianharmonic curves if $j=0$ in the sense of Hirschfeld \cite{hirsh}. If the
characteristic of $\mathbb{K}$ is $2$ or $3$ we have that $0=1728$ and the
elliptic curves over $\mathbb{K}$ with $j$-invariants equal to $0=1728$
correspond to the superharmonic curves in Hirschfeld's book.

In Section 2 we give tables summarizing the completeness of all inequivalent
absolutely irreducible singular cubic curves and non-singular cubic curves for
$q\leq81$. In Section 3 for each $q$ we give detailed tables that for each
equation present the type of cubic curve, the number of residual points, and
the $j$-invariant.

Appendix contains an Errata for \cite[Chapter 11]{hirsh}.

\section{Completeness of planar cubic curves, $q\leq81$}

In this section we present a summary of the results we obtained about the
completeness of absolutely irreducible singular cubic curves and about
non-singular cubic curves as $(n,3)$-arcs in $\mathrm{PG}(2,q)$ for
$q\leq73$ and $q=81$, cases in which Theorem 1.1 does not hold. For the sake
of completeness also the case $q=79$ is considered. The following section
gives a detailed description of such cubic curves. For the classification of
planar cubic curves see \cite{hirsh}, where all the inequivalent equations are
listed. The equations depend on the value $m=q \bmod 12$; see
\cite{hirsh}.

The results in this section have been obtained by computer. For each equation
the set $\mathcal{C}$ of $\mathbb{F}_{q}$-rational points has been computed.
Then the completeness of $\mathcal{C}$ as $(n,3)$-arc has been tested.

When necessary, the projective equivalence of couples of equations
$\mathfrak{e}_{1},$ $\mathfrak{e}_{2}$ has been verified. In this case we
proceeded in the following way.

\begin{enumerate}
\item Compute the set of $\mathbb{F}_{q}$-rational points of the two
equations, let they be $\mathcal{C}_{1},$ $\mathcal{C}_{2},$ respectively.

\item If $\mathcal{C}_{1}$ and $\mathcal{C}_{2}$ are not projectively
equivalent, then $\mathfrak{e}_{1}$ and $\mathfrak{e}_{2}$ are not projectively equivalent.

\item If $\mathcal{C}_{1}$ and $\mathcal{C}_{2}$ are projectively equivalent,
then compute a projective transformation $\pi$ sending $\mathcal{C}_{1}$ to
$\mathcal{C}_{2}$ and the stabilizer $\mathcal{S}$ of $\mathcal{C}_{2}$ in
PGL$(3,q)$.

\item Then if a projective transformation from $\mathfrak{e}_{1}$ to
$\mathfrak{e}_{2}$ exists, it can be obtained composing $\pi$\ and a
projectivity of $\mathcal{S}$.
\end{enumerate}

All computations have been performed using MAGMA; see \cite{magma}.

Table 2 summarizes the results: for $q\leq81$ the numbers of complete and
incomplete absolutely irreducible (a. i.)  singular cubic curves and non-singular cubic
curves are given.

\bigskip\noindent\begin{minipage}{\textwidth}
\centering
\textbf{Table 2}  Completeness of plane cubic curves as $(n,3)$-arcs  \medskip \\
\medskip
\begin{tabular}
{|c|c|c|c|c|c|c|c|c|c|}
\hline
$q$ &
\begin{minipage}{1 cm}
\smallskip
\centering
$\#$ of\\
inc.\\
a. i. \\
sing.\\ cubic\\ curves \\
\medskip
\end{minipage}   &
\begin{minipage}{1 cm}
\smallskip
\centering
$\#$ of\\
comp.\\
a. i. \\
sing.\\ cubic\\ curves \\
\medskip
\end{minipage}   &
\begin{minipage}{1 cm}
\smallskip
\centering
$\#$ of\\ inc.\\
non-sing.\\ cubic\\ curves \\
\medskip
\end{minipage}   &
\begin{minipage}{1 cm}
\smallskip
\centering
$\#$ of\\ comp.\\
non-sing.\\ cubic\\ curves \\
\medskip
\end{minipage}  &
$q$ &
\begin{minipage}{1 cm}
\smallskip
\centering
$\#$ of\\
inc.\\
a. i. \\
sing.\\ cubic\\ curves \\
\medskip
\end{minipage}   &
\begin{minipage}{1 cm}
\smallskip
\centering
$\#$ of\\
comp.\\
a. i. \\
sing.\\ cubic\\ curves \\
\medskip
\end{minipage}   &
\begin{minipage}{1 cm}
\smallskip
\centering
$\#$ of\\ inc.\\
non-sing.\\ cubic\\ curves \\
\medskip
\end{minipage}   &
\begin{minipage}{1 cm}
\smallskip
\centering
$\#$ of\\ comp.\\
non-sing.\\ cubic\\ curves \\
\medskip
\end{minipage}   \\
\hline
$ 2 $ & $ 4 $ & $ 0 $ & $ 6 $ & $ 0 $  &
$ 3 $ & $ 4 $ & $ 0 $ & $ 10 $ & $ 0 $  \\
$ 4 $ & $ 4 $ & $ 0 $ & $ 16 $ & $ 2 $  &
$ 5 $ & $ 4 $ & $ 0 $ & $ 15 $ & $ 1 $  \\
$ 7 $ & $ 4 $ & $ 0 $ & $ 21 $ & $ 5 $  &
$ 8 $ & $ 4 $ & $ 0 $ & $ 19 $ & $ 5 $  \\
$ 9 $ & $ 4 $ & $ 0 $ & $ 24 $ & $ 6 $  &
$ 11 $ & $ 4 $ & $ 0 $ & $ 25 $ & $ 7 $  \\
$ 13 $ & $ 4 $ & $ 0 $ & $ 32 $ & $ 14 $  &
$ 16 $ & $ 3 $ & $ 1 $ & $ 32 $ & $ 22 $  \\
$ 17 $ & $ 4 $ & $ 0 $ & $ 29 $ & $ 23 $  &
$ 19 $ & $ 4 $ & $ 0 $ & $ 32 $ & $ 30 $  \\
$ 23 $ & $ 2 $ & $ 2 $ & $ 27 $ & $ 41 $  &
$ 25 $ & $ 1 $ & $ 3 $ & $ 26 $ & $ 56 $  \\
$ 27 $ & $ 2 $ & $ 2 $ & $ 32 $ & $ 50 $  &
$ 29 $ & $ 2 $ & $ 2 $ & $ 26 $ & $ 62 $  \\
$ 31 $ & $ 1 $ & $ 3 $ & $ 27 $ & $ 71 $  &
$ 32 $ & $ 2 $ & $ 2 $ & $ 34 $ & $ 62 $  \\
$ 37 $ & $ 1 $ & $ 3 $ & $ 21 $ & $ 97 $  &
$ 41 $ & $ 2 $ & $ 2 $ & $ 16 $ & $ 108 $  \\
$ 43 $ & $ 1 $ & $ 3 $ & $ 10 $ & $ 124 $  &
$ 47 $ & $ 2 $ & $ 2 $ & $ 10 $ & $ 130 $  \\
$ 49 $ & $ 1 $ & $ 3 $ & $ 8 $ & $ 146 $  &
$ 53 $ & $ 2 $ & $ 2 $ & $ 5 $ & $ 155 $  \\
$ 59 $ & $ 2 $ & $ 2 $ & $ 5 $ & $ 171 $  &
$ 61 $ & $ 1 $ & $ 3 $ & $ 2 $ & $ 188 $  \\
$ 64 $ & $ 1 $ & $ 3 $ & $ 0 $ & $ 198 $  &
$ 67 $ & $ 1 $ & $ 3 $ & $ 0 $ & $ 206 $  \\
$ 71 $ & $ 2 $ & $ 2 $ & $ 3 $ & $ 209 $  &
$ 73 $ & $ 1 $ & $ 3 $ & $ 1 $ & $ 225 $  \\
$ 79 $ & $ 1 $ & $ 3 $ & $ 0 $ & $ 242 $  &
$ 81 $ & $ 2 $ & $ 2 $ & $ 2 $ & $ 244 $  \\
\hline
\end{tabular}
\end{minipage}

\bigskip

In the following tables the symbol $N^{c}$ denotes $c$ complete cubic curves
having $N$ $\mathbb{F}_{q}$-rational points, while the symbol $N_{r}^{c}$
denotes $c$ incomplete cubic curves having $N$ $\mathbb{F}_{q}$-rational
points and $r$ residual points, i.e. $r$ points of $\mathrm{PG}(2,q)$ not
belonging to the $\mathbb{F}_{q}$-rational points of the cubic curve and not
lying on trisecant lines of the $\mathbb{F}_{q}$-rational points of the cubic curve.

Table 3 gives detailed results for absolutely irreducible singular cubic
curves. Table 4 and Table 5 give  detailed results for incomplete non-singular cubic curves and for complete non-singular cubic curves, respectively.

\bigskip


\noindent\begin{minipage}{\textwidth}
\centering
\textbf{Table 3} Classification of absolutely irreducible (a. i.) singular cubic curves  \medskip \\
\medskip
\begin{tabular}
{|c|c|c|c|c|c|}
\hline
\begin{minipage}{0.6 cm}
\centering
\smallskip
$q $
\end{minipage}
&
\begin{minipage}{3.3 cm}
\smallskip
\centering
$\#$ of a. i. singular\\ cubic curves \\
\medskip
\end{minipage}   &
\begin{minipage}{0.6 cm}
\centering
\smallskip
$q $
\end{minipage}
&
\begin{minipage}{3.3 cm}
\smallskip
\centering
$\#$ of a. i. singular\\ cubic curves \\
\medskip
\end{minipage}   &
\begin{minipage}{0.6 cm}
\centering
\smallskip
$q $
\end{minipage}
&
\begin{minipage}{3.3 cm}
\smallskip
\centering
$\#$ of a. i. singular\\ cubic curves \\
\medskip
\end{minipage}     \\
\hline
\begin{minipage}{0.6 cm}
\centering
\smallskip
$2 $
\end{minipage}
&
\begin{minipage}{3.3 cm}
\centering
\smallskip
$ 2^{1}_{5}, 3^{1}_{4}, 4^{2}_{3}$
\end{minipage} &
\begin{minipage}{0.6 cm}
\centering
\smallskip
$3 $
\end{minipage}
&
\begin{minipage}{3.3 cm}
\centering
\smallskip
$ 3^{1}_{10}, 4^{1}_{8}, 4^{1}_{9}, 5^{1}_{7}$
\end{minipage}  &
\begin{minipage}{0.6 cm}
\centering
\smallskip
$4 $
\end{minipage}
&
\begin{minipage}{3.3 cm}
\centering
\smallskip
$ 4^{1}_{15}, 4^{1}_{17}, 5^{1}_{14}, 6^{1}_{11}$
\end{minipage}  \\
\begin{minipage}{0.6 cm}
\centering
\smallskip
$5 $
\end{minipage}
&
\begin{minipage}{3.3 cm}
\centering
\smallskip
$ 5^{1}_{23}, 6^{1}_{19}, 7^{1}_{12}, 7^{1}_{15}$
\end{minipage}  &
\begin{minipage}{0.6 cm}
\centering
\smallskip
$7 $
\end{minipage}
&
\begin{minipage}{3.3 cm}
\centering
\smallskip
$ 7^{1}_{30}, 7^{1}_{35}, 8^{1}_{25}, 9^{1}_{16}$
\end{minipage}   &
\begin{minipage}{0.6 cm}
\centering
\smallskip
$8 $
\end{minipage}
&
\begin{minipage}{3.3 cm}
\centering
\smallskip
$ 8^{1}_{36}, 9^{1}_{22}, 10^{1}_{12}, 10^{1}_{18}$
\end{minipage}  \\
\begin{minipage}{0.6 cm}
\centering
\smallskip
$9 $
\end{minipage}
&
\begin{minipage}{3.3 cm}
\centering
\smallskip
$ 9^{1}_{36}, 10^{1}_{5}, 10^{1}_{27}, 11^{1}_{14}$
\end{minipage}  &
\begin{minipage}{0.6 cm}
\centering
\smallskip
$11 $
\end{minipage}
&
\begin{minipage}{3.3 cm}
\centering
\smallskip
$ 11^{1}_{31}, 12^{1}_{21}, 13^{1}_{9}, 13^{1}_{12}$
\end{minipage}   &
\begin{minipage}{0.6 cm}
\centering
\smallskip
$13 $
\end{minipage}
&
\begin{minipage}{3.3 cm}
\centering
\smallskip
$ 13^{1}_{21}, 13^{1}_{28}, 14^{1}_{16}, 15^{1}_{5}$
\end{minipage}  \\
\begin{minipage}{0.6 cm}
\centering
\smallskip
$16 $
\end{minipage}
&
\begin{minipage}{3.3 cm}
\centering
\smallskip
$ 16^{1}_{12}, 16^{1}_{15}, 17^{1}, 18^{1}_{1}$
\end{minipage}  &
\begin{minipage}{0.6 cm}
\centering
\smallskip
$17 $
\end{minipage}
&
\begin{minipage}{3.3 cm}
\centering
\smallskip
$ 17^{1}_{8}, 18^{1}_{17}, 19^{2}_{6}$
\end{minipage}   &
\begin{minipage}{0.6 cm}
\centering
\smallskip
$19 $
\end{minipage}
&
\begin{minipage}{3.3 cm}
\centering
\smallskip
$ 19^{1}_{3}, 19^{1}_{4}, 20^{1}_{6}, 21^{1}_{3}$
\end{minipage}  \\
\begin{minipage}{0.6 cm}
\centering
\smallskip
$23 $
\end{minipage}
&
\begin{minipage}{3.3 cm}
\centering
\smallskip
$ 23^{1}, 24^{1}_{1}, 25^{1}, 25^{1}_{3}$
\end{minipage}   &
\begin{minipage}{0.6 cm}
\centering
\smallskip
$25 $
\end{minipage}
&
\begin{minipage}{3.3 cm}
\centering
\smallskip
$ 25^{2}, 26^{1}, 27^{1}_{1}$
\end{minipage}   &
\begin{minipage}{0.6 cm}
\centering
\smallskip
$27 $
\end{minipage}
&
\begin{minipage}{3.3 cm}
\centering
\smallskip
$ 27^{1}, 28^{1}, 28^{1}_{14}, 29^{1}_{1}$
\end{minipage}  \\
\begin{minipage}{0.6 cm}
\centering
\smallskip
$29 $
\end{minipage}
&
\begin{minipage}{3.3 cm}
\centering
\smallskip
$ 29^{1}, 30^{1}_{1}, 31^{1}, 31^{1}_{3}$
\end{minipage}  &
\begin{minipage}{0.6 cm}
\centering
\smallskip
$31 $
\end{minipage}
&
\begin{minipage}{3.3 cm}
\centering
\smallskip
$ 31^{2}, 32^{1}, 33^{1}_{1}$
\end{minipage}  &
\begin{minipage}{0.6 cm}
\centering
\smallskip
$32 $
\end{minipage}
&
\begin{minipage}{3.3 cm}
\centering
\smallskip
$ 32^{1}, 33^{1}_{1}, 34^{1}, 34^{1}_{3}$
\end{minipage}  \\
\begin{minipage}{0.6 cm}
\centering
\smallskip
$37 $
\end{minipage}
&
\begin{minipage}{3.3 cm}
\centering
\smallskip
$ 37^{2}, 38^{1}, 39^{1}_{1}$
\end{minipage}  &
\begin{minipage}{0.6 cm}
\centering
\smallskip
$41 $
\end{minipage}
&
\begin{minipage}{3.3 cm}
\centering
\smallskip
$ 41^{1}, 42^{1}_{1}, 43^{1}, 43^{1}_{3}$
\end{minipage}  &
\begin{minipage}{0.6 cm}
\centering
\smallskip
$43 $
\end{minipage}
&
\begin{minipage}{3.3 cm}
\centering
\smallskip
$ 43^{2}, 44^{1}, 45^{1}_{1}$
\end{minipage}  \\
\begin{minipage}{0.6 cm}
\centering
\smallskip
$47 $
\end{minipage}
&
\begin{minipage}{3.3 cm}
\centering
\smallskip
$ 47^{1}, 48^{1}_{1}, 49^{1}, 49^{1}_{3}$
\end{minipage}  &
\begin{minipage}{0.6 cm}
\centering
\smallskip
$49 $
\end{minipage}
&
\begin{minipage}{3.3 cm}
\centering
\smallskip
$ 49^{2}, 50^{1}, 51^{1}_{1}$
\end{minipage}   &
\begin{minipage}{0.6 cm}
\centering
\smallskip
$53 $
\end{minipage}
&
\begin{minipage}{3.3 cm}
\centering
\smallskip
$ 53^{1}, 54^{1}_{1}, 55^{1}, 55^{1}_{3}$
\end{minipage}  \\
\begin{minipage}{0.6 cm}
\centering
\smallskip
$59 $
\end{minipage}
&
\begin{minipage}{3.3 cm}
\centering
\smallskip
$ 59^{1}, 60^{1}_{1}, 61^{1}, 61^{1}_{3}$
\end{minipage}   &
\begin{minipage}{0.6 cm}
\centering
\smallskip
$61 $
\end{minipage}
&
\begin{minipage}{3.3 cm}
\centering
\smallskip
$ 61^{2}, 62^{1}, 63^{1}_{1}$
\end{minipage}  &
\begin{minipage}{0.6 cm}
\centering
\smallskip
$64 $
\end{minipage}
&
\begin{minipage}{3.3 cm}
\centering
\smallskip
$ 64^{2}, 65^{1}, 66^{1}_{1}$
\end{minipage}  \\
\begin{minipage}{0.6 cm}
\centering
\smallskip
$67 $
\end{minipage}
&
\begin{minipage}{3.3 cm}
\centering
\smallskip
$ 67^{2}, 68^{1}, 69^{1}_{1}$
\end{minipage}   &
\begin{minipage}{0.6 cm}
\centering
\smallskip
$71 $
\end{minipage}
&
\begin{minipage}{3.3 cm}
\centering
\smallskip
$ 71^{1}, 72^{1}_{1}, 73^{1}, 73^{1}_{3}$
\end{minipage}  &
\begin{minipage}{0.6 cm}
\centering
\smallskip
$73 $
\end{minipage}
&
\begin{minipage}{3.3 cm}
\centering
\smallskip
$ 73^{2}, 74^{1}, 75^{1}_{1}$
\end{minipage}  \\
\begin{minipage}{0.6 cm}
\centering
\smallskip
$79 $
\end{minipage}
&
\begin{minipage}{3.3 cm}
\centering
\smallskip
$ 79^{2}, 80^{1}, 81^{1}_{1}$
\end{minipage}  &
\begin{minipage}{0.6 cm}
\centering
\smallskip
$81 $
\end{minipage}
&
\begin{minipage}{3.3 cm}
\centering
\smallskip
$ 81^{1}, 82^{1}, 82^{1}_{41}, 83^{1}_{1}$
\smallskip
\end{minipage}  &
&
\\
\hline
\end{tabular}
\end{minipage}\bigskip



\noindent\begin{minipage}{\textwidth}
\centering
\textbf{Table 4}  Classification of incomplete non-singular cubic curves \medskip \\
\medskip
\begin{tabular}
{|c|c|}
\hline
$q$ &
\begin{minipage}{8  cm}
\smallskip
\centering
incomplete non-singular cubic curves \\
\medskip
\end{minipage}   \\
\hline
\begin{minipage}{1 cm}
\centering
\smallskip
$2 $
\end{minipage}
&
\begin{minipage}{11 cm}
\centering
\smallskip
$ 1^{1}_{6}, 2^{1}_{5}, 3^{2}_{4}, 4^{1}_{3}, 5^{1}_{2}$
\end{minipage}  \\
\begin{minipage}{1 cm}
\centering
\smallskip
$3 $
\end{minipage}
&
\begin{minipage}{11 cm}
\centering
\smallskip
$ 1^{1}_{12}, 2^{1}_{11}, 3^{1}_{9}, 3^{1}_{10}, 4^{2}_{8}, 5^{1}_{6}, 6^{1}_{3}, 6^{1}_{4}, 7^{1}_{2}$
\end{minipage}  \\
\begin{minipage}{1 cm}
\centering
\smallskip
$4 $
\end{minipage}
&
\begin{minipage}{11 cm}
\centering
\smallskip
$ 1^{1}_{20}, 2^{1}_{19}, 3^{2}_{16}, 3^{2}_{18}, 4^{2}_{15}, 5^{1}_{12}, 6^{2}_{7}, 6^{2}_{9}, 7^{2}_{5}, 8^{1}_{2}$
\end{minipage}  \\
\begin{minipage}{1 cm}
\centering
\smallskip
$5 $
\end{minipage}
&
\begin{minipage}{11 cm}
\centering
\smallskip
$ 2^{1}_{29}, 3^{1}_{25}, 3^{1}_{28}, 4^{2}_{24}, 5^{1}_{20}, 6^{2}_{13}, 6^{2}_{16}, 7^{1}_{10}, 8^{1}_{4}, 8^{1}_{5}, 9^{1}_{4}, 10^{1}_{1}$
\end{minipage}  \\
\begin{minipage}{1 cm}
\centering
\smallskip
$7 $
\end{minipage}
&
\begin{minipage}{11 cm}
\centering
\smallskip
$ 3^{1}_{49}, 3^{1}_{54}, 4^{2}_{48}, 5^{1}_{42}, 6^{2}_{31}, 6^{2}_{36}, 7^{2}_{26}, 8^{1}_{16}, 8^{1}_{17}, 9^{1}_{6}, 9^{1}_{9}, 9^{2}_{12}, 10^{1}_{4}, 10^{1}_{5}, 11^{1}_{3}, 13^{1}_{1}$
\end{minipage}  \\
\begin{minipage}{1 cm}
\centering
\smallskip
$8 $
\end{minipage}
&
\begin{minipage}{11 cm}
\centering
\smallskip
$ 4^{1}_{63}, 5^{1}_{56}, 6^{3}_{43}, 6^{3}_{49}, 8^{3}_{26}, 9^{1}_{13}, 9^{1}_{19}, 10^{3}_{9}, 12^{3}_{4}$
\end{minipage}  \\
\begin{minipage}{1 cm}
\centering
\smallskip
$9 $
\end{minipage}
&
\begin{minipage}{11 cm}
\centering
\smallskip
$ 4^{1}_{80}, 5^{1}_{72}, 6^{2}_{57}, 6^{2}_{64}, 7^{1}_{50}, 8^{1}_{36}, 8^{2}_{37}, 9^{2}_{21}, 9^{2}_{28}, 10^{2}_{15}, 11^{2}_{8}, 12^{1}_{1}, 12^{4}_{3}, 13^{1}_{7}$
\end{minipage}  \\
\begin{minipage}{1 cm}
\centering
\smallskip
$11 $
\end{minipage}
&
\begin{minipage}{11 cm}
\centering
\smallskip
$ 6^{1}_{91}, 6^{1}_{100}, 7^{1}_{82}, 8^{1}_{64}, 8^{1}_{65}, 9^{1}_{42}, 9^{1}_{43}, 9^{2}_{52}, 10^{2}_{33}, 11^{1}_{20}, 12^{1}_{6}, 12^{1}_{7}, 12^{1}_{9}, 12^{2}_{10},$ \\
$ 12^{1}_{12}, 12^{2}_{13}, 13^{1}_{5}, 14^{1}_{1}, 14^{1}_{3}, 15^{1}_{1}, 15^{1}_{3}$
\end{minipage}  \\
\begin{minipage}{1 cm}
\centering
\smallskip
$13 $
\end{minipage}
&
\begin{minipage}{11 cm}
\centering
\smallskip
$ 7^{1}_{122}, 8^{1}_{100}, 8^{1}_{101}, 9^{1}_{54}, 9^{1}_{73}, 9^{1}_{81}, 9^{2}_{84}, 10^{1}_{58}, 10^{2}_{59}, 11^{1}_{42}, 12^{1}_{12}, 12^{1}_{19},$\\
$ 12^{1}_{21}, 12^{3}_{24}, 12^{2}_{27}, 13^{1}_{14}, 13^{1}_{17}, 14^{1}_{6}, 14^{1}_{8}, 15^{3}_{3}, 15^{1}_{6}, 16^{2}_{2}, 17^{1}_{1}, 18^{1}_{1}$
\end{minipage}  \\
\begin{minipage}{1 cm}
\centering
\smallskip
$16 $
\end{minipage}
&
\begin{minipage}{11 cm}
\centering
\smallskip
$ 9^{1}_{108}, 9^{1}_{144}, 10^{2}_{113}, 12^{4}_{58}, 12^{4}_{63}, 13^{2}_{41}, 14^{4}_{26}, 16^{1}_{4}, 16^{4}_{9}, 18^{4}_{4}, 18^{5}_{9}$
\end{minipage}  \\
\begin{minipage}{1 cm}
\centering
\smallskip
$17 $
\end{minipage}
&
\begin{minipage}{11 cm}
\centering
\smallskip
$ 10^{1}_{135}, 11^{1}_{105}, 12^{1}_{70}, 12^{1}_{72}, 12^{1}_{73}, 12^{1}_{76}, 12^{1}_{79}, 12^{1}_{82}, 13^{1}_{53}, 14^{1}_{31}, 14^{1}_{35}, 15^{1}_{15},$\\
$ 15^{2}_{21}, 15^{2}_{25}, 15^{1}_{28}, 16^{1}_{4}, 16^{1}_{6}, 16^{1}_{13}, 16^{1}_{16}, 17^{1}_{3}, 18^{2}_{4}, 18^{2}_{6}, 19^{1}_{2}, 20^{1}_{1}, 24^{1}_{1}$
\end{minipage}  \\
\begin{minipage}{1 cm}
\centering
\smallskip
$19 $
\end{minipage}
&
\begin{minipage}{11 cm}
\centering
\smallskip
$ 12^{1}_{96}, 12^{1}_{103}, 12^{1}_{117}, 12^{1}_{120}, 13^{1}_{83}, 13^{1}_{86}, 14^{1}_{56}, 14^{1}_{61}, 15^{1}_{28}, 15^{2}_{39}, 15^{1}_{42}, 16^{3}_{22},$\\
$16^{1}_{25}, 17^{1}_{15}, 18^{2}_{3}, 18^{1}_{6}, 18^{1}_{7}, 18^{3}_{9}, 19^{1}_{2}, 19^{2}_{3}, 20^{1}_{1}, 20^{1}_{2}, 20^{1}_{3}, 22^{2}_{1}$
\end{minipage}  \\
\begin{minipage}{1 cm}
\centering
\smallskip
$23 $
\end{minipage}
&
\begin{minipage}{11 cm}
\centering
\smallskip
$ 15^{1}_{87}, 15^{1}_{100}, 16^{1}_{64}, 16^{1}_{68}, 17^{1}_{40}, 18^{1}_{16}, 18^{1}_{24}, 18^{2}_{25}, 18^{1}_{27}, 18^{1}_{31}, 18^{2}_{33}, 19^{1}_{12},$\\
$ 20^{1}_{6}, 20^{2}_{10}, 20^{1}_{12}, 21^{1}_{3}, 21^{1}_{4}, 21^{1}_{7}, 21^{1}_{9}, 22^{1}_{2}, 22^{1}_{3}, 24^{2}_{1}, 24^{1}_{3}$
\end{minipage}  \\
\begin{minipage}{1 cm}
\centering
\smallskip
$25 $
\end{minipage}
&
\begin{minipage}{11 cm}
\centering
\smallskip
$ 16^{1}_{96}, 17^{1}_{71}, 18^{4}_{39}, 18^{2}_{42}, 18^{1}_{45}, 19^{2}_{33}, 20^{2}_{16}, 20^{2}_{20}, 21^{2}_{21},$\\
$ 22^{2}_{3}, 22^{2}_{5}, 23^{2}_{5}, 24^{2}_{4}, 25^{1}_{3}$
\end{minipage}  \\
\begin{minipage}{1 cm}
\centering
\smallskip
$27 $
\end{minipage}
&
\begin{minipage}{11 cm}
\centering
\smallskip
$ 18^{1}_{54}, 18^{1}_{73}, 19^{1}_{59}, 20^{3}_{30}, 20^{1}_{32}, 21^{3}_{15}, 21^{3}_{18}, 23^{3}_{8}, 24^{9}_{3}, 26^{3}_{1}, 26^{3}_{2}, 28^{1}_{1}$
\end{minipage}  \\
\begin{minipage}{1 cm}
\centering
\smallskip
$29 $
\end{minipage}
&
\begin{minipage}{11 cm}
\centering
\smallskip
$ 20^{1}_{40}, 20^{1}_{55}, 21^{2}_{30}, 21^{1}_{31}, 21^{1}_{34}, 22^{1}_{13}, 22^{1}_{18}, 23^{1}_{16}, 24^{1}_{1}, 24^{1}_{4}, 24^{3}_{6}, 24^{2}_{9}, 24^{1}_{12},$\\
$ 25^{1}_{1}, 25^{1}_{5}, 26^{1}_{2}, 26^{1}_{3}, 27^{1}_{3}, 30^{2}_{1}, 30^{2}_{4}$
\end{minipage}  \\
\begin{minipage}{1 cm}
\centering
\smallskip
$31 $
\end{minipage}
&
\begin{minipage}{11 cm}
\centering
\smallskip
$ 21^{1}_{27}, 21^{1}_{54}, 22^{1}_{29}, 22^{1}_{30}, 23^{1}_{16}, 24^{1}_{6}, 24^{2}_{9}, 24^{1}_{16}, 24^{1}_{18}, 24^{1}_{25}, 25^{1}_{7}, 25^{1}_{8}, 25^{1}_{9},$\\
$ 26^{1}_{2}, 26^{1}_{4}, 27^{2}_{3}, 27^{1}_{6}, 28^{2}_{1}, 28^{2}_{2}, 29^{1}_{1}, 30^{1}_{3}, 30^{1}_{6}, 33^{1}_{1}$
\end{minipage}  \\
\begin{minipage}{1 cm}
\centering
\smallskip
$32 $
\end{minipage}
&
\begin{minipage}{11 cm}
\centering
\smallskip
$ 22^{1}_{33}, 24^{5}_{9}, 24^{5}_{28}, 25^{1}_{2}, 26^{5}_{5}, 28^{5}_{2}, 30^{5}_{1}, 30^{5}_{3}, 33^{1}_{1}, 33^{1}_{4}$
\end{minipage}  \\
\begin{minipage}{1 cm}
\centering
\smallskip
$37 $
\end{minipage}
&
\begin{minipage}{11 cm}
\centering
\smallskip
$ 26^{1}_{22}, 27^{1}_{9}, 27^{1}_{16}, 27^{2}_{18}, 28^{1}_{2}, 28^{1}_{4}, 28^{1}_{5}, 28^{1}_{9}, 29^{1}_{1},$\\ $30^{5}_{3}, 30^{2}_{4}, 31^{1}_{1}, 32^{1}_{2}, 33^{1}_{3}, 34^{1}_{1}$
\end{minipage}  \\
\begin{minipage}{1 cm}
\centering
\smallskip
$41 $
\end{minipage}
&
\begin{minipage}{11 cm}
\centering
\smallskip
$ 30^{1}_{3}, 30^{1}_{6}, 30^{1}_{12}, 30^{1}_{16}, 31^{1}_{1}, 32^{2}_{2}, 32^{1}_{3}, 32^{1}_{4}, 33^{1}_{1}, 34^{1}_{1}, 35^{1}_{1}, 36^{1}_{4}, 42^{2}_{1}, 42^{1}_{3}$
\end{minipage}  \\
\begin{minipage}{1 cm}
\centering
\smallskip
$43 $
\end{minipage}
&
\begin{minipage}{11 cm}
\centering
\smallskip
$ 31^{1}_{6}, 32^{1}_{5}, 33^{1}_{6}, 34^{2}_{1}, 34^{1}_{2}, 35^{2}_{1}, 37^{1}_{2}, 38^{1}_{2}$
\end{minipage}  \\[1 mm]
\hline
\end{tabular}
\end{minipage}



\noindent\begin{minipage}{\textwidth}
\centering
\textbf{Table 4}  Continue  \medskip \\
\medskip
\begin{tabular}
{|c|c|}
\hline
$q$ &
\begin{minipage}{8  cm}
\smallskip
\centering
incomplete non-singular cubic curves \\
\medskip
\end{minipage}   \\
\hline
\begin{minipage}{1 cm}
\centering
\smallskip
$47 $
\end{minipage}
&
\begin{minipage}{11 cm}
\centering
\smallskip
$ 35^{1}_{1}, 36^{3}_{3}, 36^{1}_{6}, 40^{1}_{1}, 45^{1}_{1}, 48^{2}_{1}, 48^{1}_{3}$
\end{minipage}  \\
\begin{minipage}{1 cm}
\centering
\smallskip
$49 $
\end{minipage}
&
\begin{minipage}{11 cm}
\centering
\smallskip
$ 37^{1}_{3}, 37^{1}_{5}, 40^{4}_{1}, 41^{2}_{1}$
\end{minipage}  \\
\begin{minipage}{1 cm}
\centering
\smallskip
$53 $
\end{minipage}
&
\begin{minipage}{11 cm}
\centering
\smallskip
$ 42^{1}_{1}, 43^{1}_{2}, 54^{2}_{1}, 54^{1}_{3}$
\end{minipage}  \\
\begin{minipage}{1 cm}
\centering
\smallskip
$59 $
\end{minipage}
&
\begin{minipage}{11 cm}
\centering
\smallskip
$ 45^{1}_{3}, 47^{1}_{1}, 60^{2}_{1}, 60^{1}_{3}$
\end{minipage}  \\
\begin{minipage}{1 cm}
\centering
\smallskip
$61 $
\end{minipage}
&
\begin{minipage}{11 cm}
\centering
\smallskip
$ 48^{1}_{3}, 54^{1}_{1}$
\end{minipage}  \\
\begin{minipage}{1 cm}
\centering
\smallskip
$64 $
\end{minipage}
&
\begin{minipage}{11 cm}
\centering
\smallskip
$0$
\end{minipage}  \\
\begin{minipage}{1 cm}
\centering
\smallskip
$67 $
\end{minipage}
&
\begin{minipage}{11 cm}
\centering
\smallskip
$0$
\end{minipage}  \\
\begin{minipage}{1 cm}
\centering
\smallskip
$71 $
\end{minipage}
&
\begin{minipage}{11 cm}
\centering
\smallskip
$ 72^{2}_{1}, 72^{1}_{3}$
\end{minipage}  \\
\begin{minipage}{1 cm}
\centering
\smallskip
$73 $
\end{minipage}
&
\begin{minipage}{11 cm}
\centering
\smallskip
$ 60^{1}_{1}$
\end{minipage}  \\
\begin{minipage}{1 cm}
\centering
\smallskip
$79 $
\end{minipage}
&
\begin{minipage}{11 cm}
\centering
\smallskip
$0$
\end{minipage}  \\
\begin{minipage}{1 cm}
\centering
\smallskip
$81 $
\end{minipage}
&
\begin{minipage}{11 cm}
\centering
\smallskip
$ 82^{2}_{1}$
\end{minipage}  \\[1 mm]
\hline
\end{tabular}
\end{minipage}

\bigskip



\noindent\begin{minipage}{\textwidth}
\centering
\textbf{Table 5}  Classification of complete non-singular cubic curves  \medskip \\
\medskip
\begin{tabular}
{|c|c|}
\hline
$q$ &
\begin{minipage}{8  cm}
\smallskip
\centering
complete non-singular cubic curves \\
\medskip
\end{minipage}   \\
\hline
\begin{minipage}{1 cm}
\centering
\smallskip
$2 $
\end{minipage}
&
\begin{minipage}{11 cm}
\centering
\smallskip
$0$
\end{minipage}  \\
\begin{minipage}{1 cm}
\centering
\smallskip
$3 $
\end{minipage}
&
\begin{minipage}{11 cm}
\centering
\smallskip
$0$
\end{minipage}  \\
\begin{minipage}{1 cm}
\centering
\smallskip
$4 $
\end{minipage}
&
\begin{minipage}{11 cm}
\centering
\smallskip
$ 9^{2}$
\end{minipage}  \\
\begin{minipage}{1 cm}
\centering
\smallskip
$5 $
\end{minipage}
&
\begin{minipage}{11 cm}
\centering
\smallskip
$ 9^{1}$
\end{minipage}  \\
\begin{minipage}{1 cm}
\centering
\smallskip
$7 $
\end{minipage}
&
\begin{minipage}{11 cm}
\centering
\smallskip
$ 9^{1}, 12^{4}$
\end{minipage}  \\
\begin{minipage}{1 cm}
\centering
\smallskip
$8 $
\end{minipage}
&
\begin{minipage}{11 cm}
\centering
\smallskip
$ 12^{3}, 13^{1}, 14^{1}$
\end{minipage}  \\
\begin{minipage}{1 cm}
\centering
\smallskip
$9 $
\end{minipage}
&
\begin{minipage}{11 cm}
\centering
\smallskip
$ 12^{1}, 14^{2}, 15^{2}, 16^{1}$
\end{minipage}  \\
\begin{minipage}{1 cm}
\centering
\smallskip
$11 $
\end{minipage}
&
\begin{minipage}{11 cm}
\centering
\smallskip
$ 15^{2}, 16^{2}, 17^{1}, 18^{2}$
\end{minipage}  \\
\begin{minipage}{1 cm}
\centering
\smallskip
$13 $
\end{minipage}
&
\begin{minipage}{11 cm}
\centering
\smallskip
$ 16^{2}, 18^{6}, 19^{2}, 20^{2}, 21^{2}$
\end{minipage}  \\
\begin{minipage}{1 cm}
\centering
\smallskip
$16 $
\end{minipage}
&
\begin{minipage}{11 cm}
\centering
\smallskip
$ 17^{1}, 18^{4}, 20^{4}, 21^{4}, 22^{4}, 24^{4}, 25^{1}$
\end{minipage}  \\
\begin{minipage}{1 cm}
\centering
\smallskip
$17 $
\end{minipage}
&
\begin{minipage}{11 cm}
\centering
\smallskip
$ 18^{4}, 20^{3}, 21^{6}, 22^{2}, 23^{1}, 24^{5}, 25^{1}, 26^{1}$
\end{minipage}  \\
\begin{minipage}{1 cm}
\centering
\smallskip
$19 $
\end{minipage}
&
\begin{minipage}{11 cm}
\centering
\smallskip
$ 18^{2}, 20^{1}, 21^{6}, 22^{1}, 23^{1}, 24^{8}, 25^{2}, 26^{2}, 27^{5}, 28^{2}$
\end{minipage}  \\
\begin{minipage}{1 cm}
\centering
\smallskip
$23 $
\end{minipage}
&
\begin{minipage}{11 cm}
\centering
\smallskip
$ 21^{2}, 23^{2}, 24^{9}, 25^{2}, 26^{2}, 27^{6}, 28^{4}, 29^{1}, 30^{8}, 31^{1}, 32^{2}, 33^{2}$
\end{minipage}  \\
\begin{minipage}{1 cm}
\centering
\smallskip
$25 $
\end{minipage}
&
\begin{minipage}{11 cm}
\centering
\smallskip
$ 21^{2}, 24^{10}, 25^{2}, 27^{9}, 28^{6}, 29^{2}, 30^{8}, 31^{2}, 32^{4}, 33^{4}, 34^{3}, 35^{1}, 36^{3}$
\end{minipage}  \\
\begin{minipage}{1 cm}
\centering
\smallskip
$27 $
\end{minipage}
&
\begin{minipage}{11 cm}
\centering
\smallskip
$ 24^{3}, 27^{6}, 28^{1}, 29^{3}, 30^{12}, 32^{6}, 33^{6}, 35^{3}, 36^{8}, 37^{1}, 38^{1}$
\end{minipage}  \\
\begin{minipage}{1 cm}
\centering
\smallskip
$29 $
\end{minipage}
&
\begin{minipage}{11 cm}
\centering
\smallskip
$ 24^{4}, 26^{1}, 27^{5}, 28^{4}, 29^{2}, 30^{8}, 31^{2}, 32^{4}, 33^{6}, 34^{3}, 35^{2}, 36^{12}, 37^{1}, 38^{2},$\\ $39^{4}, 40^{2}$
\end{minipage}  \\
\begin{minipage}{1 cm}
\centering
\smallskip
$31 $
\end{minipage}
&
\begin{minipage}{11 cm}
\centering
\smallskip
$ 24^{2}, 27^{6}, 28^{2}, 29^{1}, 30^{6}, 31^{2}, 32^{6}, 33^{3}, 34^{4}, 35^{2}, 36^{16}, 37^{3}, 38^{2}, 39^{6},$\\ $40^{4}, 41^{1}, 42^{4}, 43^{1}$
\end{minipage}  \\
\hline
\end{tabular}
\end{minipage}\bigskip



\noindent\begin{minipage}{\textwidth}
\centering
\textbf{Table 5}  Continue  \medskip \\
\medskip
\begin{tabular}
{|c|c|}
\hline
$q$ &
\begin{minipage}{8  cm}
\smallskip
\centering
complete non-singular cubic curves \\
\medskip
\end{minipage}   \\
\hline
\begin{minipage}{1 cm}
\centering
\smallskip
$32 $
\end{minipage}
&
\begin{minipage}{11 cm}
\centering
\smallskip
$ 30^{10}, 32^{5}, 34^{5}, 36^{20}, 38^{5}, 40^{5}, 41^{1}, 42^{10}, 44^{1}$
\end{minipage}  \\
\begin{minipage}{1 cm}
\centering
\smallskip
$37 $
\end{minipage}
&
\begin{minipage}{11 cm}
\centering
\smallskip
$ 27^{1}, 30^{1}, 31^{2}, 32^{3}, 33^{3}, 34^{3}, 35^{3}, 36^{20}, 37^{3}, 38^{2}, 39^{6}, 40^{8}, 41^{3}, 42^{8},$\\ $43^{2}, 44^{4}, 45^{9}, 46^{4}, 47^{1}, 48^{8}, 49^{2}, 50^{1}$
\end{minipage}  \\
\begin{minipage}{1 cm}
\centering
\smallskip
$41 $
\end{minipage}
&
\begin{minipage}{11 cm}
\centering
\smallskip
$ 33^{5}, 34^{2}, 35^{1}, 36^{13}, 37^{3}, 38^{2}, 39^{8}, 40^{6}, 41^{1}, 42^{13}, 43^{1}, 44^{6}, 45^{8}, 46^{2},$\\ $47^{3}, 48^{14}, 49^{2}, 50^{3}, 51^{6}, 52^{4}, 53^{1}, 54^{4}$
\end{minipage}  \\
\begin{minipage}{1 cm}
\centering
\smallskip
$43 $
\end{minipage}
&
\begin{minipage}{11 cm}
\centering
\smallskip
$ 32^{1}, 33^{3}, 36^{16}, 37^{1}, 38^{3}, 39^{6}, 40^{8}, 41^{1}, 42^{8}, 43^{5}, 44^{4}, 45^{13}, 46^{4}, 47^{1},$\\ $48^{16}, 49^{3}, 50^{4}, 51^{4}, 52^{6}, 53^{2}, 54^{9}, 55^{2}, 56^{2}, 57^{2}$
\end{minipage}  \\
\begin{minipage}{1 cm}
\centering
\smallskip
$47 $
\end{minipage}
&
\begin{minipage}{11 cm}
\centering
\smallskip
$ 36^{4}, 37^{1}, 38^{2}, 39^{6}, 40^{5}, 41^{3}, 42^{12}, 43^{1}, 44^{4}, 45^{9}, 46^{4}, 47^{2}, 48^{17}, 49^{2},$\\ $50^{4}, 51^{10}, 52^{4}, 53^{1}, 54^{12}, 55^{3}, 56^{6}, 57^{6}, 58^{2}, 59^{1}, 60^{8}, 61^{1}$
\end{minipage}  \\
\begin{minipage}{1 cm}
\centering
\smallskip
$49 $
\end{minipage}
&
\begin{minipage}{11 cm}
\centering
\smallskip
$ 36^{3}, 38^{2}, 39^{6}, 40^{2}, 42^{8}, 44^{6}, 45^{13}, 46^{6}, 47^{2}, 48^{16}, 49^{4}, 50^{2}, 51^{8}, 52^{8},$\\ $53^{2}, 54^{18}, 55^{5}, 56^{6}, 58^{4}, 59^{2}, 60^{12}, 61^{3}, 62^{2}, 63^{5}, 64^{1}$
\end{minipage}  \\
\begin{minipage}{1 cm}
\centering
\smallskip
$53 $
\end{minipage}
&
\begin{minipage}{11 cm}
\centering
\smallskip
$ 40^{2}, 41^{1}, 42^{7}, 43^{1}, 44^{4}, 45^{10}, 46^{2}, 47^{1}, 48^{20}, 49^{2}, 50^{5}, 51^{8}, 52^{6}, 53^{3},$\\ $54^{9}, 55^{3}, 56^{6}, 57^{8}, 58^{5}, 59^{2}, 60^{20}, 61^{1}, 62^{2}, 63^{10}, 64^{4}, 65^{2}, 66^{8}, 67^{1},$\\ $68^{2}$
\end{minipage}  \\
\begin{minipage}{1 cm}
\centering
\smallskip
$59 $
\end{minipage}
&
\begin{minipage}{11 cm}
\centering
\smallskip
$ 45^{1}, 46^{2}, 48^{12}, 49^{2}, 50^{4}, 51^{8}, 52^{4}, 53^{2}, 54^{14}, 55^{3}, 56^{8}, 57^{10}, 58^{2}, 59^{2},$\\ $60^{21}, 61^{2}, 62^{2}, 63^{10}, 64^{8}, 65^{3}, 66^{14}, 67^{2}, 68^{4}, 69^{8}, 70^{4}, 71^{2}, 72^{12}, 73^{1},$\\ $74^{2}, 75^{2}$
\end{minipage}  \\
\begin{minipage}{1 cm}
\centering
\smallskip
$61 $
\end{minipage}
&
\begin{minipage}{11 cm}
\centering
\smallskip
$ 47^{1}, 48^{7}, 49^{3}, 50^{3}, 51^{4}, 52^{8}, 53^{1}, 54^{17}, 55^{4}, 56^{6}, 57^{8}, 58^{4}, 59^{2}, 60^{16},$\\ $61^{5}, 62^{6}, 63^{14}, 64^{8}, 65^{2}, 66^{8}, 67^{4}, 68^{6}, 69^{8}, 70^{6}, 71^{1}, 72^{20}, 73^{2}, 74^{3},$\\ $75^{6}, 76^{4}, 77^{1}$
\end{minipage}  \\
\begin{minipage}{1 cm}
\centering
\smallskip
$64 $
\end{minipage}
&
\begin{minipage}{11 cm}
\centering
\smallskip
$ 49^{1}, 50^{3}, 52^{6}, 54^{22}, 56^{7}, 57^{4}, 58^{9}, 60^{24}, 62^{6}, 64^{12}, 65^{1}, 66^{24}, 68^{6}, 70^{12},$\\
$72^{27}, 73^{2}, 74^{7}, 76^{8}, 78^{12}, 80^{3}, 81^{2}$
\end{minipage}  \\
\begin{minipage}{1 cm}
\centering
\smallskip
$67 $
\end{minipage}
&
\begin{minipage}{11 cm}
\centering
\smallskip
$ 52^{2}, 53^{1}, 54^{9}, 55^{3}, 56^{6}, 57^{6}, 58^{4}, 59^{2}, 60^{16}, 61^{4}, 62^{2}, 63^{14}, 64^{10}, 65^{4},$\\ $66^{16}, 67^{2}, 68^{4}, 69^{4}, 70^{8}, 71^{4}, 72^{26}, 73^{5}, 74^{2}, 75^{8}, 76^{8}, 77^{2}, 78^{8}, 79^{3},$\\ $80^{6}, 81^{9}, 82^{3}, 83^{1}, 84^{4}$
\end{minipage}  \\
\begin{minipage}{1 cm}
\centering
\smallskip
$71 $
\end{minipage}
&
\begin{minipage}{11 cm}
\centering
\smallskip
$ 56^{2}, 57^{6}, 58^{2}, 59^{2}, 60^{16}, 61^{1}, 62^{4}, 63^{8}, 64^{8}, 65^{2}, 66^{16}, 67^{4}, 68^{4}, 69^{10},$\\ $70^{4}, 71^{3}, 72^{25}, 73^{3}, 74^{4}, 75^{10}, 76^{4}, 77^{4}, 78^{16}, 79^{2}, 80^{8}, 81^{8}, 82^{4}, 83^{1},$\\ $84^{16}, 85^{2}, 86^{2}, 87^{6}, 88^{2}$
\end{minipage}  \\
\begin{minipage}{1 cm}
\centering
\smallskip
$73 $
\end{minipage}
&
\begin{minipage}{11 cm}
\centering
\smallskip
$ 57^{2}, 58^{3}, 59^{1}, 60^{11}, 61^{2}, 62^{2}, 63^{13}, 64^{8}, 65^{3}, 66^{8}, 67^{5}, 68^{8}, 69^{4}, 70^{8},$\\ $71^{3}, 72^{27}, 73^{4}, 74^{4}, 75^{8}, 76^{9}, 77^{3}, 78^{16}, 79^{2}, 80^{8}, 81^{14}, 82^{4}, 83^{3}, 84^{16},$\\ $85^{5}, 86^{2}, 87^{4}, 88^{6}, 89^{1}, 90^{7}, 91^{1}$
\end{minipage}  \\
\begin{minipage}{1 cm}
\centering
\smallskip
$79 $
\end{minipage}
&
\begin{minipage}{11 cm}
\centering
\smallskip
$ 63^{5}, 64^{4}, 65^{2}, 66^{8}, 67^{3}, 68^{4}, 69^{8}, 70^{8}, 71^{2}, 72^{26}, 73^{2}, 74^{4}, 75^{8}, 76^{10},$\\ $77^{3}, 78^{8}, 79^{6}, 80^{10}, 81^{18}, 82^{4}, 83^{3}, 84^{20}, 85^{4}, 86^{4}, 87^{4}, 88^{10}, 89^{2},$\\ $90^{22}, 91^{4}, 92^{4}, 93^{6}, 94^{4}, 95^{2}, 96^{8}, 97^{2}$
\end{minipage}  \\
\begin{minipage}{1 cm}
\centering
\smallskip
$81 $
\end{minipage}
&
\begin{minipage}{11 cm}
\centering
\smallskip
$ 64^{1}, 65^{2}, 66^{8}, 68^{7}, 69^{8}, 71^{4}, 72^{24}, 73^{1}, 74^{8}, 75^{10}, 77^{8}, 78^{16}, 80^{14}, 81^{8},$\\ $83^{4}, 84^{28}, 86^{8}, 87^{16}, 89^{5}, 90^{16}, 91^{1}, 92^{12}, 93^{8}, 95^{4}, 96^{14}, 98^{4}, 99^{4}, 100^{1}$
\end{minipage}  \\
\hline
\end{tabular}
\end{minipage}\bigskip


\section{Detailed description of planar cubic curves, $q\leq81$}

In this section we give a detailed description of absolutely irreducible (a. i.)
singular cubic curves and of non-singular cubic curves for each $q\leq81$.
Complete and incomplete cubic curves are described in different tables. For
singular cubics the type, the equation and the number of $\mathbb{F}_{q}%
$-rational points are given. For incomplete curves also the number of residual
points is given. The symbol $\mathcal{N}_{i}^{j}$ denotes an absolutely
irreducible singular cubic with $i$ $\mathbb{F}_{q}$-rational inflexions and
$j$ distinct $\mathbb{F}_{q}$-rational tangents at the double point.


For non-singular cubics the type, the equation, the number of $\mathbb{F}_{q}%
$-rational points and the $j$-invariant are given. For incomplete curves also
the number of residual points is given.

The symbols $\mathcal{G}$, $\mathcal{H}$, $\mathcal{E}$, $\mathcal{S}$ denote
general, harmonic, equianharmonic, superharmonic non-singular cubic curves, respectively.

Also, the lower suffix on $\mathcal{G}$, $\mathcal{H}$, $\mathcal{E}$,
$\mathcal{S}$ is the number of rational inflexions; the upper suffix on
$\mathcal{G}$, $\mathcal{H}$, $\mathcal{E}$ is the number of rational
inflexional triangles; the upper suffix on $\mathcal{S}$ is the number of
rational points in which the harmonic polar meets the curve; $\overline
{\mathcal{E}}_{3}^{r}$ denotes equianharmonic non-singular cubic curves with
three concurrent inflexional tangents; $\overline{\mathcal{E}}_{0}^{4}$
denotes equianharmonic non-singular cubic curves with canonical equation of
the form $X^{3} + \xi Y^{3} + \xi^{2} T^{3}$, with $\xi$ a primitive element
of $\mathbb{F}_{q}$; see \cite[Cap. 11]{hirsh} for the details.



$\; $ \phantom  \newline\newline\noindent\begin{minipage}{\textwidth}
\centering
\textbf{Table 6} $ q = 2$, a. i. singular cubic curves giving incomplete $(n,3)$-arcs  \medskip \\
\medskip
\begin{tabular}
{|c|c|c|c|}
\hline
type & equation &$n$ &
\begin{minipage}{1.3 cm}
\smallskip
\centering
$\#$ of \\
residual\\
points
\end{minipage}\\
\hline
\begin{minipage}{0.7 cm}
\centering
\smallskip
$ \mathcal{N}_3^0$  \end{minipage}
&   \begin{minipage}{6 cm}
\centering
\smallskip
$ YT^2 + Y^2T + XT^2 + XYT + XY^2 $ \end{minipage}
&  \begin{minipage}{0.7 cm}
\centering
\smallskip
$4$ \end{minipage}
&  \begin{minipage}{1.3 cm}
\centering
\smallskip
$3$ \end{minipage}
\\
\begin{minipage}{0.7 cm}
\centering
\smallskip
$ \mathcal{N}_1^2$  \end{minipage}
&   \begin{minipage}{6 cm}
\centering
\smallskip
$ T^3 + Y^3 + XYT $ \end{minipage}
&  \begin{minipage}{0.7 cm}
\centering
\smallskip
$2$ \end{minipage}
&  \begin{minipage}{1.3 cm}
\centering
\smallskip
$5$ \end{minipage}
\\
\begin{minipage}{0.7 cm}
\centering
\smallskip
$ \mathcal{N}_1^1$  \end{minipage}
&   \begin{minipage}{6 cm}
\centering
\smallskip
$ T^3 + XY^2 $ \end{minipage}
&  \begin{minipage}{0.7 cm}
\centering
\smallskip
$3$ \end{minipage}
&  \begin{minipage}{1.3 cm}
\centering
\smallskip
$4$ \end{minipage}
\\
\begin{minipage}{0.7 cm}
\centering
\smallskip
$ \mathcal{N}_0^0$  \end{minipage}
&   \begin{minipage}{6 cm}
\centering
\smallskip
$ T^3 + Y^2T + XT^2 + Y^3 + XYT + XY^2 $ \end{minipage}
&  \begin{minipage}{0.7 cm}
\centering
\smallskip
$4$ \end{minipage}
&  \begin{minipage}{1.3 cm}
\centering
\smallskip
$3$ \end{minipage}
\\
\hline
\end{tabular}
\end{minipage}\newline\phantom  \newline\newline


\noindent\begin{minipage}{\textwidth}
\centering
\textbf{Table 7}  $ q = 2$, non-singular cubic curves giving incomplete $(n,3)$-arcs  \medskip \\
\medskip
\begin{tabular}
{|c|c|c|c|c|}
\hline
type & equation & $n$ &
\begin{minipage}{1.3 cm}
\smallskip
\centering
$\#$ of \\
residual\\
points
\end{minipage}   &
\begin{minipage}{1 cm}
\smallskip
\centering
$j-$ \\
inv. \\
\end{minipage}  \\
\hline
\begin{minipage}{0.7 cm}
\centering
\smallskip
$ \overline{\mathcal{E}}_3^2$  \end{minipage}
&   \begin{minipage}{5 cm}
\centering
\smallskip
$ T^3 + XY^2 + X^2Y $ \end{minipage}
&  \begin{minipage}{0.7 cm}
\centering
\smallskip
$3$ \end{minipage}
&  \begin{minipage}{1.3 cm}
\centering
\smallskip
$4$ \end{minipage}
&  \begin{minipage}{1 cm}
\centering
\smallskip
$0$ \end{minipage}
\\
\begin{minipage}{0.7 cm}
\centering
\smallskip
$ \mathcal{G}_1^0$  \end{minipage}
&   \begin{minipage}{5 cm}
\centering
\smallskip
$ YT^2 + XYT + XY^2 + X^3 $ \end{minipage}
&  \begin{minipage}{0.7 cm}
\centering
\smallskip
$4$ \end{minipage}
&  \begin{minipage}{1.3 cm}
\centering
\smallskip
$3$ \end{minipage}
&  \begin{minipage}{1 cm}
\centering
\smallskip
$1$ \end{minipage}
\\
\begin{minipage}{0.7 cm}
\centering
\smallskip
$ \mathcal{G}_1^0$  \end{minipage}
&   \begin{minipage}{5 cm}
\centering
\smallskip
$ YT^2 + XYT + XY^2 + X^2Y + X^3 $ \end{minipage}
&  \begin{minipage}{0.7 cm}
\centering
\smallskip
$2$ \end{minipage}
&  \begin{minipage}{1.3 cm}
\centering
\smallskip
$5$ \end{minipage}
&  \begin{minipage}{1 cm}
\centering
\smallskip
$1$ \end{minipage}
\\
\begin{minipage}{0.7 cm}
\centering
\smallskip
$ \mathcal{E}_1^0$  \end{minipage}
&   \begin{minipage}{5 cm}
\centering
\smallskip
$ YT^2 + Y^2T + XY^2 + X^3 $ \end{minipage}
&  \begin{minipage}{0.7 cm}
\centering
\smallskip
$5$ \end{minipage}
&  \begin{minipage}{1.3 cm}
\centering
\smallskip
$2$ \end{minipage}
&  \begin{minipage}{1 cm}
\centering
\smallskip
$0$ \end{minipage}
\\
\begin{minipage}{0.7 cm}
\centering
\smallskip
$ \mathcal{E}_1^0$  \end{minipage}
&   \begin{minipage}{5 cm}
\centering
\smallskip
$ YT^2 + Y^2T + Y^3 + XY^2 + X^3 $ \end{minipage}
&  \begin{minipage}{0.7 cm}
\centering
\smallskip
$1$ \end{minipage}
&  \begin{minipage}{1.3 cm}
\centering
\smallskip
$6$ \end{minipage}
&  \begin{minipage}{1 cm}
\centering
\smallskip
$0$ \end{minipage}
\\
\begin{minipage}{0.7 cm}
\centering
\smallskip
$ \mathcal{E}_0^2$  \end{minipage}
&   \begin{minipage}{5 cm}
\centering
\smallskip
$ T^3 + Y^3 + XY^2 + X^3 $ \end{minipage}
&  \begin{minipage}{0.7 cm}
\centering
\smallskip
$3$ \end{minipage}
&  \begin{minipage}{1.3 cm}
\centering
\smallskip
$4$ \end{minipage}
&  \begin{minipage}{1 cm}
\centering
\smallskip
$0$ \end{minipage}
\\
\hline
\end{tabular}
\end{minipage}\newline


$\; $ \phantom  \newline\newline\noindent\begin{minipage}{\textwidth}
\centering
\textbf{Table 8} $ q = 3$, a. i. singular cubic curves giving incomplete $(n,3)$-arcs  \medskip \\
\medskip
\begin{tabular}
{|c|c|c|c|}
\hline
type & equation &$n$ &
\begin{minipage}{1.3 cm}
\smallskip
\centering
$\#$ of \\
residual\\
points
\end{minipage}\\
\hline
\begin{minipage}{0.7 cm}
\centering
\smallskip
$ \mathcal{N}_q^1$  \end{minipage}
&   \begin{minipage}{3 cm}
\centering
\smallskip
$ T^3 + XY^2 $ \end{minipage}
&  \begin{minipage}{0.7 cm}
\centering
\smallskip
$4$ \end{minipage}
&  \begin{minipage}{1.3 cm}
\centering
\smallskip
$8$ \end{minipage}
\\
\begin{minipage}{0.7 cm}
\centering
\smallskip
$ \mathcal{N}_1^2$  \end{minipage}
&   \begin{minipage}{3 cm}
\centering
\smallskip
$ T^3 + Y^3 + XYT $ \end{minipage}
&  \begin{minipage}{0.7 cm}
\centering
\smallskip
$3$ \end{minipage}
&  \begin{minipage}{1.3 cm}
\centering
\smallskip
$10$ \end{minipage}
\\
\begin{minipage}{0.7 cm}
\centering
\smallskip
$ \mathcal{N}_1^0$  \end{minipage}
&   \begin{minipage}{3 cm}
\centering
\smallskip
$ XT^2 + Y^3 + XY^2 $ \end{minipage}
&  \begin{minipage}{0.7 cm}
\centering
\smallskip
$5$ \end{minipage}
&  \begin{minipage}{1.3 cm}
\centering
\smallskip
$7$ \end{minipage}
\\
\begin{minipage}{0.7 cm}
\centering
\smallskip
$ \mathcal{N}_0^1$  \end{minipage}
&   \begin{minipage}{3 cm}
\centering
\smallskip
$ T^3 + YT^2 + XY^2 $ \end{minipage}
&  \begin{minipage}{0.7 cm}
\centering
\smallskip
$4$ \end{minipage}
&  \begin{minipage}{1.3 cm}
\centering
\smallskip
$9$ \end{minipage}
\\
\hline
\end{tabular}
\end{minipage}\newline\phantom  \newline\newline


\noindent\begin{minipage}{\textwidth}
\centering
\textbf{Table 9}  $ q = 3$, non-singular cubic curves giving incomplete $(n,3)$-arcs  \medskip \\
\medskip
\begin{tabular}
{|c|c|c|c|c|}
\hline
type & equation &$n$ &
\begin{minipage}{1.3 cm}
\smallskip
\centering
$\#$ of \\
residual\\
points
\end{minipage}   &
\begin{minipage}{1 cm}
\smallskip
\centering
$j-$ \\
inv. \\
\end{minipage}  \\
\hline
\begin{minipage}{0.7 cm}
\centering
\smallskip
$ \mathcal{G}_3$  \end{minipage}
&   \begin{minipage}{7 cm}
\centering
\smallskip
$ T^3 + Y^3 + XYT + X^3 $ \end{minipage}
&  \begin{minipage}{0.7 cm}
\centering
\smallskip
$6$ \end{minipage}
&  \begin{minipage}{1.3 cm}
\centering
\smallskip
$3$ \end{minipage}
&  \begin{minipage}{1 cm}
\centering
\smallskip
$2$ \end{minipage}
\\
\begin{minipage}{0.7 cm}
\centering
\smallskip
$ \mathcal{G}_3$  \end{minipage}
&   \begin{minipage}{7 cm}
\centering
\smallskip
$ 2T^3 + 2Y^3 + XYT + 2X^3 $ \end{minipage}
&  \begin{minipage}{0.7 cm}
\centering
\smallskip
$3$ \end{minipage}
&  \begin{minipage}{1.3 cm}
\centering
\smallskip
$9$ \end{minipage}
&  \begin{minipage}{1 cm}
\centering
\smallskip
$1$ \end{minipage}
\\
\begin{minipage}{0.7 cm}
\centering
\smallskip
$ \mathcal{S}_1^3$  \end{minipage}
&   \begin{minipage}{7 cm}
\centering
\smallskip
$ YT^2 + 2XY^2 + X^3 $ \end{minipage}
&  \begin{minipage}{0.7 cm}
\centering
\smallskip
$4$ \end{minipage}
&  \begin{minipage}{1.3 cm}
\centering
\smallskip
$8$ \end{minipage}
&  \begin{minipage}{1 cm}
\centering
\smallskip
$0$ \end{minipage}
\\
\begin{minipage}{0.7 cm}
\centering
\smallskip
$ \mathcal{S}_1^1$  \end{minipage}
&   \begin{minipage}{7 cm}
\centering
\smallskip
$ YT^2 + XY^2 + X^3 $ \end{minipage}
&  \begin{minipage}{0.7 cm}
\centering
\smallskip
$4$ \end{minipage}
&  \begin{minipage}{1.3 cm}
\centering
\smallskip
$8$ \end{minipage}
&  \begin{minipage}{1 cm}
\centering
\smallskip
$0$ \end{minipage}
\\
\begin{minipage}{0.7 cm}
\centering
\smallskip
$ \mathcal{S}_1^0$  \end{minipage}
&   \begin{minipage}{7 cm}
\centering
\smallskip
$ YT^2 + 2Y^3 + 2XY^2 + X^3 $ \end{minipage}
&  \begin{minipage}{0.7 cm}
\centering
\smallskip
$7$ \end{minipage}
&  \begin{minipage}{1.3 cm}
\centering
\smallskip
$2$ \end{minipage}
&  \begin{minipage}{1 cm}
\centering
\smallskip
$0$ \end{minipage}
\\
\begin{minipage}{0.7 cm}
\centering
\smallskip
$ \mathcal{S}_1^0$  \end{minipage}
&   \begin{minipage}{7 cm}
\centering
\smallskip
$ YT^2 + Y^3 + 2XY^2 + X^3 $ \end{minipage}
&  \begin{minipage}{0.7 cm}
\centering
\smallskip
$1$ \end{minipage}
&  \begin{minipage}{1.3 cm}
\centering
\smallskip
$12$ \end{minipage}
&  \begin{minipage}{1 cm}
\centering
\smallskip
$0$ \end{minipage}
\\
\begin{minipage}{0.7 cm}
\centering
\smallskip
$ \mathcal{G}_1$  \end{minipage}
&   \begin{minipage}{7 cm}
\centering
\smallskip
$ YT^2 + 2Y^3 + X^2Y + X^3 $ \end{minipage}
&  \begin{minipage}{0.7 cm}
\centering
\smallskip
$5$ \end{minipage}
&  \begin{minipage}{1.3 cm}
\centering
\smallskip
$6$ \end{minipage}
&  \begin{minipage}{1 cm}
\centering
\smallskip
$1$ \end{minipage}
\\
\begin{minipage}{0.7 cm}
\centering
\smallskip
$ \mathcal{G}_1$  \end{minipage}
&   \begin{minipage}{7 cm}
\centering
\smallskip
$ YT^2 + Y^3 + X^2Y + X^3 $ \end{minipage}
&  \begin{minipage}{0.7 cm}
\centering
\smallskip
$2$ \end{minipage}
&  \begin{minipage}{1.3 cm}
\centering
\smallskip
$11$ \end{minipage}
&  \begin{minipage}{1 cm}
\centering
\smallskip
$2$ \end{minipage}
\\
\begin{minipage}{0.7 cm}
\centering
\smallskip
$ \mathcal{G}_0$  \end{minipage}
&   \begin{minipage}{7 cm}
\centering
\smallskip
$ 2YT^2 + Y^3 + 2XY^2 + 2X^2T + X^2Y + X^3 $ \end{minipage}
&  \begin{minipage}{0.7 cm}
\centering
\smallskip
$6$ \end{minipage}
&  \begin{minipage}{1.3 cm}
\centering
\smallskip
$4$ \end{minipage}
&  \begin{minipage}{1 cm}
\centering
\smallskip
$2$ \end{minipage}
\\
\begin{minipage}{0.7 cm}
\centering
\smallskip
$ \mathcal{G}_0$  \end{minipage}
&   \begin{minipage}{7 cm}
\centering
\smallskip
$ 2T^3 + 2YT^2 + Y^3 + 2XY^2 + 2X^2T + X^2Y + X^3 $ \end{minipage}
&  \begin{minipage}{0.7 cm}
\centering
\smallskip
$3$ \end{minipage}
&  \begin{minipage}{1.3 cm}
\centering
\smallskip
$10$ \end{minipage}
&  \begin{minipage}{1 cm}
\centering
\smallskip
$1$ \end{minipage}
\\
\hline
\end{tabular}
\end{minipage}\newline\phantom    \newline\phantom   \newline The sets of $\mathbb{F}_q$-rational points of cubics $YT^{2} + XY^{2} + X^{3} $ and $YT^{2} + 2XY^{2} +
X^{3} $ are equivalent.\newline


$\; $ \phantom  \newline\newline\noindent\begin{minipage}{\textwidth}
\centering
\textbf{Table 10} $ q = 4$, a. i. singular cubic curves giving incomplete $(n,3)$-arcs  \medskip \\
\medskip
\begin{tabular}
{|c|c|c|c|}
\hline
type & equation &$n$ &
\begin{minipage}{1.3 cm}
\smallskip
\centering
$\#$ of \\
residual\\
points
\end{minipage}\\
\hline
\begin{minipage}{0.7 cm}
\centering
\smallskip
$ \mathcal{N}_3^2$  \end{minipage}
&   \begin{minipage}{8 cm}
\centering
\smallskip
$ T^3 + Y^3 + XYT $ \end{minipage}
&  \begin{minipage}{0.7 cm}
\centering
\smallskip
$4$ \end{minipage}
&  \begin{minipage}{1.3 cm}
\centering
\smallskip
$15$ \end{minipage}
\\
\begin{minipage}{0.7 cm}
\centering
\smallskip
$ \mathcal{N}_1^1$  \end{minipage}
&   \begin{minipage}{6 cm}
\centering
\smallskip
$ T^3 + XY^2 $ \end{minipage}
&  \begin{minipage}{0.7 cm}
\centering
\smallskip
$5$ \end{minipage}
&  \begin{minipage}{1.3 cm}
\centering
\smallskip
$14$ \end{minipage}
\\
\begin{minipage}{0.7 cm}
\centering
\smallskip
$ \mathcal{N}_1^0$  \end{minipage}
&   \begin{minipage}{8 cm}
\centering
\smallskip
$ \xi^{2}YT^2 + \xi Y^2T + \xi XT^2 + \xi^{2}Y^3 + XYT + XY^2 $ \end{minipage}
&  \begin{minipage}{0.7 cm}
\centering
\smallskip
$6$ \end{minipage}
&  \begin{minipage}{1.3 cm}
\centering
\smallskip
$11$ \end{minipage}
\\
\begin{minipage}{0.7 cm}
\centering
\smallskip
$ \mathcal{N}_0^2$  \end{minipage}
&   \begin{minipage}{6 cm}
\centering
\smallskip
$ \xi T^3 + Y^3 + XYT $ \end{minipage}
&  \begin{minipage}{0.7 cm}
\centering
\smallskip
$4$ \end{minipage}
&  \begin{minipage}{1.3 cm}
\centering
\smallskip
$17$ \end{minipage}
\\
\hline
\end{tabular}
\end{minipage}\newline\phantom  \newline\newline


\noindent\begin{minipage}{\textwidth}
\centering
\textbf{Table 11}  $ q = 4$, non-singular cubic curves giving incomplete $(n,3)$-arcs  \medskip \\
\medskip
\begin{tabular}
{|c|c|c|c|c|}
\hline
type & equation &$n$ &
\begin{minipage}{1.3 cm}
\smallskip
\centering
$\#$ of \\
residual\\
points
\end{minipage}   &
\begin{minipage}{1 cm}
\smallskip
\centering
$j-$ \\
inv. \\
\end{minipage}  \\
\hline
\begin{minipage}{0.7 cm}
\centering
\smallskip
$ \mathcal{G}_3^1$  \end{minipage}
&   \begin{minipage}{8 cm}
\centering
\smallskip
$ \xi^{2}T^3 + \xi^{2}YT^2 + \xi^{2}Y^2T + \xi^{2}XT^2 + \xi^{2}Y^3 + XYT + \xi^{2}XY^2 + \xi^{2}X^2T + \xi^{2}X^2Y + \xi^{2}X^3 $ \end{minipage}
&  \begin{minipage}{0.7 cm}
\centering
\smallskip
$6$ \end{minipage}
&  \begin{minipage}{1.3 cm}
\centering
\smallskip
$7$ \end{minipage}
&  \begin{minipage}{1 cm}
\centering
\smallskip
$\xi^{2}$ \end{minipage}
\\
\begin{minipage}{0.7 cm}
\centering
\smallskip
$ \mathcal{G}_3^1$  \end{minipage}
&   \begin{minipage}{8 cm}
\centering
\smallskip
$ \xi T^3 + \xi YT^2 + \xi Y^2T + \xi XT^2 + \xi Y^3 + XYT + \xi XY^2 + \xi X^2T + \xi X^2Y + \xi X^3 $ \end{minipage}
&  \begin{minipage}{0.7 cm}
\centering
\smallskip
$6$ \end{minipage}
&  \begin{minipage}{1.3 cm}
\centering
\smallskip
$7$ \end{minipage}
&  \begin{minipage}{1 cm}
\centering
\smallskip
$\xi$ \end{minipage}
\\
\begin{minipage}{0.7 cm}
\centering
\smallskip
$ \overline{\mathcal{E}}_3^1$  \end{minipage}
&   \begin{minipage}{8 cm}
\centering
\smallskip
$ \xi^{2}T^3 + XY^2 + X^2Y $ \end{minipage}
&  \begin{minipage}{0.7 cm}
\centering
\smallskip
$3$ \end{minipage}
&  \begin{minipage}{1.3 cm}
\centering
\smallskip
$16$ \end{minipage}
&  \begin{minipage}{1 cm}
\centering
\smallskip
$0$ \end{minipage}
\\
\begin{minipage}{0.7 cm}
\centering
\smallskip
$ \overline{\mathcal{E}}_3^1$  \end{minipage}
&   \begin{minipage}{8 cm}
\centering
\smallskip
$ \xi T^3 + XY^2 + X^2Y $ \end{minipage}
&  \begin{minipage}{0.7 cm}
\centering
\smallskip
$3$ \end{minipage}
&  \begin{minipage}{1.3 cm}
\centering
\smallskip
$16$ \end{minipage}
&  \begin{minipage}{1 cm}
\centering
\smallskip
$0$ \end{minipage}
\\
\begin{minipage}{0.7 cm}
\centering
\smallskip
$ \mathcal{E}_1^4$  \end{minipage}
&   \begin{minipage}{8 cm}
\centering
\smallskip
$ YT^2 + Y^2T + \xi Y^3 + X^3 $ \end{minipage}
&  \begin{minipage}{0.7 cm}
\centering
\smallskip
$1$ \end{minipage}
&  \begin{minipage}{1.3 cm}
\centering
\smallskip
$20$ \end{minipage}
&  \begin{minipage}{1 cm}
\centering
\smallskip
$0$ \end{minipage}
\\
\begin{minipage}{0.7 cm}
\centering
\smallskip
$ \mathcal{G}_1^1$  \end{minipage}
&   \begin{minipage}{8 cm}
\centering
\smallskip
$ YT^2 + XYT + \xi^{2}XY^2 + X^3 $ \end{minipage}
&  \begin{minipage}{0.7 cm}
\centering
\smallskip
$4$ \end{minipage}
&  \begin{minipage}{1.3 cm}
\centering
\smallskip
$15$ \end{minipage}
&  \begin{minipage}{1 cm}
\centering
\smallskip
$\xi^{2}$ \end{minipage}
\\
\begin{minipage}{0.7 cm}
\centering
\smallskip
$ \mathcal{G}_1^1$  \end{minipage}
&   \begin{minipage}{8 cm}
\centering
\smallskip
$ YT^2 + XYT + \xi XY^2 + X^3 $ \end{minipage}
&  \begin{minipage}{0.7 cm}
\centering
\smallskip
$4$ \end{minipage}
&  \begin{minipage}{1.3 cm}
\centering
\smallskip
$15$ \end{minipage}
&  \begin{minipage}{1 cm}
\centering
\smallskip
$\xi$ \end{minipage}
\\
\begin{minipage}{0.7 cm}
\centering
\smallskip
$ \mathcal{E}_1^1$  \end{minipage}
&   \begin{minipage}{8 cm}
\centering
\smallskip
$ YT^2 + Y^2T + XY^2 + \xi^{2}X^3 $ \end{minipage}
&  \begin{minipage}{0.7 cm}
\centering
\smallskip
$7$ \end{minipage}
&  \begin{minipage}{1.3 cm}
\centering
\smallskip
$5$ \end{minipage}
&  \begin{minipage}{1 cm}
\centering
\smallskip
$0$ \end{minipage}
\\
\begin{minipage}{0.7 cm}
\centering
\smallskip
$ \mathcal{E}_1^1$  \end{minipage}
&   \begin{minipage}{8 cm}
\centering
\smallskip
$ YT^2 + Y^2T + XY^2 + \xi X^3 $ \end{minipage}
&  \begin{minipage}{0.7 cm}
\centering
\smallskip
$7$ \end{minipage}
&  \begin{minipage}{1.3 cm}
\centering
\smallskip
$5$ \end{minipage}
&  \begin{minipage}{1 cm}
\centering
\smallskip
$0$ \end{minipage}
\\
\begin{minipage}{0.7 cm}
\centering
\smallskip
$ \mathcal{G}_1^0$  \end{minipage}
&   \begin{minipage}{8 cm}
\centering
\smallskip
$ YT^2 + XYT + XY^2 + X^3 $ \end{minipage}
&  \begin{minipage}{0.7 cm}
\centering
\smallskip
$8$ \end{minipage}
&  \begin{minipage}{1.3 cm}
\centering
\smallskip
$2$ \end{minipage}
&  \begin{minipage}{1 cm}
\centering
\smallskip
$1$ \end{minipage}
\\
\begin{minipage}{0.7 cm}
\centering
\smallskip
$ \mathcal{G}_1^0$  \end{minipage}
&   \begin{minipage}{8 cm}
\centering
\smallskip
$ YT^2 + XYT + XY^2 + \xi X^2Y + X^3 $ \end{minipage}
&  \begin{minipage}{0.7 cm}
\centering
\smallskip
$2$ \end{minipage}
&  \begin{minipage}{1.3 cm}
\centering
\smallskip
$19$ \end{minipage}
&  \begin{minipage}{1 cm}
\centering
\smallskip
$1$ \end{minipage}
\\
\begin{minipage}{0.7 cm}
\centering
\smallskip
$ \mathcal{E}_1^0$  \end{minipage}
&   \begin{minipage}{8 cm}
\centering
\smallskip
$ YT^2 + Y^2T + \xi Y^3 + XY^2 + X^3 $ \end{minipage}
&  \begin{minipage}{0.7 cm}
\centering
\smallskip
$5$ \end{minipage}
&  \begin{minipage}{1.3 cm}
\centering
\smallskip
$12$ \end{minipage}
&  \begin{minipage}{1 cm}
\centering
\smallskip
$0$ \end{minipage}
\\
\begin{minipage}{0.7 cm}
\centering
\smallskip
$ \mathcal{G}_0^1$  \end{minipage}
&   \begin{minipage}{8 cm}
\centering
\smallskip
$ \xi^{2}T^3 + \xi YT^2 + \xi Y^3 + \xi XYT + XY^2 + X^2T + X^3 $ \end{minipage}
&  \begin{minipage}{0.7 cm}
\centering
\smallskip
$6$ \end{minipage}
&  \begin{minipage}{1.3 cm}
\centering
\smallskip
$9$ \end{minipage}
&  \begin{minipage}{1 cm}
\centering
\smallskip
$\xi^{2}$ \end{minipage}
\\
\begin{minipage}{0.7 cm}
\centering
\smallskip
$ \mathcal{G}_0^1$  \end{minipage}
&   \begin{minipage}{8 cm}
\centering
\smallskip
$ \xi T^3 + \xi^{2}YT^2 + \xi^{2}Y^3 + \xi^{2}XYT + XY^2 + X^2T + X^3 $ \end{minipage}
&  \begin{minipage}{0.7 cm}
\centering
\smallskip
$6$ \end{minipage}
&  \begin{minipage}{1.3 cm}
\centering
\smallskip
$9$ \end{minipage}
&  \begin{minipage}{1 cm}
\centering
\smallskip
$\xi$ \end{minipage}
\\
\begin{minipage}{0.7 cm}
\centering
\smallskip
$ \mathcal{E}_0^1$  \end{minipage}
&   \begin{minipage}{8 cm}
\centering
\smallskip
$ \xi^{2}YT^2 + XY^2 + X^2T $ \end{minipage}
&  \begin{minipage}{0.7 cm}
\centering
\smallskip
$3$ \end{minipage}
&  \begin{minipage}{1.3 cm}
\centering
\smallskip
$18$ \end{minipage}
&  \begin{minipage}{1 cm}
\centering
\smallskip
$0$ \end{minipage}
\\
\begin{minipage}{0.7 cm}
\centering
\smallskip
$ \mathcal{E}_0^1$  \end{minipage}
&   \begin{minipage}{8 cm}
\centering
\smallskip
$ \xi YT^2 + XY^2 + X^2T $ \end{minipage}
&  \begin{minipage}{0.7 cm}
\centering
\smallskip
$3$ \end{minipage}
&  \begin{minipage}{1.3 cm}
\centering
\smallskip
$18$ \end{minipage}
&  \begin{minipage}{1 cm}
\centering
\smallskip
$0$ \end{minipage}
\\
\hline
\end{tabular}
\end{minipage}\newline\phantom    \newline\phantom   \newline The sets of $\mathbb{F}_q$-rational points of cubics $\xi YT^{2} + XY^{2} + X^{2}T $ and $\xi^{2}YT^{2} +
XY^{2} + X^{2}T $ are equivalent.\newline The sets of $\mathbb{F}_q$-rational points of
cubics $\xi T^{3} + \xi YT^{2} + \xi Y^{2}T + \xi XT^{2} + \xi Y^{3} + XYT +
\xi XY^{2} + \xi X^{2}T + \xi X^{2}Y + \xi X^{3} $ and $\xi^{2}T^{3} + \xi
^{2}YT^{2} + \xi^{2}Y^{2}T + \xi^{2}XT^{2} + \xi^{2}Y^{3} + XYT + \xi
^{2}XY^{2} + \xi^{2}X^{2}T + \xi^{2}X^{2}Y + \xi^{2}X^{3} $ are
equivalent.\newline The sets of $\mathbb{F}_q$-rational points of cubics $\xi T^{3} + XY^{2}
+ X^{2}Y $ and $\xi^{2}T^{3} + XY^{2} + X^{2}Y $ are equivalent.\newline The
sets of rational points of cubics $\xi^{2}T^{3} + \xi YT^{2} + \xi Y^{3} + \xi
XYT + XY^{2} + X^{2}T + X^{3} $ and $\xi T^{3} + \xi^{2}YT^{2} + \xi^{2}Y^{3}
+ \xi^{2}XYT + XY^{2} + X^{2}T + X^{3} $ are equivalent.\newline The sets of $\mathbb{F}_q$-rational points of cubics $YT^{2} + XYT + \xi XY^{2} + X^{3} $ and $YT^{2} +
XYT + \xi^{2}XY^{2} + X^{3} $ are equivalent.\newline The sets of $\mathbb{F}_q$-rational
points of cubics $YT^{2} + Y^{2}T + XY^{2} + \xi X^{3} $ and $YT^{2} + Y^{2}T
+ XY^{2} + \xi^{2}X^{3} $ are equivalent.\newline\phantom  \newline\newline


\noindent\begin{minipage}{\textwidth}
\centering
\textbf{Table 12}  $ q = 4$, non-singular cubic curves giving complete $(n,3)$-arcs  \medskip \\
\medskip
\begin{tabular}
{|c|c|c|c|}
\hline
type & equation &$n$ &
\begin{minipage}{1 cm}
\smallskip
\centering
$j-$ \\
inv. \\
\smallskip
\end{minipage}  \\
\hline
\begin{minipage}{0.7 cm}
\centering
\smallskip
$ \mathcal{E}_9^4$  \end{minipage}
&   \begin{minipage}{3 cm}
\centering
\smallskip
$ T^3 + Y^3 + X^3 $ \end{minipage}
&  \begin{minipage}{0.7 cm}
\centering
\smallskip
$9$ \end{minipage}
&  \begin{minipage}{1 cm}
\centering
\smallskip
$0$ \end{minipage}
\\
\begin{minipage}{0.7 cm}
\centering
\smallskip
$ \overline{\mathcal{E}}_0^4$  \end{minipage}
&   \begin{minipage}{3 cm}
\centering
\smallskip
$ \xi^{2}T^3 + \xi Y^3 + X^3 $ \end{minipage}
&  \begin{minipage}{0.7 cm}
\centering
\smallskip
$9$ \end{minipage}
&  \begin{minipage}{1 cm}
\centering
\smallskip
$0$ \end{minipage}
\\
\hline
\end{tabular}
\end{minipage}\newline


$\; $ \phantom  \newline\newline\noindent\begin{minipage}{\textwidth}
\centering
\textbf{Table 13} $ q = 5$, a. i. singular cubic curves giving incomplete $(n,3)$-arcs  \medskip \\
\medskip
\begin{tabular}
{|c|c|c|c|}
\hline
type & equation &$n$ &
\begin{minipage}{1.3 cm}
\smallskip
\centering
$\#$ of \\
residual\\
points
\end{minipage}\\
\hline
\begin{minipage}{0.7 cm}
\centering
\smallskip
$ \mathcal{N}_3^0$  \end{minipage}
&   \begin{minipage}{6 cm}
\centering
\smallskip
$ YT^2 + 3XT^2 + Y^3 + XY^2 $ \end{minipage}
&  \begin{minipage}{0.7 cm}
\centering
\smallskip
$7$ \end{minipage}
&  \begin{minipage}{1.3 cm}
\centering
\smallskip
$12$ \end{minipage}
\\
\begin{minipage}{0.7 cm}
\centering
\smallskip
$ \mathcal{N}_1^2$  \end{minipage}
&   \begin{minipage}{6 cm}
\centering
\smallskip
$ T^3 + Y^3 + XYT $ \end{minipage}
&  \begin{minipage}{0.7 cm}
\centering
\smallskip
$5$ \end{minipage}
&  \begin{minipage}{1.3 cm}
\centering
\smallskip
$23$ \end{minipage}
\\
\begin{minipage}{0.7 cm}
\centering
\smallskip
$ \mathcal{N}_1^1$  \end{minipage}
&   \begin{minipage}{6 cm}
\centering
\smallskip
$ T^3 + XY^2 $ \end{minipage}
&  \begin{minipage}{0.7 cm}
\centering
\smallskip
$6$ \end{minipage}
&  \begin{minipage}{1.3 cm}
\centering
\smallskip
$19$ \end{minipage}
\\
\begin{minipage}{0.7 cm}
\centering
\smallskip
$ \mathcal{N}_0^0$  \end{minipage}
&   \begin{minipage}{6 cm}
\centering
\smallskip
$ 4T^3 + YT^2 + Y^2T + 3XT^2 + Y^3 + XY^2 $ \end{minipage}
&  \begin{minipage}{0.7 cm}
\centering
\smallskip
$7$ \end{minipage}
&  \begin{minipage}{1.3 cm}
\centering
\smallskip
$15$ \end{minipage}
\\
\hline
\end{tabular}
\end{minipage}\newline\phantom  \newline\newline


\noindent\begin{minipage}{\textwidth}
\centering
\textbf{Table 14}  $ q = 5$, non-singular cubic curves giving incomplete $(n,3)$-arcs  \medskip \\
\medskip
\begin{tabular}
{|c|c|c|c|c|}
\hline
type & equation &$n$ &
\begin{minipage}{1.3 cm}
\smallskip
\centering
$\#$ of \\
residual\\
points
\end{minipage}   &
\begin{minipage}{1 cm}
\smallskip
\centering
$j-$ \\
inv. \\
\end{minipage}  \\
\hline
\begin{minipage}{0.7 cm}
\centering
\smallskip
$ \mathcal{G}_3^2$  \end{minipage}
&   \begin{minipage}{8 cm}
\centering
\smallskip
$ 3T^3 + 4YT^2 + 4Y^2T + 4XT^2 + 3Y^3 + 4XYT + 4XY^2 + 4X^2T + 4X^2Y + 3X^3 $ \end{minipage}
&  \begin{minipage}{0.7 cm}
\centering
\smallskip
$3$ \end{minipage}
&  \begin{minipage}{1.3 cm}
\centering
\smallskip
$25$ \end{minipage}
&  \begin{minipage}{1 cm}
\centering
\smallskip
$2$ \end{minipage}
\\
\begin{minipage}{0.7 cm}
\centering
\smallskip
$ \overline{\mathcal{E}}_3^2$  \end{minipage}
&   \begin{minipage}{8 cm}
\centering
\smallskip
$ 2T^3 + XY^2 + X^2Y $ \end{minipage}
&  \begin{minipage}{0.7 cm}
\centering
\smallskip
$6$ \end{minipage}
&  \begin{minipage}{1.3 cm}
\centering
\smallskip
$13$ \end{minipage}
&  \begin{minipage}{1 cm}
\centering
\smallskip
$0$ \end{minipage}
\\
\begin{minipage}{0.7 cm}
\centering
\smallskip
$ \mathcal{E}_3^2$  \end{minipage}
&   \begin{minipage}{8 cm}
\centering
\smallskip
$ T^3 + 3YT^2 + 3Y^2T + 3XT^2 + Y^3 + 2XYT + 3XY^2 + 3X^2T + 3X^2Y + X^3 $ \end{minipage}
&  \begin{minipage}{0.7 cm}
\centering
\smallskip
$6$ \end{minipage}
&  \begin{minipage}{1.3 cm}
\centering
\smallskip
$13$ \end{minipage}
&  \begin{minipage}{1 cm}
\centering
\smallskip
$0$ \end{minipage}
\\
\begin{minipage}{0.7 cm}
\centering
\smallskip
$ \mathcal{H}_1^0$  \end{minipage}
&   \begin{minipage}{8 cm}
\centering
\smallskip
$ YT^2 + 3XY^2 + X^3 $ \end{minipage}
&  \begin{minipage}{0.7 cm}
\centering
\smallskip
$10$ \end{minipage}
&  \begin{minipage}{1.3 cm}
\centering
\smallskip
$1$ \end{minipage}
&  \begin{minipage}{1 cm}
\centering
\smallskip
$3$ \end{minipage}
\\
\begin{minipage}{0.7 cm}
\centering
\smallskip
$ \mathcal{H}_1^0$  \end{minipage}
&   \begin{minipage}{8 cm}
\centering
\smallskip
$ YT^2 + 4XY^2 + X^3 $ \end{minipage}
&  \begin{minipage}{0.7 cm}
\centering
\smallskip
$8$ \end{minipage}
&  \begin{minipage}{1.3 cm}
\centering
\smallskip
$4$ \end{minipage}
&  \begin{minipage}{1 cm}
\centering
\smallskip
$3$ \end{minipage}
\\
\begin{minipage}{0.7 cm}
\centering
\smallskip
$ \mathcal{H}_1^0$  \end{minipage}
&   \begin{minipage}{8 cm}
\centering
\smallskip
$ YT^2 + XY^2 + X^3 $ \end{minipage}
&  \begin{minipage}{0.7 cm}
\centering
\smallskip
$4$ \end{minipage}
&  \begin{minipage}{1.3 cm}
\centering
\smallskip
$24$ \end{minipage}
&  \begin{minipage}{1 cm}
\centering
\smallskip
$3$ \end{minipage}
\\
\begin{minipage}{0.7 cm}
\centering
\smallskip
$ \mathcal{H}_1^0$  \end{minipage}
&   \begin{minipage}{8 cm}
\centering
\smallskip
$ YT^2 + 2XY^2 + X^3 $ \end{minipage}
&  \begin{minipage}{0.7 cm}
\centering
\smallskip
$2$ \end{minipage}
&  \begin{minipage}{1.3 cm}
\centering
\smallskip
$29$ \end{minipage}
&  \begin{minipage}{1 cm}
\centering
\smallskip
$3$ \end{minipage}
\\
\begin{minipage}{0.7 cm}
\centering
\smallskip
$ \mathcal{G}_1$  \end{minipage}
&   \begin{minipage}{8 cm}
\centering
\smallskip
$ YT^2 + Y^3 + 4XY^2 + X^3 $ \end{minipage}
&  \begin{minipage}{0.7 cm}
\centering
\smallskip
$8$ \end{minipage}
&  \begin{minipage}{1.3 cm}
\centering
\smallskip
$5$ \end{minipage}
&  \begin{minipage}{1 cm}
\centering
\smallskip
$1$ \end{minipage}
\\
\begin{minipage}{0.7 cm}
\centering
\smallskip
$ \mathcal{G}_1$  \end{minipage}
&   \begin{minipage}{8 cm}
\centering
\smallskip
$ YT^2 + Y^3 + 2XY^2 + X^3 $ \end{minipage}
&  \begin{minipage}{0.7 cm}
\centering
\smallskip
$7$ \end{minipage}
&  \begin{minipage}{1.3 cm}
\centering
\smallskip
$10$ \end{minipage}
&  \begin{minipage}{1 cm}
\centering
\smallskip
$4$ \end{minipage}
\\
\begin{minipage}{0.7 cm}
\centering
\smallskip
$ \mathcal{G}_1$  \end{minipage}
&   \begin{minipage}{8 cm}
\centering
\smallskip
$ YT^2 + 2Y^3 + 3XY^2 + X^3 $ \end{minipage}
&  \begin{minipage}{0.7 cm}
\centering
\smallskip
$5$ \end{minipage}
&  \begin{minipage}{1.3 cm}
\centering
\smallskip
$20$ \end{minipage}
&  \begin{minipage}{1 cm}
\centering
\smallskip
$4$ \end{minipage}
\\
\begin{minipage}{0.7 cm}
\centering
\smallskip
$ \mathcal{G}_1$  \end{minipage}
&   \begin{minipage}{8 cm}
\centering
\smallskip
$ YT^2 + 2Y^3 + XY^2 + X^3 $ \end{minipage}
&  \begin{minipage}{0.7 cm}
\centering
\smallskip
$4$ \end{minipage}
&  \begin{minipage}{1.3 cm}
\centering
\smallskip
$24$ \end{minipage}
&  \begin{minipage}{1 cm}
\centering
\smallskip
$1$ \end{minipage}
\\
\begin{minipage}{0.7 cm}
\centering
\smallskip
$ \mathcal{G}_0^2$  \end{minipage}
&   \begin{minipage}{8 cm}
\centering
\smallskip
$ T^3 + 4Y^2T + 4Y^3 + XYT + 3XY^2 + 4X^2T + 4X^3 $ \end{minipage}
&  \begin{minipage}{0.7 cm}
\centering
\smallskip
$9$ \end{minipage}
&  \begin{minipage}{1.3 cm}
\centering
\smallskip
$4$ \end{minipage}
&  \begin{minipage}{1 cm}
\centering
\smallskip
$2$ \end{minipage}
\\
\begin{minipage}{0.7 cm}
\centering
\smallskip
$ \mathcal{G}_0^2$  \end{minipage}
&   \begin{minipage}{8 cm}
\centering
\smallskip
$ T^3 + 2Y^2T + 4Y^3 + 3XYT + 3XY^2 + 2X^2T + 4X^3 $ \end{minipage}
&  \begin{minipage}{0.7 cm}
\centering
\smallskip
$3$ \end{minipage}
&  \begin{minipage}{1.3 cm}
\centering
\smallskip
$28$ \end{minipage}
&  \begin{minipage}{1 cm}
\centering
\smallskip
$2$ \end{minipage}
\\
\begin{minipage}{0.7 cm}
\centering
\smallskip
$ \mathcal{E}_0^2$  \end{minipage}
&   \begin{minipage}{8 cm}
\centering
\smallskip
$ T^3 + 4Y^3 + 3XY^2 + 4X^3 $ \end{minipage}
&  \begin{minipage}{0.7 cm}
\centering
\smallskip
$6$ \end{minipage}
&  \begin{minipage}{1.3 cm}
\centering
\smallskip
$16$ \end{minipage}
&  \begin{minipage}{1 cm}
\centering
\smallskip
$0$ \end{minipage}
\\
\begin{minipage}{0.7 cm}
\centering
\smallskip
$ \mathcal{E}_0^2$  \end{minipage}
&   \begin{minipage}{8 cm}
\centering
\smallskip
$ T^3 + 3Y^2T + 4Y^3 + 2XYT + 3XY^2 + 3X^2T + 4X^3 $ \end{minipage}
&  \begin{minipage}{0.7 cm}
\centering
\smallskip
$6$ \end{minipage}
&  \begin{minipage}{1.3 cm}
\centering
\smallskip
$16$ \end{minipage}
&  \begin{minipage}{1 cm}
\centering
\smallskip
$0$ \end{minipage}
\\
\hline
\end{tabular}
\end{minipage}\newline\phantom    \newline\phantom   \newline The sets of $\mathbb{F}_q$-rational points of cubics $2T^{3} + XY^{2} + X^{2}Y $ and $T^{3} + 3YT^{2} +
3Y^{2}T + 3XT^{2} + Y^{3} + 2XYT + 3XY^{2} + 3X^{2}T + 3X^{2}Y + X^{3} $ are
equivalent.\newline The sets of $\mathbb{F}_q$-rational points of cubics $YT^{2} + 2Y^{3} +
XY^{2} + X^{3} $ and $YT^{2} + XY^{2} + X^{3} $ are equivalent.\newline%
\phantom  \newline\newline


\noindent\begin{minipage}{\textwidth}
\centering
\textbf{Table 15}  $ q = 5$, non-singular cubic curves giving complete $(n,3)$-arcs  \medskip \\
\medskip
\begin{tabular}
{|c|c|c|c|}
\hline
type & equation &$n$ &
\begin{minipage}{1 cm}
\smallskip
\centering
$j-$ \\
inv. \\
\smallskip
\end{minipage}  \\
\hline
\begin{minipage}{0.7 cm}
\centering
\smallskip
$ \mathcal{G}_3^2$  \end{minipage}
&   \begin{minipage}{9 cm}
\centering
\smallskip
$ 4T^3 + 2YT^2 + 2Y^2T + 2XT^2 + 4Y^3 + 2XY^2 + 2X^2T + 2X^2Y + 4X^3 $ \end{minipage}
&  \begin{minipage}{0.7 cm}
\centering
\smallskip
$9$ \end{minipage}
&  \begin{minipage}{1 cm}
\centering
\smallskip
$2$ \end{minipage}
\\
\hline
\end{tabular}
\end{minipage}\newline


$\; $ \phantom  \newline\newline\noindent\begin{minipage}{\textwidth}
\centering
\textbf{Table 16} $ q = 7$, a. i. singular cubic curves giving incomplete $(n,3)$-arcs  \medskip \\
\medskip
\begin{tabular}
{|c|c|c|c|}
\hline
type & equation &$n$ &
\begin{minipage}{1.3 cm}
\smallskip
\centering
$\#$ of \\
residual\\
points
\end{minipage}\\
\hline
\begin{minipage}{0.7 cm}
\centering
\smallskip
$ \mathcal{N}_3^2$  \end{minipage}
&   \begin{minipage}{6 cm}
\centering
\smallskip
$ T^3 + Y^3 + XYT $ \end{minipage}
&  \begin{minipage}{0.7 cm}
\centering
\smallskip
$7$ \end{minipage}
&  \begin{minipage}{1.3 cm}
\centering
\smallskip
$30$ \end{minipage}
\\
\begin{minipage}{0.7 cm}
\centering
\smallskip
$ \mathcal{N}_1^1$  \end{minipage}
&   \begin{minipage}{4 cm}
\centering
\smallskip
$ T^3 + XY^2 $ \end{minipage}
&  \begin{minipage}{0.7 cm}
\centering
\smallskip
$8$ \end{minipage}
&  \begin{minipage}{1.3 cm}
\centering
\smallskip
$25$ \end{minipage}
\\
\begin{minipage}{0.7 cm}
\centering
\smallskip
$ \mathcal{N}_1^0$  \end{minipage}
&   \begin{minipage}{6 cm}
\centering
\smallskip
$ 2YT^2 + 4XT^2 + Y^3 + XY^2 $ \end{minipage}
&  \begin{minipage}{0.7 cm}
\centering
\smallskip
$9$ \end{minipage}
&  \begin{minipage}{1.3 cm}
\centering
\smallskip
$16$ \end{minipage}
\\
\begin{minipage}{0.7 cm}
\centering
\smallskip
$ \mathcal{N}_0^2$  \end{minipage}
&   \begin{minipage}{6 cm}
\centering
\smallskip
$ 2T^3 + Y^3 + XYT $ \end{minipage}
&  \begin{minipage}{0.7 cm}
\centering
\smallskip
$7$ \end{minipage}
&  \begin{minipage}{1.3 cm}
\centering
\smallskip
$35$ \end{minipage}
\\
\hline
\end{tabular}
\end{minipage}\newline\phantom  \newline\newline


\noindent\begin{minipage}{\textwidth}
\centering
\textbf{Table 17}  $ q = 7$, non-singular cubic curves giving incomplete $(n,3)$-arcs  \medskip \\
\medskip
\begin{tabular}
{|c|c|c|c|c|}
\hline
type & equation &$n$ &
\begin{minipage}{1.3 cm}
\smallskip
\centering
$\#$ of \\
residual\\
points
\end{minipage}   &
\begin{minipage}{1 cm}
\smallskip
\centering
$j-$ \\
inv. \\
\end{minipage}  \\
\hline
\begin{minipage}{0.7 cm}
\centering
\smallskip
$ \mathcal{G}_3^1$  \end{minipage}
&   \begin{minipage}{8 cm}
\centering
\smallskip
$ 6T^3 + 4YT^2 + 4Y^2T + 4XT^2 + 6Y^3 + 2XYT + 4XY^2 + 4X^2T + 4X^2Y + 6X^3 $ \end{minipage}
&  \begin{minipage}{0.7 cm}
\centering
\smallskip
$9$ \end{minipage}
&  \begin{minipage}{1.3 cm}
\centering
\smallskip
$6$ \end{minipage}
&  \begin{minipage}{1 cm}
\centering
\smallskip
$3$ \end{minipage}
\\
\begin{minipage}{0.7 cm}
\centering
\smallskip
$ \mathcal{G}_3^1$  \end{minipage}
&   \begin{minipage}{8 cm}
\centering
\smallskip
$ 5T^3 + YT^2 + Y^2T + XT^2 + 5Y^3 + 3XYT + XY^2 + X^2T + X^2Y + 5X^3 $ \end{minipage}
&  \begin{minipage}{0.7 cm}
\centering
\smallskip
$6$ \end{minipage}
&  \begin{minipage}{1.3 cm}
\centering
\smallskip
$31$ \end{minipage}
&  \begin{minipage}{1 cm}
\centering
\smallskip
$5$ \end{minipage}
\\
\begin{minipage}{0.7 cm}
\centering
\smallskip
$ \mathcal{G}_3^1$  \end{minipage}
&   \begin{minipage}{8 cm}
\centering
\smallskip
$ 3T^3 + 2YT^2 + 2Y^2T + 2XT^2 + 3Y^3 + 5XYT + 2XY^2 + 2X^2T + 2X^2Y + 3X^3 $ \end{minipage}
&  \begin{minipage}{0.7 cm}
\centering
\smallskip
$6$ \end{minipage}
&  \begin{minipage}{1.3 cm}
\centering
\smallskip
$31$ \end{minipage}
&  \begin{minipage}{1 cm}
\centering
\smallskip
$4$ \end{minipage}
\\
\begin{minipage}{0.7 cm}
\centering
\smallskip
$ \overline{\mathcal{E}}_3^1$  \end{minipage}
&   \begin{minipage}{8 cm}
\centering
\smallskip
$ 3T^3 + XY^2 + X^2Y $ \end{minipage}
&  \begin{minipage}{0.7 cm}
\centering
\smallskip
$3$ \end{minipage}
&  \begin{minipage}{1.3 cm}
\centering
\smallskip
$49$ \end{minipage}
&  \begin{minipage}{1 cm}
\centering
\smallskip
$0$ \end{minipage}
\\
\begin{minipage}{0.7 cm}
\centering
\smallskip
$ \mathcal{E}_1^4$  \end{minipage}
&   \begin{minipage}{8 cm}
\centering
\smallskip
$ YT^2 + 2Y^3 + X^3 $ \end{minipage}
&  \begin{minipage}{0.7 cm}
\centering
\smallskip
$7$ \end{minipage}
&  \begin{minipage}{1.3 cm}
\centering
\smallskip
$26$ \end{minipage}
&  \begin{minipage}{1 cm}
\centering
\smallskip
$0$ \end{minipage}
\\
\begin{minipage}{0.7 cm}
\centering
\smallskip
$ \mathcal{E}_1^1$  \end{minipage}
&   \begin{minipage}{8 cm}
\centering
\smallskip
$ YT^2 + 4Y^3 + X^3 $ \end{minipage}
&  \begin{minipage}{0.7 cm}
\centering
\smallskip
$13$ \end{minipage}
&  \begin{minipage}{1.3 cm}
\centering
\smallskip
$1$ \end{minipage}
&  \begin{minipage}{1 cm}
\centering
\smallskip
$0$ \end{minipage}
\\
\begin{minipage}{0.7 cm}
\centering
\smallskip
$ \mathcal{E}_1^1$  \end{minipage}
&   \begin{minipage}{8 cm}
\centering
\smallskip
$ YT^2 + Y^3 + X^3 $ \end{minipage}
&  \begin{minipage}{0.7 cm}
\centering
\smallskip
$4$ \end{minipage}
&  \begin{minipage}{1.3 cm}
\centering
\smallskip
$48$ \end{minipage}
&  \begin{minipage}{1 cm}
\centering
\smallskip
$0$ \end{minipage}
\\
\begin{minipage}{0.7 cm}
\centering
\smallskip
$ \mathcal{H}_1^0$  \end{minipage}
&   \begin{minipage}{8 cm}
\centering
\smallskip
$ YT^2 + XY^2 + X^3 $ \end{minipage}
&  \begin{minipage}{0.7 cm}
\centering
\smallskip
$8$ \end{minipage}
&  \begin{minipage}{1.3 cm}
\centering
\smallskip
$17$ \end{minipage}
&  \begin{minipage}{1 cm}
\centering
\smallskip
$6$ \end{minipage}
\\
\begin{minipage}{0.7 cm}
\centering
\smallskip
$ \mathcal{H}_1^0$  \end{minipage}
&   \begin{minipage}{8 cm}
\centering
\smallskip
$ YT^2 + 3XY^2 + X^3 $ \end{minipage}
&  \begin{minipage}{0.7 cm}
\centering
\smallskip
$8$ \end{minipage}
&  \begin{minipage}{1.3 cm}
\centering
\smallskip
$16$ \end{minipage}
&  \begin{minipage}{1 cm}
\centering
\smallskip
$6$ \end{minipage}
\\
\begin{minipage}{0.7 cm}
\centering
\smallskip
$ \mathcal{G}_1$  \end{minipage}
&   \begin{minipage}{8 cm}
\centering
\smallskip
$ YT^2 + Y^3 + XY^2 + X^3 $ \end{minipage}
&  \begin{minipage}{0.7 cm}
\centering
\smallskip
$11$ \end{minipage}
&  \begin{minipage}{1.3 cm}
\centering
\smallskip
$3$ \end{minipage}
&  \begin{minipage}{1 cm}
\centering
\smallskip
$1$ \end{minipage}
\\
\begin{minipage}{0.7 cm}
\centering
\smallskip
$ \mathcal{G}_1$  \end{minipage}
&   \begin{minipage}{8 cm}
\centering
\smallskip
$ YT^2 + 3Y^3 + XY^2 + X^3 $ \end{minipage}
&  \begin{minipage}{0.7 cm}
\centering
\smallskip
$10$ \end{minipage}
&  \begin{minipage}{1.3 cm}
\centering
\smallskip
$5$ \end{minipage}
&  \begin{minipage}{1 cm}
\centering
\smallskip
$5$ \end{minipage}
\\
\begin{minipage}{0.7 cm}
\centering
\smallskip
$ \mathcal{G}_1$  \end{minipage}
&   \begin{minipage}{8 cm}
\centering
\smallskip
$ YT^2 + 3Y^3 + 3XY^2 + X^3 $ \end{minipage}
&  \begin{minipage}{0.7 cm}
\centering
\smallskip
$10$ \end{minipage}
&  \begin{minipage}{1.3 cm}
\centering
\smallskip
$4$ \end{minipage}
&  \begin{minipage}{1 cm}
\centering
\smallskip
$4$ \end{minipage}
\\
\begin{minipage}{0.7 cm}
\centering
\smallskip
$ \mathcal{G}_1$  \end{minipage}
&   \begin{minipage}{8 cm}
\centering
\smallskip
$ YT^2 + 2Y^3 + 3XY^2 + X^3 $ \end{minipage}
&  \begin{minipage}{0.7 cm}
\centering
\smallskip
$7$ \end{minipage}
&  \begin{minipage}{1.3 cm}
\centering
\smallskip
$26$ \end{minipage}
&  \begin{minipage}{1 cm}
\centering
\smallskip
$3$ \end{minipage}
\\
\begin{minipage}{0.7 cm}
\centering
\smallskip
$ \mathcal{G}_1$  \end{minipage}
&   \begin{minipage}{8 cm}
\centering
\smallskip
$ YT^2 + 6Y^3 + XY^2 + X^3 $ \end{minipage}
&  \begin{minipage}{0.7 cm}
\centering
\smallskip
$5$ \end{minipage}
&  \begin{minipage}{1.3 cm}
\centering
\smallskip
$42$ \end{minipage}
&  \begin{minipage}{1 cm}
\centering
\smallskip
$1$ \end{minipage}
\\
\begin{minipage}{0.7 cm}
\centering
\smallskip
$ \mathcal{G}_1$  \end{minipage}
&   \begin{minipage}{8 cm}
\centering
\smallskip
$ YT^2 + Y^3 + 3XY^2 + X^3 $ \end{minipage}
&  \begin{minipage}{0.7 cm}
\centering
\smallskip
$4$ \end{minipage}
&  \begin{minipage}{1.3 cm}
\centering
\smallskip
$48$ \end{minipage}
&  \begin{minipage}{1 cm}
\centering
\smallskip
$2$ \end{minipage}
\\
\begin{minipage}{0.7 cm}
\centering
\smallskip
$ \overline{\mathcal{E}}_0^4$  \end{minipage}
&   \begin{minipage}{8 cm}
\centering
\smallskip
$ 2T^3 + 3Y^3 + X^3 $ \end{minipage}
&  \begin{minipage}{0.7 cm}
\centering
\smallskip
$9$ \end{minipage}
&  \begin{minipage}{1.3 cm}
\centering
\smallskip
$9$ \end{minipage}
&  \begin{minipage}{1 cm}
\centering
\smallskip
$0$ \end{minipage}
\\
\begin{minipage}{0.7 cm}
\centering
\smallskip
$ \mathcal{E}_0^4$  \end{minipage}
&   \begin{minipage}{8 cm}
\centering
\smallskip
$ 2T^3 + 3Y^3 + 4XYT + X^3 $ \end{minipage}
&  \begin{minipage}{0.7 cm}
\centering
\smallskip
$9$ \end{minipage}
&  \begin{minipage}{1.3 cm}
\centering
\smallskip
$12$ \end{minipage}
&  \begin{minipage}{1 cm}
\centering
\smallskip
$0$ \end{minipage}
\\
\begin{minipage}{0.7 cm}
\centering
\smallskip
$ \mathcal{G}_0^1$  \end{minipage}
&   \begin{minipage}{8 cm}
\centering
\smallskip
$ 2T^3 + 2YT^2 + Y^3 + 4XYT + XY^2 + X^2T + 4X^3 $ \end{minipage}
&  \begin{minipage}{0.7 cm}
\centering
\smallskip
$9$ \end{minipage}
&  \begin{minipage}{1.3 cm}
\centering
\smallskip
$12$ \end{minipage}
&  \begin{minipage}{1 cm}
\centering
\smallskip
$3$ \end{minipage}
\\
\begin{minipage}{0.7 cm}
\centering
\smallskip
$ \mathcal{G}_0^1$  \end{minipage}
&   \begin{minipage}{8 cm}
\centering
\smallskip
$ 5T^3 + 3YT^2 + 4Y^3 + 2XYT + XY^2 + X^2T + 6X^3 $ \end{minipage}
&  \begin{minipage}{0.7 cm}
\centering
\smallskip
$6$ \end{minipage}
&  \begin{minipage}{1.3 cm}
\centering
\smallskip
$36$ \end{minipage}
&  \begin{minipage}{1 cm}
\centering
\smallskip
$5$ \end{minipage}
\\
\begin{minipage}{0.7 cm}
\centering
\smallskip
$ \mathcal{G}_0^1$  \end{minipage}
&   \begin{minipage}{8 cm}
\centering
\smallskip
$ 3T^3 + 2YT^2 + 5Y^3 + 6XYT + XY^2 + X^2T + 6X^3 $ \end{minipage}
&  \begin{minipage}{0.7 cm}
\centering
\smallskip
$6$ \end{minipage}
&  \begin{minipage}{1.3 cm}
\centering
\smallskip
$36$ \end{minipage}
&  \begin{minipage}{1 cm}
\centering
\smallskip
$4$ \end{minipage}
\\
\begin{minipage}{0.7 cm}
\centering
\smallskip
$ \mathcal{E}_0^1$  \end{minipage}
&   \begin{minipage}{8 cm}
\centering
\smallskip
$ 3YT^2 + XY^2 + X^2T $ \end{minipage}
&  \begin{minipage}{0.7 cm}
\centering
\smallskip
$3$ \end{minipage}
&  \begin{minipage}{1.3 cm}
\centering
\smallskip
$54$ \end{minipage}
&  \begin{minipage}{1 cm}
\centering
\smallskip
$0$ \end{minipage}
\\
\hline
\end{tabular}
\end{minipage}\newline\phantom    \newline\phantom   \newline The sets of $\mathbb{F}_q$-rational points of cubics $3T^{3} + 2YT^{2} + 2Y^{2}T + 2XT^{2} + 3Y^{3} +
5XYT + 2XY^{2} + 2X^{2}T + 2X^{2}Y + 3X^{3} $ and $5T^{3} + YT^{2} + Y^{2}T +
XT^{2} + 5Y^{3} + 3XYT + XY^{2} + X^{2}T + X^{2}Y + 5X^{3} $ are
equivalent.\newline The sets of $\mathbb{F}_q$-rational points of cubics $YT^{2} + Y^{3} +
3XY^{2} + X^{3} $ and $YT^{2} + Y^{3} + X^{3} $ are equivalent.\newline%
\phantom  \newline\newline


\noindent\begin{minipage}{\textwidth}
\centering
\textbf{Table 18}  $ q = 7$, non-singular cubic curves giving complete $(n,3)$-arcs  \medskip \\
\medskip
\begin{tabular}
{|c|c|c|c|}
\hline
type & equation &$n$ &
\begin{minipage}{1 cm}
\smallskip
\centering
$j-$ \\
inv. \\
\smallskip
\end{minipage}  \\
\hline
\begin{minipage}{0.7 cm}
\centering
\smallskip
$ \mathcal{E}_9^4$  \end{minipage}
&   \begin{minipage}{10 cm}
\centering
\smallskip
$ T^3 + Y^3 + X^3 $ \end{minipage}
&  \begin{minipage}{0.7 cm}
\centering
\smallskip
$9$ \end{minipage}
&  \begin{minipage}{1 cm}
\centering
\smallskip
$0$ \end{minipage}
\\
\begin{minipage}{0.7 cm}
\centering
\smallskip
$ \mathcal{G}_3^1$  \end{minipage}
&   \begin{minipage}{10 cm}
\centering
\smallskip
$ 4T^3 + 5YT^2 + 5Y^2T + 5XT^2 + 4Y^3 + 4XYT + 5XY^2 + 5X^2T + 5X^2Y + 4X^3 $ \end{minipage}
&  \begin{minipage}{0.7 cm}
\centering
\smallskip
$12$ \end{minipage}
&  \begin{minipage}{1 cm}
\centering
\smallskip
$2$ \end{minipage}
\\
\begin{minipage}{0.7 cm}
\centering
\smallskip
$ \overline{\mathcal{E}}_3^1$  \end{minipage}
&   \begin{minipage}{10 cm}
\centering
\smallskip
$ 2T^3 + XY^2 + X^2Y $ \end{minipage}
&  \begin{minipage}{0.7 cm}
\centering
\smallskip
$12$ \end{minipage}
&  \begin{minipage}{1 cm}
\centering
\smallskip
$0$ \end{minipage}
\\
\begin{minipage}{0.7 cm}
\centering
\smallskip
$ \mathcal{G}_0^1$  \end{minipage}
&   \begin{minipage}{10 cm}
\centering
\smallskip
$ T^3 + 3YT^2 + 5Y^3 + 6XYT + XY^2 + X^2T + 4X^3 $ \end{minipage}
&  \begin{minipage}{0.7 cm}
\centering
\smallskip
$12$ \end{minipage}
&  \begin{minipage}{1 cm}
\centering
\smallskip
$2$ \end{minipage}
\\
\begin{minipage}{0.7 cm}
\centering
\smallskip
$ \mathcal{E}_0^1$  \end{minipage}
&   \begin{minipage}{10 cm}
\centering
\smallskip
$ 2YT^2 + XY^2 + X^2T $ \end{minipage}
&  \begin{minipage}{0.7 cm}
\centering
\smallskip
$12$ \end{minipage}
&  \begin{minipage}{1 cm}
\centering
\smallskip
$0$ \end{minipage}
\\
\hline
\end{tabular}
\end{minipage}\newline


$\; $ \phantom  \newline\newline\noindent\begin{minipage}{\textwidth}
\centering
\textbf{Table 19} $ q = 8$, a. i. singular cubic curves giving incomplete $(n,3)$-arcs  \medskip \\
\medskip
\begin{tabular}
{|c|c|c|c|}
\hline
type & equation &$n$ &
\begin{minipage}{1.3 cm}
\smallskip
\centering
$\#$ of \\
residual\\
points
\end{minipage}\\
\hline
\begin{minipage}{0.7 cm}
\centering
\smallskip
$ \mathcal{N}_3^0$  \end{minipage}
&   \begin{minipage}{9 cm}
\centering
\smallskip
$ YT^2 + Y^2T + XT^2 + XYT + XY^2 $ \end{minipage}
&  \begin{minipage}{0.7 cm}
\centering
\smallskip
$10$ \end{minipage}
&  \begin{minipage}{1.3 cm}
\centering
\smallskip
$12$ \end{minipage}
\\
\begin{minipage}{0.7 cm}
\centering
\smallskip
$ \mathcal{N}_1^2$  \end{minipage}
&   \begin{minipage}{7 cm}
\centering
\smallskip
$ T^3 + Y^3 + XYT $ \end{minipage}
&  \begin{minipage}{0.7 cm}
\centering
\smallskip
$8$ \end{minipage}
&  \begin{minipage}{1.3 cm}
\centering
\smallskip
$36$ \end{minipage}
\\
\begin{minipage}{0.7 cm}
\centering
\smallskip
$ \mathcal{N}_1^1$  \end{minipage}
&   \begin{minipage}{9 cm}
\centering
\smallskip
$ T^3 + XY^2 $ \end{minipage}
&  \begin{minipage}{0.7 cm}
\centering
\smallskip
$9$ \end{minipage}
&  \begin{minipage}{1.3 cm}
\centering
\smallskip
$22$ \end{minipage}
\\
\begin{minipage}{0.7 cm}
\centering
\smallskip
$ \mathcal{N}_0^0$  \end{minipage}
&   \begin{minipage}{9 cm}
\centering
\smallskip
$ \xi^{2}T^3 + \xi^{6}YT^2 + Y^2T + XT^2 + \xi^{2}Y^3 + XYT + XY^2 $ \end{minipage}
&  \begin{minipage}{0.7 cm}
\centering
\smallskip
$10$ \end{minipage}
&  \begin{minipage}{1.3 cm}
\centering
\smallskip
$18$ \end{minipage}
\\
\hline
\end{tabular}
\end{minipage}\newline\phantom  \newline\newline


\noindent\begin{minipage}{\textwidth}
\centering
\textbf{Table 20}  $ q = 8$, non-singular cubic curves giving incomplete $(n,3)$-arcs  \medskip \\
\medskip
\begin{tabular}
{|c|c|c|c|c|}
\hline
type & equation &$n$ &
\begin{minipage}{1.3 cm}
\smallskip
\centering
$\#$ of \\
residual\\
points
\end{minipage}   &
\begin{minipage}{1 cm}
\smallskip
\centering
$j-$ \\
inv. \\
\end{minipage}  \\
\hline
\begin{minipage}{0.7 cm}
\centering
\smallskip
$ \mathcal{G}_3^2$  \end{minipage}
&   \begin{minipage}{8 cm}
\centering
\smallskip
$ \xi^{4}T^3 + \xi^{4}YT^2 + \xi^{4}Y^2T + \xi^{4}XT^2 + \xi^{4}Y^3 + XYT + \xi^{4}XY^2 + \xi^{4}X^2T + \xi^{4}X^2Y + \xi^{4}X^3 $ \end{minipage}
&  \begin{minipage}{0.7 cm}
\centering
\smallskip
$12$ \end{minipage}
&  \begin{minipage}{1.3 cm}
\centering
\smallskip
$4$ \end{minipage}
&  \begin{minipage}{1 cm}
\centering
\smallskip
$\xi^{4}$ \end{minipage}
\\
\begin{minipage}{0.7 cm}
\centering
\smallskip
$ \mathcal{G}_3^2$  \end{minipage}
&   \begin{minipage}{8 cm}
\centering
\smallskip
$ \xi^{2}T^3 + \xi^{2}YT^2 + \xi^{2}Y^2T + \xi^{2}XT^2 + \xi^{2}Y^3 + XYT + \xi^{2}XY^2 + \xi^{2}X^2T + \xi^{2}X^2Y + \xi^{2}X^3 $ \end{minipage}
&  \begin{minipage}{0.7 cm}
\centering
\smallskip
$12$ \end{minipage}
&  \begin{minipage}{1.3 cm}
\centering
\smallskip
$4$ \end{minipage}
&  \begin{minipage}{1 cm}
\centering
\smallskip
$\xi^{2}$ \end{minipage}
\\
\begin{minipage}{0.7 cm}
\centering
\smallskip
$ \mathcal{G}_3^2$  \end{minipage}
&   \begin{minipage}{8 cm}
\centering
\smallskip
$ \xi T^3 + \xi YT^2 + \xi Y^2T + \xi XT^2 + \xi Y^3 + XYT + \xi XY^2 + \xi X^2T + \xi X^2Y + \xi X^3 $ \end{minipage}
&  \begin{minipage}{0.7 cm}
\centering
\smallskip
$12$ \end{minipage}
&  \begin{minipage}{1.3 cm}
\centering
\smallskip
$4$ \end{minipage}
&  \begin{minipage}{1 cm}
\centering
\smallskip
$\xi$ \end{minipage}
\\
\begin{minipage}{0.7 cm}
\centering
\smallskip
$ \mathcal{G}_3^2$  \end{minipage}
&   \begin{minipage}{8 cm}
\centering
\smallskip
$ \xi^{6}T^3 + \xi^{6}YT^2 + \xi^{6}Y^2T + \xi^{6}XT^2 + \xi^{6}Y^3 + XYT + \xi^{6}XY^2 + \xi^{6}X^2T + \xi^{6}X^2Y + \xi^{6}X^3 $ \end{minipage}
&  \begin{minipage}{0.7 cm}
\centering
\smallskip
$6$ \end{minipage}
&  \begin{minipage}{1.3 cm}
\centering
\smallskip
$43$ \end{minipage}
&  \begin{minipage}{1 cm}
\centering
\smallskip
$\xi$ \end{minipage}
\\
\begin{minipage}{0.7 cm}
\centering
\smallskip
$ \mathcal{G}_3^2$  \end{minipage}
&   \begin{minipage}{8 cm}
\centering
\smallskip
$ \xi^{5}T^3 + \xi^{5}YT^2 + \xi^{5}Y^2T + \xi^{5}XT^2 + \xi^{5}Y^3 + XYT + \xi^{5}XY^2 + \xi^{5}X^2T + \xi^{5}X^2Y + \xi^{5}X^3 $ \end{minipage}
&  \begin{minipage}{0.7 cm}
\centering
\smallskip
$6$ \end{minipage}
&  \begin{minipage}{1.3 cm}
\centering
\smallskip
$43$ \end{minipage}
&  \begin{minipage}{1 cm}
\centering
\smallskip
$\xi^{2}$ \end{minipage}
\\
\begin{minipage}{0.7 cm}
\centering
\smallskip
$ \mathcal{G}_3^2$  \end{minipage}
&   \begin{minipage}{8 cm}
\centering
\smallskip
$ \xi^{3}T^3 + \xi^{3}YT^2 + \xi^{3}Y^2T + \xi^{3}XT^2 + \xi^{3}Y^3 + XYT + \xi^{3}XY^2 + \xi^{3}X^2T + \xi^{3}X^2Y + \xi^{3}X^3 $ \end{minipage}
&  \begin{minipage}{0.7 cm}
\centering
\smallskip
$6$ \end{minipage}
&  \begin{minipage}{1.3 cm}
\centering
\smallskip
$43$ \end{minipage}
&  \begin{minipage}{1 cm}
\centering
\smallskip
$\xi^{4}$ \end{minipage}
\\
\begin{minipage}{0.7 cm}
\centering
\smallskip
$ \overline{\mathcal{E}}_3^2$  \end{minipage}
&   \begin{minipage}{8 cm}
\centering
\smallskip
$ \xi T^3 + XY^2 + X^2Y $ \end{minipage}
&  \begin{minipage}{0.7 cm}
\centering
\smallskip
$9$ \end{minipage}
&  \begin{minipage}{1.3 cm}
\centering
\smallskip
$13$ \end{minipage}
&  \begin{minipage}{1 cm}
\centering
\smallskip
$0$ \end{minipage}
\\
\begin{minipage}{0.7 cm}
\centering
\smallskip
$ \mathcal{G}_1^0$  \end{minipage}
&   \begin{minipage}{8 cm}
\centering
\smallskip
$ YT^2 + XYT + \xi^{4}XY^2 + X^2Y + X^3 $ \end{minipage}
&  \begin{minipage}{0.7 cm}
\centering
\smallskip
$10$ \end{minipage}
&  \begin{minipage}{1.3 cm}
\centering
\smallskip
$9$ \end{minipage}
&  \begin{minipage}{1 cm}
\centering
\smallskip
$\xi^{6}$ \end{minipage}
\\
\begin{minipage}{0.7 cm}
\centering
\smallskip
$ \mathcal{G}_1^0$  \end{minipage}
&   \begin{minipage}{8 cm}
\centering
\smallskip
$ YT^2 + XYT + \xi^{2}XY^2 + X^2Y + X^3 $ \end{minipage}
&  \begin{minipage}{0.7 cm}
\centering
\smallskip
$10$ \end{minipage}
&  \begin{minipage}{1.3 cm}
\centering
\smallskip
$9$ \end{minipage}
&  \begin{minipage}{1 cm}
\centering
\smallskip
$\xi^{3}$ \end{minipage}
\\
\begin{minipage}{0.7 cm}
\centering
\smallskip
$ \mathcal{G}_1^0$  \end{minipage}
&   \begin{minipage}{8 cm}
\centering
\smallskip
$ YT^2 + XYT + \xi XY^2 + X^2Y + X^3 $ \end{minipage}
&  \begin{minipage}{0.7 cm}
\centering
\smallskip
$10$ \end{minipage}
&  \begin{minipage}{1.3 cm}
\centering
\smallskip
$9$ \end{minipage}
&  \begin{minipage}{1 cm}
\centering
\smallskip
$\xi^{5}$ \end{minipage}
\\
\begin{minipage}{0.7 cm}
\centering
\smallskip
$ \mathcal{G}_1^0$  \end{minipage}
&   \begin{minipage}{8 cm}
\centering
\smallskip
$ YT^2 + XYT + \xi^{4}XY^2 + X^3 $ \end{minipage}
&  \begin{minipage}{0.7 cm}
\centering
\smallskip
$8$ \end{minipage}
&  \begin{minipage}{1.3 cm}
\centering
\smallskip
$26$ \end{minipage}
&  \begin{minipage}{1 cm}
\centering
\smallskip
$\xi^{6}$ \end{minipage}
\\
\begin{minipage}{0.7 cm}
\centering
\smallskip
$ \mathcal{G}_1^0$  \end{minipage}
&   \begin{minipage}{8 cm}
\centering
\smallskip
$ YT^2 + XYT + \xi^{2}XY^2 + X^3 $ \end{minipage}
&  \begin{minipage}{0.7 cm}
\centering
\smallskip
$8$ \end{minipage}
&  \begin{minipage}{1.3 cm}
\centering
\smallskip
$26$ \end{minipage}
&  \begin{minipage}{1 cm}
\centering
\smallskip
$\xi^{3}$ \end{minipage}
\\
\begin{minipage}{0.7 cm}
\centering
\smallskip
$ \mathcal{G}_1^0$  \end{minipage}
&   \begin{minipage}{8 cm}
\centering
\smallskip
$ YT^2 + XYT + \xi XY^2 + X^3 $ \end{minipage}
&  \begin{minipage}{0.7 cm}
\centering
\smallskip
$8$ \end{minipage}
&  \begin{minipage}{1.3 cm}
\centering
\smallskip
$26$ \end{minipage}
&  \begin{minipage}{1 cm}
\centering
\smallskip
$\xi^{5}$ \end{minipage}
\\
\begin{minipage}{0.7 cm}
\centering
\smallskip
$ \mathcal{G}_1^0$  \end{minipage}
&   \begin{minipage}{8 cm}
\centering
\smallskip
$ YT^2 + XYT + XY^2 + X^3 $ \end{minipage}
&  \begin{minipage}{0.7 cm}
\centering
\smallskip
$4$ \end{minipage}
&  \begin{minipage}{1.3 cm}
\centering
\smallskip
$63$ \end{minipage}
&  \begin{minipage}{1 cm}
\centering
\smallskip
$1$ \end{minipage}
\\
\begin{minipage}{0.7 cm}
\centering
\smallskip
$ \mathcal{E}_1^0$  \end{minipage}
&   \begin{minipage}{8 cm}
\centering
\smallskip
$ YT^2 + Y^2T + XY^2 + X^3 $ \end{minipage}
&  \begin{minipage}{0.7 cm}
\centering
\smallskip
$5$ \end{minipage}
&  \begin{minipage}{1.3 cm}
\centering
\smallskip
$56$ \end{minipage}
&  \begin{minipage}{1 cm}
\centering
\smallskip
$0$ \end{minipage}
\\
\begin{minipage}{0.7 cm}
\centering
\smallskip
$ \mathcal{G}_0^2$  \end{minipage}
&   \begin{minipage}{8 cm}
\centering
\smallskip
$ T^3 + \xi^{6}Y^2T + \xi^{4}Y^3 + \xi^{3}XYT + XY^2 + \xi^{6}X^2T + X^3 $ \end{minipage}
&  \begin{minipage}{0.7 cm}
\centering
\smallskip
$6$ \end{minipage}
&  \begin{minipage}{1.3 cm}
\centering
\smallskip
$49$ \end{minipage}
&  \begin{minipage}{1 cm}
\centering
\smallskip
$\xi$ \end{minipage}
\\
\begin{minipage}{0.7 cm}
\centering
\smallskip
$ \mathcal{G}_0^2$  \end{minipage}
&   \begin{minipage}{8 cm}
\centering
\smallskip
$ T^3 + \xi^{6}Y^2T + \xi Y^3 + XYT + XY^2 + \xi^{6}X^2T + X^3 $ \end{minipage}
&  \begin{minipage}{0.7 cm}
\centering
\smallskip
$6$ \end{minipage}
&  \begin{minipage}{1.3 cm}
\centering
\smallskip
$49$ \end{minipage}
&  \begin{minipage}{1 cm}
\centering
\smallskip
$\xi^{4}$ \end{minipage}
\\
\begin{minipage}{0.7 cm}
\centering
\smallskip
$ \mathcal{G}_0^2$  \end{minipage}
&   \begin{minipage}{8 cm}
\centering
\smallskip
$ T^3 + \xi^{5}Y^2T + \xi Y^3 + \xi^{6}XYT + XY^2 + \xi^{5}X^2T + X^3 $ \end{minipage}
&  \begin{minipage}{0.7 cm}
\centering
\smallskip
$6$ \end{minipage}
&  \begin{minipage}{1.3 cm}
\centering
\smallskip
$49$ \end{minipage}
&  \begin{minipage}{1 cm}
\centering
\smallskip
$\xi^{2}$ \end{minipage}
\\
\begin{minipage}{0.7 cm}
\centering
\smallskip
$ \mathcal{E}_0^2$  \end{minipage}
&   \begin{minipage}{8 cm}
\centering
\smallskip
$ T^3 + \xi^{4}Y^3 + XY^2 + X^3 $ \end{minipage}
&  \begin{minipage}{0.7 cm}
\centering
\smallskip
$9$ \end{minipage}
&  \begin{minipage}{1.3 cm}
\centering
\smallskip
$19$ \end{minipage}
&  \begin{minipage}{1 cm}
\centering
\smallskip
$0$ \end{minipage}
\\
\hline
\end{tabular}
\end{minipage}\newline\phantom    \newline\phantom   \newline The sets of $\mathbb{F}_q$-rational points of cubics $\xi^{6}T^{3} + \xi^{6}YT^{2} + \xi^{6}Y^{2}T +
\xi^{6}XT^{2} + \xi^{6}Y^{3} + XYT + \xi^{6}XY^{2} + \xi^{6}X^{2}T + \xi
^{6}X^{2}Y + \xi^{6}X^{3} $ and $\xi^{3}T^{3} + \xi^{3}YT^{2} + \xi^{3}Y^{2}T
+ \xi^{3}XT^{2} + \xi^{3}Y^{3} + XYT + \xi^{3}XY^{2} + \xi^{3}X^{2}T + \xi
^{3}X^{2}Y + \xi^{3}X^{3} $ are equivalent.\newline The sets of $\mathbb{F}_q$-rational
points of cubics $\xi^{6}T^{3} + \xi^{6}YT^{2} + \xi^{6}Y^{2}T + \xi^{6}XT^{2}
+ \xi^{6}Y^{3} + XYT + \xi^{6}XY^{2} + \xi^{6}X^{2}T + \xi^{6}X^{2}Y + \xi
^{6}X^{3} $ and $\xi^{5}T^{3} + \xi^{5}YT^{2} + \xi^{5}Y^{2}T + \xi^{5}XT^{2}
+ \xi^{5}Y^{3} + XYT + \xi^{5}XY^{2} + \xi^{5}X^{2}T + \xi^{5}X^{2}Y + \xi
^{5}X^{3} $ are equivalent.\newline\phantom  \newline\newline\newline\newline


\noindent\begin{minipage}{\textwidth}
\centering
\textbf{Table 21}  $ q = 8$, non-singular cubic curves giving complete $(n,3)$-arcs  \medskip \\
\medskip
\begin{tabular}
{|c|c|c|c|}
\hline
type & equation &$n$ &
\begin{minipage}{1 cm}
\smallskip
\centering
$j-$ \\
inv. \\
\smallskip
\end{minipage}  \\
\hline
\begin{minipage}{0.7 cm}
\centering
\smallskip
$ \mathcal{G}_1^0$  \end{minipage}
&   \begin{minipage}{7 cm}
\centering
\smallskip
$ YT^2 + XYT + XY^2 + X^2Y + X^3 $ \end{minipage}
&  \begin{minipage}{0.7 cm}
\centering
\smallskip
$14$ \end{minipage}
&  \begin{minipage}{1 cm}
\centering
\smallskip
$1$ \end{minipage}
\\
\begin{minipage}{0.7 cm}
\centering
\smallskip
$ \mathcal{E}_1^0$  \end{minipage}
&   \begin{minipage}{7 cm}
\centering
\smallskip
$ YT^2 + Y^2T + Y^3 + XY^2 + X^3 $ \end{minipage}
&  \begin{minipage}{0.7 cm}
\centering
\smallskip
$13$ \end{minipage}
&  \begin{minipage}{1 cm}
\centering
\smallskip
$0$ \end{minipage}
\\
\begin{minipage}{0.7 cm}
\centering
\smallskip
$ \mathcal{G}_0^2$  \end{minipage}
&   \begin{minipage}{7 cm}
\centering
\smallskip
$ T^3 + Y^2T + \xi^{4}Y^3 + \xi^{4}XYT + XY^2 + X^2T + X^3 $ \end{minipage}
&  \begin{minipage}{0.7 cm}
\centering
\smallskip
$12$ \end{minipage}
&  \begin{minipage}{1 cm}
\centering
\smallskip
$\xi^{2}$ \end{minipage}
\\
\begin{minipage}{0.7 cm}
\centering
\smallskip
$ \mathcal{G}_0^2$  \end{minipage}
&   \begin{minipage}{7 cm}
\centering
\smallskip
$ T^3 + Y^2T + \xi^{2}Y^3 + \xi^{2}XYT + XY^2 + X^2T + X^3 $ \end{minipage}
&  \begin{minipage}{0.7 cm}
\centering
\smallskip
$12$ \end{minipage}
&  \begin{minipage}{1 cm}
\centering
\smallskip
$\xi$ \end{minipage}
\\
\begin{minipage}{0.7 cm}
\centering
\smallskip
$ \mathcal{G}_0^2$  \end{minipage}
&   \begin{minipage}{7 cm}
\centering
\smallskip
$ T^3 + Y^2T + \xi Y^3 + \xi XYT + XY^2 + X^2T + X^3 $ \end{minipage}
&  \begin{minipage}{0.7 cm}
\centering
\smallskip
$12$ \end{minipage}
&  \begin{minipage}{1 cm}
\centering
\smallskip
$\xi^{4}$ \end{minipage}
\\
\hline
\end{tabular}
\end{minipage}\newline


$\; $ \phantom  \newline\newline\newline\newline\noindent
\begin{minipage}{\textwidth}
\centering
\textbf{Table 22} $ q = 9$, a. i. singular cubic curves giving incomplete $(n,3)$-arcs  \medskip \\
\medskip
\begin{tabular}
{|c|c|c|c|}
\hline
type & equation &$n$ &
\begin{minipage}{1.3 cm}
\smallskip
\centering
$\#$ of \\
residual\\
points
\end{minipage}\\
\hline
\begin{minipage}{0.7 cm}
\centering
\smallskip
$ \mathcal{N}_q^1$  \end{minipage}
&   \begin{minipage}{3 cm}
\centering
\smallskip
$ T^3 + XY^2 $ \end{minipage}
&  \begin{minipage}{0.7 cm}
\centering
\smallskip
$10$ \end{minipage}
&  \begin{minipage}{1.3 cm}
\centering
\smallskip
$5$ \end{minipage}
\\
\begin{minipage}{0.7 cm}
\centering
\smallskip
$ \mathcal{N}_1^2$  \end{minipage}
&   \begin{minipage}{3 cm}
\centering
\smallskip
$ T^3 + Y^3 + XYT $ \end{minipage}
&  \begin{minipage}{0.7 cm}
\centering
\smallskip
$9$ \end{minipage}
&  \begin{minipage}{1.3 cm}
\centering
\smallskip
$36$ \end{minipage}
\\
\begin{minipage}{0.7 cm}
\centering
\smallskip
$ \mathcal{N}_1^0$  \end{minipage}
&   \begin{minipage}{3 cm}
\centering
\smallskip
$ \xi XT^2 + Y^3 + XY^2 $ \end{minipage}
&  \begin{minipage}{0.7 cm}
\centering
\smallskip
$11$ \end{minipage}
&  \begin{minipage}{1.3 cm}
\centering
\smallskip
$14$ \end{minipage}
\\
\begin{minipage}{0.7 cm}
\centering
\smallskip
$ \mathcal{N}_0^1$  \end{minipage}
&   \begin{minipage}{3 cm}
\centering
\smallskip
$ T^3 + YT^2 + XY^2 $ \end{minipage}
&  \begin{minipage}{0.7 cm}
\centering
\smallskip
$10$ \end{minipage}
&  \begin{minipage}{1.3 cm}
\centering
\smallskip
$27$ \end{minipage}
\\
\hline
\end{tabular}
\end{minipage}\newline\phantom  \newline\newline\newline\newline


\noindent\begin{minipage}{\textwidth}
\centering
\textbf{Table 23}  $ q = 9$, non-singular cubic curves giving incomplete $(n,3)$-arcs  \medskip \\
\medskip
\begin{tabular}
{|c|c|c|c|c|}
\hline
type & equation &$n$ &
\begin{minipage}{1.3 cm}
\smallskip
\centering
$\#$ of \\
residual\\
points
\end{minipage}   &
\begin{minipage}{1 cm}
\smallskip
\centering
$j-$ \\
inv. \\
\end{minipage}  \\
\hline
\begin{minipage}{0.7 cm}
\centering
\smallskip
$ \mathcal{G}_3$  \end{minipage}
&   \begin{minipage}{7 cm}
\centering
\smallskip
$ \xi^{7}T^3 + \xi^{7}Y^3 + XYT + \xi^{7}X^3 $ \end{minipage}
&  \begin{minipage}{0.7 cm}
\centering
\smallskip
$12$ \end{minipage}
&  \begin{minipage}{1.3 cm}
\centering
\smallskip
$3$ \end{minipage}
&  \begin{minipage}{1 cm}
\centering
\smallskip
$\xi^{7}$ \end{minipage}
\\
\begin{minipage}{0.7 cm}
\centering
\smallskip
$ \mathcal{G}_3$  \end{minipage}
&   \begin{minipage}{7 cm}
\centering
\smallskip
$ \xi^{5}T^3 + \xi^{5}Y^3 + XYT + \xi^{5}X^3 $ \end{minipage}
&  \begin{minipage}{0.7 cm}
\centering
\smallskip
$12$ \end{minipage}
&  \begin{minipage}{1.3 cm}
\centering
\smallskip
$3$ \end{minipage}
&  \begin{minipage}{1 cm}
\centering
\smallskip
$\xi^{5}$ \end{minipage}
\\
\begin{minipage}{0.7 cm}
\centering
\smallskip
$ \mathcal{G}_3$  \end{minipage}
&   \begin{minipage}{7 cm}
\centering
\smallskip
$ \xi^{6}T^3 + \xi^{6}Y^3 + XYT + \xi^{6}X^3 $ \end{minipage}
&  \begin{minipage}{0.7 cm}
\centering
\smallskip
$9$ \end{minipage}
&  \begin{minipage}{1.3 cm}
\centering
\smallskip
$21$ \end{minipage}
&  \begin{minipage}{1 cm}
\centering
\smallskip
$\xi^{2}$ \end{minipage}
\\
\begin{minipage}{0.7 cm}
\centering
\smallskip
$ \mathcal{G}_3$  \end{minipage}
&   \begin{minipage}{7 cm}
\centering
\smallskip
$ \xi^{2}T^3 + \xi^{2}Y^3 + XYT + \xi^{2}X^3 $ \end{minipage}
&  \begin{minipage}{0.7 cm}
\centering
\smallskip
$9$ \end{minipage}
&  \begin{minipage}{1.3 cm}
\centering
\smallskip
$21$ \end{minipage}
&  \begin{minipage}{1 cm}
\centering
\smallskip
$\xi^{6}$ \end{minipage}
\\
\begin{minipage}{0.7 cm}
\centering
\smallskip
$ \mathcal{G}_3$  \end{minipage}
&   \begin{minipage}{7 cm}
\centering
\smallskip
$ \xi^{3}T^3 + \xi^{3}Y^3 + XYT + \xi^{3}X^3 $ \end{minipage}
&  \begin{minipage}{0.7 cm}
\centering
\smallskip
$6$ \end{minipage}
&  \begin{minipage}{1.3 cm}
\centering
\smallskip
$57$ \end{minipage}
&  \begin{minipage}{1 cm}
\centering
\smallskip
$\xi^{3}$ \end{minipage}
\\
\begin{minipage}{0.7 cm}
\centering
\smallskip
$ \mathcal{G}_3$  \end{minipage}
&   \begin{minipage}{7 cm}
\centering
\smallskip
$ \xi T^3 + \xi Y^3 + XYT + \xi X^3 $ \end{minipage}
&  \begin{minipage}{0.7 cm}
\centering
\smallskip
$6$ \end{minipage}
&  \begin{minipage}{1.3 cm}
\centering
\smallskip
$57$ \end{minipage}
&  \begin{minipage}{1 cm}
\centering
\smallskip
$\xi$ \end{minipage}
\\
\begin{minipage}{0.7 cm}
\centering
\smallskip
$ \mathcal{S}_1^3$  \end{minipage}
&   \begin{minipage}{7 cm}
\centering
\smallskip
$ YT^2 + \xi^{7}Y^3 + \xi^{2}XY^2 + X^3 $ \end{minipage}
&  \begin{minipage}{0.7 cm}
\centering
\smallskip
$4$ \end{minipage}
&  \begin{minipage}{1.3 cm}
\centering
\smallskip
$80$ \end{minipage}
&  \begin{minipage}{1 cm}
\centering
\smallskip
$0$ \end{minipage}
\\
\begin{minipage}{0.7 cm}
\centering
\smallskip
$ \mathcal{S}_1^1$  \end{minipage}
&   \begin{minipage}{7 cm}
\centering
\smallskip
$ YT^2 + Y^3 + \xi^{5}XY^2 + X^3 $ \end{minipage}
&  \begin{minipage}{0.7 cm}
\centering
\smallskip
$10$ \end{minipage}
&  \begin{minipage}{1.3 cm}
\centering
\smallskip
$15$ \end{minipage}
&  \begin{minipage}{1 cm}
\centering
\smallskip
$0$ \end{minipage}
\\
\begin{minipage}{0.7 cm}
\centering
\smallskip
$ \mathcal{S}_1^1$  \end{minipage}
&   \begin{minipage}{7 cm}
\centering
\smallskip
$ YT^2 + Y^3 + \xi^{3}XY^2 + X^3 $ \end{minipage}
&  \begin{minipage}{0.7 cm}
\centering
\smallskip
$10$ \end{minipage}
&  \begin{minipage}{1.3 cm}
\centering
\smallskip
$15$ \end{minipage}
&  \begin{minipage}{1 cm}
\centering
\smallskip
$0$ \end{minipage}
\\
\begin{minipage}{0.7 cm}
\centering
\smallskip
$ \mathcal{S}_1^0$  \end{minipage}
&   \begin{minipage}{7 cm}
\centering
\smallskip
$ YT^2 + Y^3 + \xi^{2}XY^2 + X^3 $ \end{minipage}
&  \begin{minipage}{0.7 cm}
\centering
\smallskip
$13$ \end{minipage}
&  \begin{minipage}{1.3 cm}
\centering
\smallskip
$7$ \end{minipage}
&  \begin{minipage}{1 cm}
\centering
\smallskip
$0$ \end{minipage}
\\
\begin{minipage}{0.7 cm}
\centering
\smallskip
$ \mathcal{S}_1^0$  \end{minipage}
&   \begin{minipage}{7 cm}
\centering
\smallskip
$ YT^2 + Y^3 + \xi^{4}XY^2 + X^3 $ \end{minipage}
&  \begin{minipage}{0.7 cm}
\centering
\smallskip
$7$ \end{minipage}
&  \begin{minipage}{1.3 cm}
\centering
\smallskip
$50$ \end{minipage}
&  \begin{minipage}{1 cm}
\centering
\smallskip
$0$ \end{minipage}
\\
\begin{minipage}{0.7 cm}
\centering
\smallskip
$ \mathcal{G}_1$  \end{minipage}
&   \begin{minipage}{7 cm}
\centering
\smallskip
$ YT^2 + \xi^{5}Y^3 + \xi^{5}X^2Y + X^3 $ \end{minipage}
&  \begin{minipage}{0.7 cm}
\centering
\smallskip
$11$ \end{minipage}
&  \begin{minipage}{1.3 cm}
\centering
\smallskip
$8$ \end{minipage}
&  \begin{minipage}{1 cm}
\centering
\smallskip
$\xi^{6}$ \end{minipage}
\\
\begin{minipage}{0.7 cm}
\centering
\smallskip
$ \mathcal{G}_1$  \end{minipage}
&   \begin{minipage}{7 cm}
\centering
\smallskip
$ YT^2 + \xi Y^3 + \xi^{5}X^2Y + X^3 $ \end{minipage}
&  \begin{minipage}{0.7 cm}
\centering
\smallskip
$11$ \end{minipage}
&  \begin{minipage}{1.3 cm}
\centering
\smallskip
$8$ \end{minipage}
&  \begin{minipage}{1 cm}
\centering
\smallskip
$\xi^{2}$ \end{minipage}
\\
\begin{minipage}{0.7 cm}
\centering
\smallskip
$ \mathcal{G}_1$  \end{minipage}
&   \begin{minipage}{7 cm}
\centering
\smallskip
$ YT^2 + \xi^{7}Y^3 + \xi^{5}X^2Y + X^3 $ \end{minipage}
&  \begin{minipage}{0.7 cm}
\centering
\smallskip
$8$ \end{minipage}
&  \begin{minipage}{1.3 cm}
\centering
\smallskip
$36$ \end{minipage}
&  \begin{minipage}{1 cm}
\centering
\smallskip
$\xi^{4}$ \end{minipage}
\\
\begin{minipage}{0.7 cm}
\centering
\smallskip
$ \mathcal{G}_1$  \end{minipage}
&   \begin{minipage}{7 cm}
\centering
\smallskip
$ YT^2 + \xi^{6}Y^3 + \xi^{5}X^2Y + X^3 $ \end{minipage}
&  \begin{minipage}{0.7 cm}
\centering
\smallskip
$8$ \end{minipage}
&  \begin{minipage}{1.3 cm}
\centering
\smallskip
$37$ \end{minipage}
&  \begin{minipage}{1 cm}
\centering
\smallskip
$\xi^{5}$ \end{minipage}
\\
\begin{minipage}{0.7 cm}
\centering
\smallskip
$ \mathcal{G}_1$  \end{minipage}
&   \begin{minipage}{7 cm}
\centering
\smallskip
$ YT^2 + \xi^{4}Y^3 + \xi^{5}X^2Y + X^3 $ \end{minipage}
&  \begin{minipage}{0.7 cm}
\centering
\smallskip
$8$ \end{minipage}
&  \begin{minipage}{1.3 cm}
\centering
\smallskip
$37$ \end{minipage}
&  \begin{minipage}{1 cm}
\centering
\smallskip
$\xi^{7}$ \end{minipage}
\\
\begin{minipage}{0.7 cm}
\centering
\smallskip
$ \mathcal{G}_1$  \end{minipage}
&   \begin{minipage}{7 cm}
\centering
\smallskip
$ YT^2 + \xi^{3}Y^3 + \xi^{5}X^2Y + X^3 $ \end{minipage}
&  \begin{minipage}{0.7 cm}
\centering
\smallskip
$5$ \end{minipage}
&  \begin{minipage}{1.3 cm}
\centering
\smallskip
$72$ \end{minipage}
&  \begin{minipage}{1 cm}
\centering
\smallskip
$1$ \end{minipage}
\\
\begin{minipage}{0.7 cm}
\centering
\smallskip
$ \mathcal{G}_0$  \end{minipage}
&   \begin{minipage}{7 cm}
\centering
\smallskip
$ \xi^{6}T^3 + \xi^{2}YT^2 + Y^3 + \xi^{2}XY^2 + \xi^{2}X^2T + \xi^{4}X^2Y + X^3 $ \end{minipage}
&  \begin{minipage}{0.7 cm}
\centering
\smallskip
$12$ \end{minipage}
&  \begin{minipage}{1.3 cm}
\centering
\smallskip
$3$ \end{minipage}
&  \begin{minipage}{1 cm}
\centering
\smallskip
$\xi^{5}$ \end{minipage}
\\
\begin{minipage}{0.7 cm}
\centering
\smallskip
$ \mathcal{G}_0$  \end{minipage}
&   \begin{minipage}{7 cm}
\centering
\smallskip
$ \xi^{5}T^3 + \xi^{2}YT^2 + Y^3 + \xi^{2}XY^2 + \xi^{2}X^2T + \xi^{4}X^2Y + X^3 $ \end{minipage}
&  \begin{minipage}{0.7 cm}
\centering
\smallskip
$12$ \end{minipage}
&  \begin{minipage}{1.3 cm}
\centering
\smallskip
$1$ \end{minipage}
&  \begin{minipage}{1 cm}
\centering
\smallskip
$\xi^{4}$ \end{minipage}
\\
\begin{minipage}{0.7 cm}
\centering
\smallskip
$ \mathcal{G}_0$  \end{minipage}
&   \begin{minipage}{7 cm}
\centering
\smallskip
$ \xi^{3}T^3 + \xi^{2}YT^2 + Y^3 + \xi^{2}XY^2 + \xi^{2}X^2T + \xi^{4}X^2Y + X^3 $ \end{minipage}
&  \begin{minipage}{0.7 cm}
\centering
\smallskip
$12$ \end{minipage}
&  \begin{minipage}{1.3 cm}
\centering
\smallskip
$3$ \end{minipage}
&  \begin{minipage}{1 cm}
\centering
\smallskip
$\xi^{7}$ \end{minipage}
\\
\begin{minipage}{0.7 cm}
\centering
\smallskip
$ \mathcal{G}_0$  \end{minipage}
&   \begin{minipage}{7 cm}
\centering
\smallskip
$ \xi^{4}T^3 + \xi^{2}YT^2 + Y^3 + \xi^{2}XY^2 + \xi^{2}X^2T + \xi^{4}X^2Y + X^3 $ \end{minipage}
&  \begin{minipage}{0.7 cm}
\centering
\smallskip
$9$ \end{minipage}
&  \begin{minipage}{1.3 cm}
\centering
\smallskip
$28$ \end{minipage}
&  \begin{minipage}{1 cm}
\centering
\smallskip
$\xi^{6}$ \end{minipage}
\\
\begin{minipage}{0.7 cm}
\centering
\smallskip
$ \mathcal{G}_0$  \end{minipage}
&   \begin{minipage}{7 cm}
\centering
\smallskip
$ \xi^{2}YT^2 + Y^3 + \xi^{2}XY^2 + \xi^{2}X^2T + \xi^{4}X^2Y + X^3 $ \end{minipage}
&  \begin{minipage}{0.7 cm}
\centering
\smallskip
$9$ \end{minipage}
&  \begin{minipage}{1.3 cm}
\centering
\smallskip
$28$ \end{minipage}
&  \begin{minipage}{1 cm}
\centering
\smallskip
$\xi^{2}$ \end{minipage}
\\
\begin{minipage}{0.7 cm}
\centering
\smallskip
$ \mathcal{G}_0$  \end{minipage}
&   \begin{minipage}{7 cm}
\centering
\smallskip
$ \xi^{2}T^3 + \xi^{2}YT^2 + Y^3 + \xi^{2}XY^2 + \xi^{2}X^2T + \xi^{4}X^2Y + X^3 $ \end{minipage}
&  \begin{minipage}{0.7 cm}
\centering
\smallskip
$6$ \end{minipage}
&  \begin{minipage}{1.3 cm}
\centering
\smallskip
$64$ \end{minipage}
&  \begin{minipage}{1 cm}
\centering
\smallskip
$\xi^{3}$ \end{minipage}
\\
\begin{minipage}{0.7 cm}
\centering
\smallskip
$ \mathcal{G}_0$  \end{minipage}
&   \begin{minipage}{7 cm}
\centering
\smallskip
$ \xi T^3 + \xi^{2}YT^2 + Y^3 + \xi^{2}XY^2 + \xi^{2}X^2T + \xi^{4}X^2Y + X^3 $ \end{minipage}
&  \begin{minipage}{0.7 cm}
\centering
\smallskip
$6$ \end{minipage}
&  \begin{minipage}{1.3 cm}
\centering
\smallskip
$64$ \end{minipage}
&  \begin{minipage}{1 cm}
\centering
\smallskip
$\xi$ \end{minipage}
\\
\hline
\end{tabular}
\end{minipage}\newline\phantom    \newline\phantom   \newline The sets of $\mathbb{F}_q$-rational points of cubics $\xi T^{3} + \xi Y^{3} + XYT + \xi X^{3} $ and
$\xi^{3}T^{3} + \xi^{3}Y^{3} + XYT + \xi^{3}X^{3} $ are equivalent.\newline
The sets of $\mathbb{F}_q$-rational points of cubics $\xi^{2}T^{3} + \xi^{2}Y^{3} + XYT +
\xi^{2}X^{3} $ and $\xi^{6}T^{3} + \xi^{6}Y^{3} + XYT + \xi^{6}X^{3} $ are
equivalent.\newline The sets of $\mathbb{F}_q$-rational points of cubics $\xi^{2}T^{3} +
\xi^{2}YT^{2} + Y^{3} + \xi^{2}XY^{2} + \xi^{2}X^{2}T + \xi^{4}X^{2}Y + X^{3}
$ and $\xi T^{3} + \xi^{2}YT^{2} + Y^{3} + \xi^{2}XY^{2} + \xi^{2}X^{2}T +
\xi^{4}X^{2}Y + X^{3} $ are equivalent.\newline\phantom  \newline%
\newline\newline\newline


\noindent\begin{minipage}{\textwidth}
\centering
\textbf{Table 24}  $ q = 9$, non-singular cubic curves giving complete $(n,3)$-arcs  \medskip \\
\medskip
\begin{tabular}
{|c|c|c|c|}
\hline
type & equation &$n$ &
\begin{minipage}{1 cm}
\smallskip
\centering
$j-$ \\
inv. \\
\smallskip
\end{minipage}  \\
\hline
\begin{minipage}{0.7 cm}
\centering
\smallskip
$ \mathcal{G}_3$  \end{minipage}
&   \begin{minipage}{9 cm}
\centering
\smallskip
$ \xi^{4}T^3 + \xi^{4}Y^3 + XYT + \xi^{4}X^3 $ \end{minipage}
&  \begin{minipage}{0.7 cm}
\centering
\smallskip
$15$ \end{minipage}
&  \begin{minipage}{1 cm}
\centering
\smallskip
$1$ \end{minipage}
\\
\begin{minipage}{0.7 cm}
\centering
\smallskip
$ \mathcal{G}_3$  \end{minipage}
&   \begin{minipage}{7 cm}
\centering
\smallskip
$ T^3 + Y^3 + XYT + X^3 $ \end{minipage}
&  \begin{minipage}{0.7 cm}
\centering
\smallskip
$12$ \end{minipage}
&  \begin{minipage}{1 cm}
\centering
\smallskip
$\xi^{4}$ \end{minipage}
\\
\begin{minipage}{0.7 cm}
\centering
\smallskip
$ \mathcal{S}_1^3$  \end{minipage}
&   \begin{minipage}{7 cm}
\centering
\smallskip
$ YT^2 + \xi^{2}Y^3 + \xi^{4}XY^2 + X^3 $ \end{minipage}
&  \begin{minipage}{0.7 cm}
\centering
\smallskip
$16$ \end{minipage}
&  \begin{minipage}{1 cm}
\centering
\smallskip
$0$ \end{minipage}
\\
\begin{minipage}{0.7 cm}
\centering
\smallskip
$ \mathcal{G}_1$  \end{minipage}
&   \begin{minipage}{7 cm}
\centering
\smallskip
$ YT^2 + Y^3 + \xi^{5}X^2Y + X^3 $ \end{minipage}
&  \begin{minipage}{0.7 cm}
\centering
\smallskip
$14$ \end{minipage}
&  \begin{minipage}{1 cm}
\centering
\smallskip
$\xi^{3}$ \end{minipage}
\\
\begin{minipage}{0.7 cm}
\centering
\smallskip
$ \mathcal{G}_1$  \end{minipage}
&   \begin{minipage}{7 cm}
\centering
\smallskip
$ YT^2 + \xi^{2}Y^3 + \xi^{5}X^2Y + X^3 $ \end{minipage}
&  \begin{minipage}{0.7 cm}
\centering
\smallskip
$14$ \end{minipage}
&  \begin{minipage}{1 cm}
\centering
\smallskip
$\xi$ \end{minipage}
\\
\hline
\end{tabular}
\end{minipage}\newline


$\; $ \phantom  \newline\newline\newline\newline\noindent
\begin{minipage}{\textwidth}
\centering
\textbf{Table 25} $ q = 11$, a. i. singular cubic curves giving incomplete $(n,3)$-arcs  \medskip \\
\medskip
\begin{tabular}
{|c|c|c|c|}
\hline
type & equation &$n$ &
\begin{minipage}{1.3 cm}
\smallskip
\centering
$\#$ of \\
residual\\
points
\end{minipage}\\
\hline
\begin{minipage}{0.7 cm}
\centering
\smallskip
$ \mathcal{N}_3^0$  \end{minipage}
&   \begin{minipage}{8 cm}
\centering
\smallskip
$ 2YT^2 + 3XT^2 + Y^3 + XY^2 $ \end{minipage}
&  \begin{minipage}{0.7 cm}
\centering
\smallskip
$13$ \end{minipage}
&  \begin{minipage}{1.3 cm}
\centering
\smallskip
$9$ \end{minipage}
\\
\begin{minipage}{0.7 cm}
\centering
\smallskip
$ \mathcal{N}_1^2$  \end{minipage}
&   \begin{minipage}{6 cm}
\centering
\smallskip
$ T^3 + Y^3 + XYT $ \end{minipage}
&  \begin{minipage}{0.7 cm}
\centering
\smallskip
$11$ \end{minipage}
&  \begin{minipage}{1.3 cm}
\centering
\smallskip
$31$ \end{minipage}
\\
\begin{minipage}{0.7 cm}
\centering
\smallskip
$ \mathcal{N}_1^1$  \end{minipage}
&   \begin{minipage}{6 cm}
\centering
\smallskip
$ T^3 + XY^2 $ \end{minipage}
&  \begin{minipage}{0.7 cm}
\centering
\smallskip
$12$ \end{minipage}
&  \begin{minipage}{1.3 cm}
\centering
\smallskip
$21$ \end{minipage}
\\
\begin{minipage}{0.7 cm}
\centering
\smallskip
$ \mathcal{N}_0^0$  \end{minipage}
&   \begin{minipage}{8 cm}
\centering
\smallskip
$ 10T^3 + 2YT^2 + Y^2T + 3XT^2 + Y^3 + XY^2 $ \end{minipage}
&  \begin{minipage}{0.7 cm}
\centering
\smallskip
$13$ \end{minipage}
&  \begin{minipage}{1.3 cm}
\centering
\smallskip
$12$ \end{minipage}
\\
\hline
\end{tabular}
\end{minipage}\newline\phantom  \newline\newline\newline\newline


\noindent\begin{minipage}{\textwidth}
\centering
\textbf{Table 26}  $ q = 11$, non-singular cubic curves giving incomplete $(n,3)$-arcs  \medskip \\
\medskip
\begin{tabular}
{|c|c|c|c|c|}
\hline
type & equation &$n$ &
\begin{minipage}{1.3 cm}
\smallskip
\centering
$\#$ of \\
residual\\
points
\end{minipage}   &
\begin{minipage}{1 cm}
\smallskip
\centering
$j-$ \\
inv. \\
\end{minipage}  \\
\hline
\begin{minipage}{0.7 cm}
\centering
\smallskip
$ \mathcal{H}_3^2$  \end{minipage}
&   \begin{minipage}{8 cm}
\centering
\smallskip
$ T^3 + 3YT^2 + 3Y^2T + 3XT^2 + Y^3 + 7XYT + 3XY^2 + 3X^2T + 3X^2Y + X^3 $ \end{minipage}
&  \begin{minipage}{0.7 cm}
\centering
\smallskip
$12$ \end{minipage}
&  \begin{minipage}{1.3 cm}
\centering
\smallskip
$6$ \end{minipage}
&  \begin{minipage}{1 cm}
\centering
\smallskip
$1$ \end{minipage}
\\
\begin{minipage}{0.7 cm}
\centering
\smallskip
$ \mathcal{H}_3^2$  \end{minipage}
&   \begin{minipage}{8 cm}
\centering
\smallskip
$ 8T^3 + 2YT^2 + 2Y^2T + 2XT^2 + 8Y^3 + 5XYT + 2XY^2 + 2X^2T + 2X^2Y + 8X^3 $ \end{minipage}
&  \begin{minipage}{0.7 cm}
\centering
\smallskip
$12$ \end{minipage}
&  \begin{minipage}{1.3 cm}
\centering
\smallskip
$7$ \end{minipage}
&  \begin{minipage}{1 cm}
\centering
\smallskip
$1$ \end{minipage}
\\
\begin{minipage}{0.7 cm}
\centering
\smallskip
$ \mathcal{G}_3^2$  \end{minipage}
&   \begin{minipage}{8 cm}
\centering
\smallskip
$ 9T^3 + 5YT^2 + 5Y^2T + 5XT^2 + 9Y^3 + 5XY^2 + 5X^2T + 5X^2Y + 9X^3 $ \end{minipage}
&  \begin{minipage}{0.7 cm}
\centering
\smallskip
$15$ \end{minipage}
&  \begin{minipage}{1.3 cm}
\centering
\smallskip
$1$ \end{minipage}
&  \begin{minipage}{1 cm}
\centering
\smallskip
$10$ \end{minipage}
\\
\begin{minipage}{0.7 cm}
\centering
\smallskip
$ \mathcal{G}_3^2$  \end{minipage}
&   \begin{minipage}{8 cm}
\centering
\smallskip
$ 6T^3 + 7YT^2 + 7Y^2T + 7XT^2 + 6Y^3 + 4XYT + 7XY^2 + 7X^2T + 7X^2Y + 6X^3 $ \end{minipage}
&  \begin{minipage}{0.7 cm}
\centering
\smallskip
$9$ \end{minipage}
&  \begin{minipage}{1.3 cm}
\centering
\smallskip
$42$ \end{minipage}
&  \begin{minipage}{1 cm}
\centering
\smallskip
$10$ \end{minipage}
\\
\hline
\end{tabular}
\end{minipage}\newline\phantom  \newline\newline\newline\newline


\noindent\begin{minipage}{\textwidth}
\centering
\textbf{Table 27}  $ q = 11$, non-singular cubic curves giving complete $(n,3)$-arcs  \medskip \\
\medskip
\begin{tabular}
{|c|c|c|c|}
\hline
type & equation &$n$ &
\begin{minipage}{1 cm}
\smallskip
\centering
$j-$ \\
inv. \\
\smallskip
\end{minipage}  \\
\hline
\begin{minipage}{0.7 cm}
\centering
\smallskip
$ \mathcal{G}_3^2$  \end{minipage}
&   \begin{minipage}{10 cm}
\centering
\smallskip
$ 3T^3 + 9YT^2 + 9Y^2T + 9XT^2 + 3Y^3 + 8XYT + 9XY^2 + 9X^2T + 9X^2Y + 3X^3 $ \end{minipage}
&  \begin{minipage}{0.7 cm}
\centering
\smallskip
$18$ \end{minipage}
&  \begin{minipage}{1 cm}
\centering
\smallskip
$3$ \end{minipage}
\\
\begin{minipage}{0.7 cm}
\centering
\smallskip
$ \mathcal{G}_3^2$  \end{minipage}
&   \begin{minipage}{10 cm}
\centering
\smallskip
$ 7T^3 + 10YT^2 + 10Y^2T + 10XT^2 + 7Y^3 + 10XYT + 10XY^2 + 10X^2T + 10X^2Y + 7X^3 $ \end{minipage}
&  \begin{minipage}{0.7 cm}
\centering
\smallskip
$15$ \end{minipage}
&  \begin{minipage}{1 cm}
\centering
\smallskip
$4$ \end{minipage}
\\
\begin{minipage}{0.7 cm}
\centering
\smallskip
$ \mathcal{G}_1$  \end{minipage}
&   \begin{minipage}{10 cm}
\centering
\smallskip
$ YT^2 + 7Y^3 + 2XY^2 + X^3 $ \end{minipage}
&  \begin{minipage}{0.7 cm}
\centering
\smallskip
$17$ \end{minipage}
&  \begin{minipage}{1 cm}
\centering
\smallskip
$5$ \end{minipage}
\\
\begin{minipage}{0.7 cm}
\centering
\smallskip
$ \mathcal{G}_1$  \end{minipage}
&   \begin{minipage}{10 cm}
\centering
\smallskip
$ YT^2 + 9Y^3 + XY^2 + X^3 $ \end{minipage}
&  \begin{minipage}{0.7 cm}
\centering
\smallskip
$16$ \end{minipage}
&  \begin{minipage}{1 cm}
\centering
\smallskip
$2$ \end{minipage}
\\
\begin{minipage}{0.7 cm}
\centering
\smallskip
$ \mathcal{G}_1$  \end{minipage}
&   \begin{minipage}{10 cm}
\centering
\smallskip
$ YT^2 + 10Y^3 + 2XY^2 + X^3 $ \end{minipage}
&  \begin{minipage}{0.7 cm}
\centering
\smallskip
$16$ \end{minipage}
&  \begin{minipage}{1 cm}
\centering
\smallskip
$8$ \end{minipage}
\\
\hline
\end{tabular}
\end{minipage}\newline


$\; $ \phantom  \newline\newline\newline\newline\noindent
\begin{minipage}{\textwidth}
\centering
\textbf{Table 28} $ q = 13$, a. i. singular cubic curves giving incomplete $(n,3)$-arcs  \medskip \\
\medskip
\begin{tabular}
{|c|c|c|c|}
\hline
type & equation &$n$ &
\begin{minipage}{1.3 cm}
\smallskip
\centering
$\#$ of \\
residual\\
points
\end{minipage}\\
\hline
\begin{minipage}{0.7 cm}
\centering
\smallskip
$ \mathcal{N}_3^2$  \end{minipage}
&   \begin{minipage}{6 cm}
\centering
\smallskip
$ T^3 + Y^3 + XYT $ \end{minipage}
&  \begin{minipage}{0.7 cm}
\centering
\smallskip
$13$ \end{minipage}
&  \begin{minipage}{1.3 cm}
\centering
\smallskip
$21$ \end{minipage}
\\
\begin{minipage}{0.7 cm}
\centering
\smallskip
$ \mathcal{N}_1^1$  \end{minipage}
&   \begin{minipage}{4 cm}
\centering
\smallskip
$ T^3 + XY^2 $ \end{minipage}
&  \begin{minipage}{0.7 cm}
\centering
\smallskip
$14$ \end{minipage}
&  \begin{minipage}{1.3 cm}
\centering
\smallskip
$16$ \end{minipage}
\\
\begin{minipage}{0.7 cm}
\centering
\smallskip
$ \mathcal{N}_1^0$  \end{minipage}
&   \begin{minipage}{6 cm}
\centering
\smallskip
$ 6YT^2 + 11XT^2 + Y^3 + XY^2 $ \end{minipage}
&  \begin{minipage}{0.7 cm}
\centering
\smallskip
$15$ \end{minipage}
&  \begin{minipage}{1.3 cm}
\centering
\smallskip
$5$ \end{minipage}
\\
\begin{minipage}{0.7 cm}
\centering
\smallskip
$ \mathcal{N}_0^2$  \end{minipage}
&   \begin{minipage}{4 cm}
\centering
\smallskip
$ 2T^3 + Y^3 + XYT $ \end{minipage}
&  \begin{minipage}{0.7 cm}
\centering
\smallskip
$13$ \end{minipage}
&  \begin{minipage}{1.3 cm}
\centering
\smallskip
$28$ \end{minipage}
\\
\hline
\end{tabular}
\end{minipage}\newline\phantom  \newline\newline\newline\newline


\noindent\begin{minipage}{\textwidth}
\centering
\textbf{Table 29}  $ q = 13$, non-singular cubic curves giving incomplete $(n,3)$-arcs  \medskip \\
\medskip
\begin{tabular}
{|c|c|c|c|c|}
\hline
type & equation &$n$ &
\begin{minipage}{1.3 cm}
\smallskip
\centering
$\#$ of \\
residual\\
points
\end{minipage}   &
\begin{minipage}{1 cm}
\smallskip
\centering
$j-$ \\
inv. \\
\end{minipage}  \\
\hline
\begin{minipage}{0.7 cm}
\centering
\smallskip
$ \mathcal{E}_9^4$  \end{minipage}
&   \begin{minipage}{8 cm}
\centering
\smallskip
$ T^3 + Y^3 + X^3 $ \end{minipage}
&  \begin{minipage}{0.7 cm}
\centering
\smallskip
$9$ \end{minipage}
&  \begin{minipage}{1.3 cm}
\centering
\smallskip
$54$ \end{minipage}
&  \begin{minipage}{1 cm}
\centering
\smallskip
$0$ \end{minipage}
\\
\begin{minipage}{0.7 cm}
\centering
\smallskip
$ \mathcal{G}_3^1$  \end{minipage}
&   \begin{minipage}{8 cm}
\centering
\smallskip
$ T^3 + 3YT^2 + 3Y^2T + 3XT^2 + Y^3 + 7XYT + 3XY^2 + 3X^2T + 3X^2Y + X^3 $ \end{minipage}
&  \begin{minipage}{0.7 cm}
\centering
\smallskip
$18$ \end{minipage}
&  \begin{minipage}{1.3 cm}
\centering
\smallskip
$1$ \end{minipage}
&  \begin{minipage}{1 cm}
\centering
\smallskip
$7$ \end{minipage}
\\
\begin{minipage}{0.7 cm}
\centering
\smallskip
$ \mathcal{G}_3^1$  \end{minipage}
&   \begin{minipage}{8 cm}
\centering
\smallskip
$ 9T^3 + YT^2 + Y^2T + XT^2 + 9Y^3 + 3XYT + XY^2 + X^2T + X^2Y + 9X^3 $ \end{minipage}
&  \begin{minipage}{0.7 cm}
\centering
\smallskip
$15$ \end{minipage}
&  \begin{minipage}{1.3 cm}
\centering
\smallskip
$3$ \end{minipage}
&  \begin{minipage}{1 cm}
\centering
\smallskip
$9$ \end{minipage}
\\
\begin{minipage}{0.7 cm}
\centering
\smallskip
$ \mathcal{G}_3^1$  \end{minipage}
&   \begin{minipage}{8 cm}
\centering
\smallskip
$ 3T^3 + 9YT^2 + 9Y^2T + 9XT^2 + 3Y^3 + 6XYT + 9XY^2 + 9X^2T + 9X^2Y + 3X^3 $ \end{minipage}
&  \begin{minipage}{0.7 cm}
\centering
\smallskip
$15$ \end{minipage}
&  \begin{minipage}{1.3 cm}
\centering
\smallskip
$3$ \end{minipage}
&  \begin{minipage}{1 cm}
\centering
\smallskip
$2$ \end{minipage}
\\
\hline
\end{tabular}
\end{minipage}\newline\newline\newline\newline\noindent
\begin{minipage}{\textwidth}
\centering
\textbf{Table 29} Continue $1$  \medskip \\
\medskip
\begin{tabular}
{|c|c|c|c|c|}
\hline
type & equation &$n$ &
\begin{minipage}{1.3 cm}
\smallskip
\centering
$\#$ of \\
residual\\
points
\end{minipage}
& $j-$inv. \\
\hline
\begin{minipage}{0.7 cm}
\centering
\smallskip
$ \mathcal{G}_0^1$  \end{minipage}
&   \begin{minipage}{8 cm}
\centering
\smallskip
$ T^3 + 4YT^2 + 10Y^3 + 9XYT + XY^2 + X^2T + 9X^3 $ \end{minipage}
&  \begin{minipage}{0.7 cm}
\centering
\smallskip
$15$ \end{minipage}
&  \begin{minipage}{1.3 cm}
\centering
\smallskip
$3$ \end{minipage}
&  \begin{minipage}{1 cm}
\centering
\smallskip
$2$ \end{minipage}
\\
\begin{minipage}{0.7 cm}
\centering
\smallskip
$ \mathcal{G}_0^1$  \end{minipage}
&   \begin{minipage}{8 cm}
\centering
\smallskip
$ 9T^3 + 2YT^2 + 11Y^3 + 6XYT + XY^2 + X^2T + 12X^3 $ \end{minipage}
&  \begin{minipage}{0.7 cm}
\centering
\smallskip
$15$ \end{minipage}
&  \begin{minipage}{1.3 cm}
\centering
\smallskip
$6$ \end{minipage}
&  \begin{minipage}{1 cm}
\centering
\smallskip
$9$ \end{minipage}
\\
\begin{minipage}{0.7 cm}
\centering
\smallskip
$ \mathcal{G}_0^1$  \end{minipage}
&   \begin{minipage}{8 cm}
\centering
\smallskip
$ 5T^3 + 4YT^2 + 11Y^3 + 6XYT + XY^2 + X^2T + 6X^3 $ \end{minipage}
&  \begin{minipage}{0.7 cm}
\centering
\smallskip
$12$ \end{minipage}
&  \begin{minipage}{1.3 cm}
\centering
\smallskip
$24$ \end{minipage}
&  \begin{minipage}{1 cm}
\centering
\smallskip
$11$ \end{minipage}
\\
\begin{minipage}{0.7 cm}
\centering
\smallskip
$ \mathcal{G}_0^1$  \end{minipage}
&   \begin{minipage}{8 cm}
\centering
\smallskip
$ 5T^3 + 2YT^2 + 9Y^3 + 12XYT + XY^2 + X^2T + 11X^3 $ \end{minipage}
&  \begin{minipage}{0.7 cm}
\centering
\smallskip
$12$ \end{minipage}
&  \begin{minipage}{1.3 cm}
\centering
\smallskip
$24$ \end{minipage}
&  \begin{minipage}{1 cm}
\centering
\smallskip
$4$ \end{minipage}
\\
\hline
\end{tabular}
\end{minipage}\newline\phantom  \newline\newline\newline\newline


\noindent\begin{minipage}{\textwidth}
\centering
\textbf{Table 30}  $ q = 13$, non-singular cubic curves giving complete $(n,3)$-arcs  \medskip \\
\medskip
\begin{tabular}
{|c|c|c|c|}
\hline
type & equation &$n$ &
\begin{minipage}{1 cm}
\smallskip
\centering
$j-$ \\
inv. \\
\smallskip
\end{minipage}  \\
\hline
\begin{minipage}{0.7 cm}
\centering
\smallskip
$ \mathcal{H}_9^4$  \end{minipage}
&   \begin{minipage}{10 cm}
\centering
\smallskip
$ T^3 + Y^3 + 11XYT + X^3 $ \end{minipage}
&  \begin{minipage}{0.7 cm}
\centering
\smallskip
$18$ \end{minipage}
&  \begin{minipage}{1 cm}
\centering
\smallskip
$12$ \end{minipage}
\\
\begin{minipage}{0.7 cm}
\centering
\smallskip
$ \mathcal{G}_3^1$  \end{minipage}
&   \begin{minipage}{10 cm}
\centering
\smallskip
$ 8T^3 + 11YT^2 + 11Y^2T + 11XT^2 + 8Y^3 + 10XYT + 11XY^2 + 11X^2T + 11X^2Y + 8X^3 $ \end{minipage}
&  \begin{minipage}{0.7 cm}
\centering
\smallskip
$18$ \end{minipage}
&  \begin{minipage}{1 cm}
\centering
\smallskip
$10$ \end{minipage}
\\
\begin{minipage}{0.7 cm}
\centering
\smallskip
$ \overline{\mathcal{E}}_3^1$  \end{minipage}
&   \begin{minipage}{10 cm}
\centering
\smallskip
$ 4T^3 + XY^2 + X^2Y $ \end{minipage}
&  \begin{minipage}{0.7 cm}
\centering
\smallskip
$21$ \end{minipage}
&  \begin{minipage}{1 cm}
\centering
\smallskip
$0$ \end{minipage}
\\
\begin{minipage}{0.7 cm}
\centering
\smallskip
$ \mathcal{E}_1^4$  \end{minipage}
&   \begin{minipage}{10 cm}
\centering
\smallskip
$ YT^2 + 11Y^3 + X^3 $ \end{minipage}
&  \begin{minipage}{0.7 cm}
\centering
\smallskip
$19$ \end{minipage}
&  \begin{minipage}{1 cm}
\centering
\smallskip
$0$ \end{minipage}
\\
\begin{minipage}{0.7 cm}
\centering
\smallskip
$ \mathcal{H}_1^0$  \end{minipage}
&   \begin{minipage}{10 cm}
\centering
\smallskip
$ YT^2 + XY^2 + X^3 $ \end{minipage}
&  \begin{minipage}{0.7 cm}
\centering
\smallskip
$20$ \end{minipage}
&  \begin{minipage}{1 cm}
\centering
\smallskip
$12$ \end{minipage}
\\
\hline
\end{tabular}
\end{minipage}\newline


$\; $ \phantom  \newline\newline\newline\newline\noindent
\begin{minipage}{\textwidth}
\centering
\textbf{Table 31} $ q = 16$, a. i. singular cubic curves giving incomplete $(n,3)$-arcs  \medskip \\
\medskip
\begin{tabular}
{|c|c|c|c|}
\hline
type & equation &$n$ &
\begin{minipage}{1.3 cm}
\smallskip
\centering
$\#$ of \\
residual\\
points
\end{minipage}\\
\hline
\begin{minipage}{0.7 cm}
\centering
\smallskip
$ \mathcal{N}_3^2$  \end{minipage}
&   \begin{minipage}{9 cm}
\centering
\smallskip
$ T^3 + Y^3 + XYT $ \end{minipage}
&  \begin{minipage}{0.7 cm}
\centering
\smallskip
$16$ \end{minipage}
&  \begin{minipage}{1.3 cm}
\centering
\smallskip
$12$ \end{minipage}
\\
\begin{minipage}{0.7 cm}
\centering
\smallskip
$ \mathcal{N}_1^0$  \end{minipage}
&   \begin{minipage}{9 cm}
\centering
\smallskip
$ \xi^{14}YT^2 + \xi^{7}Y^2T + \xi^{7}XT^2 + \xi^{9}Y^3 + XYT + XY^2 $ \end{minipage}
&  \begin{minipage}{0.7 cm}
\centering
\smallskip
$18$ \end{minipage}
&  \begin{minipage}{1.3 cm}
\centering
\smallskip
$1$ \end{minipage}
\\
\begin{minipage}{0.7 cm}
\centering
\smallskip
$ \mathcal{N}_0^2$  \end{minipage}
&   \begin{minipage}{6 cm}
\centering
\smallskip
$ \xi^{7}T^3 + Y^3 + XYT $ \end{minipage}
&  \begin{minipage}{0.7 cm}
\centering
\smallskip
$16$ \end{minipage}
&  \begin{minipage}{1.3 cm}
\centering
\smallskip
$15$ \end{minipage}
\\
\hline
\end{tabular}
\end{minipage}\newline


$\; $ \phantom  \newline\newline\newline\newline\noindent
\begin{minipage}{\textwidth}
\centering
\textbf{Table 32} $ q = 16$, a. i. singular cubic curves giving complete $(n,3)$-arcs  \medskip \\
\medskip
\begin{tabular}
{|c|c|c|}
\hline
type & equation &$n$ \\
\hline
\begin{minipage}{0.7 cm}
\centering
\smallskip
$ \mathcal{N}_1^1$  \end{minipage}
&   \begin{minipage}{2 cm}
\centering
\smallskip
$ T^3 + XY^2 $ \end{minipage}
&  \begin{minipage}{0.7 cm}
\centering
\smallskip
$17$ \end{minipage}
\\
\hline
\end{tabular}
\end{minipage}\newline\phantom  \newline\newline\newline\newline


\noindent\begin{minipage}{\textwidth}
\centering
\textbf{Table 33}  $ q = 16$, non-singular cubic curves giving incomplete $(n,3)$-arcs  \medskip \\
\medskip
\begin{tabular}
{|c|c|c|c|c|}
\hline
type & equation &$n$ &
\begin{minipage}{1.3 cm}
\smallskip
\centering
$\#$ of \\
residual\\
points
\end{minipage}   &
\begin{minipage}{1 cm}
\smallskip
\centering
$j-$ \\
inv. \\
\end{minipage}  \\
\hline
\begin{minipage}{0.7 cm}
\centering
\smallskip
$ \mathcal{G}_9^4$  \end{minipage}
&   \begin{minipage}{8 cm}
\centering
\smallskip
$ T^3 + Y^3 + \xi^{7}XYT + X^3 $ \end{minipage}
&  \begin{minipage}{0.7 cm}
\centering
\smallskip
$18$ \end{minipage}
&  \begin{minipage}{1.3 cm}
\centering
\smallskip
$9$ \end{minipage}
&  \begin{minipage}{1 cm}
\centering
\smallskip
$1$ \end{minipage}
\\
\begin{minipage}{0.7 cm}
\centering
\smallskip
$ \mathcal{E}_9^4$  \end{minipage}
&   \begin{minipage}{8 cm}
\centering
\smallskip
$ T^3 + Y^3 + X^3 $ \end{minipage}
&  \begin{minipage}{0.7 cm}
\centering
\smallskip
$9$ \end{minipage}
&  \begin{minipage}{1.3 cm}
\centering
\smallskip
$108$ \end{minipage}
&  \begin{minipage}{1 cm}
\centering
\smallskip
$0$ \end{minipage}
\\
\begin{minipage}{0.7 cm}
\centering
\smallskip
$ \mathcal{G}_3^1$  \end{minipage}
&   \begin{minipage}{8 cm}
\centering
\smallskip
$ \xi^{9}T^3 + \xi^{9}YT^2 + \xi^{9}Y^2T + \xi^{9}XT^2 + \xi^{9}Y^3 + XYT + \xi^{9}XY^2 + \xi^{9}X^2T + \xi^{9}X^2Y + \xi^{9}X^3 $ \end{minipage}
&  \begin{minipage}{0.7 cm}
\centering
\smallskip
$18$ \end{minipage}
&  \begin{minipage}{1.3 cm}
\centering
\smallskip
$4$ \end{minipage}
&  \begin{minipage}{1 cm}
\centering
\smallskip
$\xi^{11}$ \end{minipage}
\\
\begin{minipage}{0.7 cm}
\centering
\smallskip
$ \mathcal{G}_3^1$  \end{minipage}
&   \begin{minipage}{8 cm}
\centering
\smallskip
$ \xi^{6}T^3 + \xi^{6}YT^2 + \xi^{6}Y^2T + \xi^{6}XT^2 + \xi^{6}Y^3 + XYT + \xi^{6}XY^2 + \xi^{6}X^2T + \xi^{6}X^2Y + \xi^{6}X^3 $ \end{minipage}
&  \begin{minipage}{0.7 cm}
\centering
\smallskip
$18$ \end{minipage}
&  \begin{minipage}{1.3 cm}
\centering
\smallskip
$4$ \end{minipage}
&  \begin{minipage}{1 cm}
\centering
\smallskip
$\xi^{14}$ \end{minipage}
\\
\hline
\end{tabular}
\end{minipage}\newline\newline\newline\newline\noindent
\begin{minipage}{\textwidth}
\centering
\textbf{Table 33} Continue $1$  \medskip \\
\medskip
\begin{tabular}
{|c|c|c|c|c|}
\hline
type & equation &$n$ &
\begin{minipage}{1.3 cm}
\smallskip
\centering
$\#$ of \\
residual\\
points
\end{minipage}
& $j-$inv. \\
\hline
\begin{minipage}{0.7 cm}
\centering
\smallskip
$ \mathcal{G}_0^4$  \end{minipage}
&   \begin{minipage}{8 cm}
\centering
\smallskip
$ \xi^{2}T^3 + \xi Y^3 + XYT + X^3 $ \end{minipage}
&  \begin{minipage}{0.7 cm}
\centering
\smallskip
$18$ \end{minipage}
&  \begin{minipage}{1.3 cm}
\centering
\smallskip
$9$ \end{minipage}
&  \begin{minipage}{1 cm}
\centering
\smallskip
$1$ \end{minipage}
\\
\begin{minipage}{0.7 cm}
\centering
\smallskip
$ \mathcal{G}_0^4$  \end{minipage}
&   \begin{minipage}{8 cm}
\centering
\smallskip
$ \xi^{2}T^3 + \xi Y^3 + \xi^{7}XYT + X^3 $ \end{minipage}
&  \begin{minipage}{0.7 cm}
\centering
\smallskip
$18$ \end{minipage}
&  \begin{minipage}{1.3 cm}
\centering
\smallskip
$9$ \end{minipage}
&  \begin{minipage}{1 cm}
\centering
\smallskip
$1$ \end{minipage}
\\
\begin{minipage}{0.7 cm}
\centering
\smallskip
$ \mathcal{G}_0^4$  \end{minipage}
&   \begin{minipage}{8 cm}
\centering
\smallskip
$ \xi^{2}T^3 + \xi Y^3 + \xi^{14}XYT + X^3 $ \end{minipage}
&  \begin{minipage}{0.7 cm}
\centering
\smallskip
$18$ \end{minipage}
&  \begin{minipage}{1.3 cm}
\centering
\smallskip
$9$ \end{minipage}
&  \begin{minipage}{1 cm}
\centering
\smallskip
$1$ \end{minipage}
\\
\begin{minipage}{0.7 cm}
\centering
\smallskip
$ \mathcal{G}_0^4$  \end{minipage}
&   \begin{minipage}{8 cm}
\centering
\smallskip
$ \xi^{2}T^3 + \xi Y^3 + \xi^{13}XYT + X^3 $ \end{minipage}
&  \begin{minipage}{0.7 cm}
\centering
\smallskip
$18$ \end{minipage}
&  \begin{minipage}{1.3 cm}
\centering
\smallskip
$9$ \end{minipage}
&  \begin{minipage}{1 cm}
\centering
\smallskip
$1$ \end{minipage}
\\
\hline
\end{tabular}
\end{minipage}\newline\phantom  \newline\newline\newline\newline


\noindent\begin{minipage}{\textwidth}
\centering
\textbf{Table 34}  $ q = 16$, non-singular cubic curves giving complete $(n,3)$-arcs  \medskip \\
\medskip
\begin{tabular}
{|c|c|c|c|}
\hline
type & equation &$n$ &
\begin{minipage}{1 cm}
\smallskip
\centering
$j-$ \\
inv. \\
\smallskip
\end{minipage}  \\
\hline
\begin{minipage}{0.7 cm}
\centering
\smallskip
$ \mathcal{G}_3^1$  \end{minipage}
&   \begin{minipage}{10 cm}
\centering
\smallskip
$ \xi^{5}T^3 + \xi^{5}YT^2 + \xi^{5}Y^2T + \xi^{5}XT^2 + \xi^{5}Y^3 + XYT + \xi^{5}XY^2 + \xi^{5}X^2T + \xi^{5}X^2Y + \xi^{5}X^3 $ \end{minipage}
&  \begin{minipage}{0.7 cm}
\centering
\smallskip
$24$ \end{minipage}
&  \begin{minipage}{1 cm}
\centering
\smallskip
$\xi^{5}$ \end{minipage}
\\
\begin{minipage}{0.7 cm}
\centering
\smallskip
$ \mathcal{G}_3^1$  \end{minipage}
&   \begin{minipage}{10 cm}
\centering
\smallskip
$ \xi^{10}T^3 + \xi^{10}YT^2 + \xi^{10}Y^2T + \xi^{10}XT^2 + \xi^{10}Y^3 + XYT + \xi^{10}XY^2 + \xi^{10}X^2T + \xi^{10}X^2Y + \xi^{10}X^3 $ \end{minipage}
&  \begin{minipage}{0.7 cm}
\centering
\smallskip
$24$ \end{minipage}
&  \begin{minipage}{1 cm}
\centering
\smallskip
$\xi^{10}$ \end{minipage}
\\
\begin{minipage}{0.7 cm}
\centering
\smallskip
$ \overline{\mathcal{E}}_3^1$  \end{minipage}
&   \begin{minipage}{10 cm}
\centering
\smallskip
$ \xi^{2}T^3 + XY^2 + X^2Y $ \end{minipage}
&  \begin{minipage}{0.7 cm}
\centering
\smallskip
$21$ \end{minipage}
&  \begin{minipage}{1 cm}
\centering
\smallskip
$0$ \end{minipage}
\\
\begin{minipage}{0.7 cm}
\centering
\smallskip
$ \overline{\mathcal{E}}_3^1$  \end{minipage}
&   \begin{minipage}{10 cm}
\centering
\smallskip
$ \xi T^3 + XY^2 + X^2Y $ \end{minipage}
&  \begin{minipage}{0.7 cm}
\centering
\smallskip
$21$ \end{minipage}
&  \begin{minipage}{1 cm}
\centering
\smallskip
$0$ \end{minipage}
\\
\begin{minipage}{0.7 cm}
\centering
\smallskip
$ \mathcal{E}_1^4$  \end{minipage}
&   \begin{minipage}{10 cm}
\centering
\smallskip
$ YT^2 + Y^2T + XY^2 + X^3 $ \end{minipage}
&  \begin{minipage}{0.7 cm}
\centering
\smallskip
$25$ \end{minipage}
&  \begin{minipage}{1 cm}
\centering
\smallskip
$0$ \end{minipage}
\\
\hline
\end{tabular}
\end{minipage}\newline


$\; $ \phantom  \newline\newline\newline\newline\noindent
\begin{minipage}{\textwidth}
\centering
\textbf{Table 35} $ q = 17$, a. i. singular cubic curves giving incomplete $(n,3)$-arcs  \medskip \\
\medskip
\begin{tabular}
{|c|c|c|c|}
\hline
type & equation &$n$ &
\begin{minipage}{1.3 cm}
\smallskip
\centering
$\#$ of \\
residual\\
points
\end{minipage}\\
\hline
\begin{minipage}{0.7 cm}
\centering
\smallskip
$ \mathcal{N}_3^0$  \end{minipage}
&   \begin{minipage}{8 cm}
\centering
\smallskip
$ 8YT^2 + 3XT^2 + Y^3 + XY^2 $ \end{minipage}
&  \begin{minipage}{0.7 cm}
\centering
\smallskip
$19$ \end{minipage}
&  \begin{minipage}{1.3 cm}
\centering
\smallskip
$6$ \end{minipage}
\\
\begin{minipage}{0.7 cm}
\centering
\smallskip
$ \mathcal{N}_1^2$  \end{minipage}
&   \begin{minipage}{6 cm}
\centering
\smallskip
$ T^3 + Y^3 + XYT $ \end{minipage}
&  \begin{minipage}{0.7 cm}
\centering
\smallskip
$17$ \end{minipage}
&  \begin{minipage}{1.3 cm}
\centering
\smallskip
$8$ \end{minipage}
\\
\begin{minipage}{0.7 cm}
\centering
\smallskip
$ \mathcal{N}_1^1$  \end{minipage}
&   \begin{minipage}{8 cm}
\centering
\smallskip
$ T^3 + XY^2 $ \end{minipage}
&  \begin{minipage}{0.7 cm}
\centering
\smallskip
$18$ \end{minipage}
&  \begin{minipage}{1.3 cm}
\centering
\smallskip
$17$ \end{minipage}
\\
\begin{minipage}{0.7 cm}
\centering
\smallskip
$ \mathcal{N}_0^0$  \end{minipage}
&   \begin{minipage}{8 cm}
\centering
\smallskip
$ 16T^3 + 8YT^2 + Y^2T + 3XT^2 + Y^3 + XY^2 $ \end{minipage}
&  \begin{minipage}{0.7 cm}
\centering
\smallskip
$19$ \end{minipage}
&  \begin{minipage}{1.3 cm}
\centering
\smallskip
$6$ \end{minipage}
\\
\hline
\end{tabular}
\end{minipage}\newline\phantom  \newline\newline\newline\newline


\noindent\begin{minipage}{\textwidth}
\centering
\textbf{Table 36}  $ q = 17$, non-singular cubic curves giving incomplete $(n,3)$-arcs  \medskip \\
\medskip
\begin{tabular}
{|c|c|c|c|c|}
\hline
type & equation &$n$ &
\begin{minipage}{1.3 cm}
\smallskip
\centering
$\#$ of \\
residual\\
points
\end{minipage}   &
\begin{minipage}{1 cm}
\smallskip
\centering
$j-$ \\
inv. \\
\end{minipage}  \\
\hline
\begin{minipage}{0.7 cm}
\centering
\smallskip
$ \mathcal{G}_3^2$  \end{minipage}
&   \begin{minipage}{8 cm}
\centering
\smallskip
$ 4T^3 + 12YT^2 + 12Y^2T + 12XT^2 + 4Y^3 + 8XYT + 12XY^2 + 12X^2T + 12X^2Y + 4X^3 $ \end{minipage}
&  \begin{minipage}{0.7 cm}
\centering
\smallskip
$24$ \end{minipage}
&  \begin{minipage}{1.3 cm}
\centering
\smallskip
$1$ \end{minipage}
&  \begin{minipage}{1 cm}
\centering
\smallskip
$1$ \end{minipage}
\\
\begin{minipage}{0.7 cm}
\centering
\smallskip
$ \mathcal{G}_3^2$  \end{minipage}
&   \begin{minipage}{8 cm}
\centering
\smallskip
$ 10T^3 + 13YT^2 + 13Y^2T + 13XT^2 + 10Y^3 + 10XYT + 13XY^2 + 13X^2T + 13X^2Y + 10X^3 $ \end{minipage}
&  \begin{minipage}{0.7 cm}
\centering
\smallskip
$18$ \end{minipage}
&  \begin{minipage}{1.3 cm}
\centering
\smallskip
$6$ \end{minipage}
&  \begin{minipage}{1 cm}
\centering
\smallskip
$8$ \end{minipage}
\\
\begin{minipage}{0.7 cm}
\centering
\smallskip
$ \mathcal{G}_3^2$  \end{minipage}
&   \begin{minipage}{8 cm}
\centering
\smallskip
$ 6T^3 + YT^2 + Y^2T + XT^2 + 6Y^3 + 3XYT + XY^2 + X^2T + X^2Y + 6X^3 $ \end{minipage}
&  \begin{minipage}{0.7 cm}
\centering
\smallskip
$15$ \end{minipage}
&  \begin{minipage}{1.3 cm}
\centering
\smallskip
$25$ \end{minipage}
&  \begin{minipage}{1 cm}
\centering
\smallskip
$2$ \end{minipage}
\\
\begin{minipage}{0.7 cm}
\centering
\smallskip
$ \mathcal{G}_3^2$  \end{minipage}
&   \begin{minipage}{8 cm}
\centering
\smallskip
$ 3T^3 + 9YT^2 + 9Y^2T + 9XT^2 + 3Y^3 + 2XYT + 9XY^2 + 9X^2T + 9X^2Y + 3X^3 $ \end{minipage}
&  \begin{minipage}{0.7 cm}
\centering
\smallskip
$15$ \end{minipage}
&  \begin{minipage}{1.3 cm}
\centering
\smallskip
$25$ \end{minipage}
&  \begin{minipage}{1 cm}
\centering
\smallskip
$7$ \end{minipage}
\\
\hline
\end{tabular}
\end{minipage}\newline\newline\newline\newline\noindent
\begin{minipage}{\textwidth}
\centering
\textbf{Table 36} Continue $1$  \medskip \\
\medskip
\begin{tabular}
{|c|c|c|c|c|}
\hline
type & equation &$n$ &
\begin{minipage}{1.3 cm}
\smallskip
\centering
$\#$ of \\
residual\\
points
\end{minipage}
& $j-$inv. \\
\hline
\begin{minipage}{0.7 cm}
\centering
\smallskip
$ \mathcal{G}_1$  \end{minipage}
&   \begin{minipage}{8 cm}
\centering
\smallskip
$ YT^2 + 6Y^3 + XY^2 + X^3 $ \end{minipage}
&  \begin{minipage}{0.7 cm}
\centering
\smallskip
$20$ \end{minipage}
&  \begin{minipage}{1.3 cm}
\centering
\smallskip
$1$ \end{minipage}
&  \begin{minipage}{1 cm}
\centering
\smallskip
$16$ \end{minipage}
\\
\begin{minipage}{0.7 cm}
\centering
\smallskip
$ \mathcal{G}_1$  \end{minipage}
&   \begin{minipage}{8 cm}
\centering
\smallskip
$ YT^2 + 2Y^3 + 2XY^2 + X^3 $ \end{minipage}
&  \begin{minipage}{0.7 cm}
\centering
\smallskip
$19$ \end{minipage}
&  \begin{minipage}{1.3 cm}
\centering
\smallskip
$2$ \end{minipage}
&  \begin{minipage}{1 cm}
\centering
\smallskip
$3$ \end{minipage}
\\
\begin{minipage}{0.7 cm}
\centering
\smallskip
$ \mathcal{G}_1$  \end{minipage}
&   \begin{minipage}{8 cm}
\centering
\smallskip
$ YT^2 + 3Y^3 + XY^2 + X^3 $ \end{minipage}
&  \begin{minipage}{0.7 cm}
\centering
\smallskip
$17$ \end{minipage}
&  \begin{minipage}{1.3 cm}
\centering
\smallskip
$3$ \end{minipage}
&  \begin{minipage}{1 cm}
\centering
\smallskip
$3$ \end{minipage}
\\
\begin{minipage}{0.7 cm}
\centering
\smallskip
$ \mathcal{G}_1$  \end{minipage}
&   \begin{minipage}{8 cm}
\centering
\smallskip
$ YT^2 + 4Y^3 + 2XY^2 + X^3 $ \end{minipage}
&  \begin{minipage}{0.7 cm}
\centering
\smallskip
$16$ \end{minipage}
&  \begin{minipage}{1.3 cm}
\centering
\smallskip
$13$ \end{minipage}
&  \begin{minipage}{1 cm}
\centering
\smallskip
$16$ \end{minipage}
\\
\hline
\end{tabular}
\end{minipage}\newline\phantom  \newline\newline\newline\newline


\noindent\begin{minipage}{\textwidth}
\centering
\textbf{Table 37}  $ q = 17$, non-singular cubic curves giving complete $(n,3)$-arcs  \medskip \\
\medskip
\begin{tabular}
{|c|c|c|c|}
\hline
type & equation &$n$ &
\begin{minipage}{1 cm}
\smallskip
\centering
$j-$ \\
inv. \\
\smallskip
\end{minipage}  \\
\hline
\begin{minipage}{0.7 cm}
\centering
\smallskip
$ \mathcal{G}_3^2$  \end{minipage}
&   \begin{minipage}{10 cm}
\centering
\smallskip
$ 9T^3 + 10YT^2 + 10Y^2T + 10XT^2 + 9Y^3 + 4XYT + 10XY^2 + 10X^2T + 10X^2Y + 9X^3 $ \end{minipage}
&  \begin{minipage}{0.7 cm}
\centering
\smallskip
$24$ \end{minipage}
&  \begin{minipage}{1 cm}
\centering
\smallskip
$6$ \end{minipage}
\\
\begin{minipage}{0.7 cm}
\centering
\smallskip
$ \mathcal{G}_3^2$  \end{minipage}
&   \begin{minipage}{10 cm}
\centering
\smallskip
$ 8T^3 + 7YT^2 + 7Y^2T + 7XT^2 + 8Y^3 + 15XYT + 7XY^2 + 7X^2T + 7X^2Y + 8X^3 $ \end{minipage}
&  \begin{minipage}{0.7 cm}
\centering
\smallskip
$24$ \end{minipage}
&  \begin{minipage}{1 cm}
\centering
\smallskip
$10$ \end{minipage}
\\
\begin{minipage}{0.7 cm}
\centering
\smallskip
$ \mathcal{G}_3^2$  \end{minipage}
&   \begin{minipage}{10 cm}
\centering
\smallskip
$ 7T^3 + 4YT^2 + 4Y^2T + 4XT^2 + 7Y^3 + 9XYT + 4XY^2 + 4X^2T + 4X^2Y + 7X^3 $ \end{minipage}
&  \begin{minipage}{0.7 cm}
\centering
\smallskip
$21$ \end{minipage}
&  \begin{minipage}{1 cm}
\centering
\smallskip
$2$ \end{minipage}
\\
\begin{minipage}{0.7 cm}
\centering
\smallskip
$ \mathcal{G}_3^2$  \end{minipage}
&   \begin{minipage}{10 cm}
\centering
\smallskip
$ 2T^3 + 6YT^2 + 6Y^2T + 6XT^2 + 2Y^3 + 13XYT + 6XY^2 + 6X^2T + 6X^2Y + 2X^3 $ \end{minipage}
&  \begin{minipage}{0.7 cm}
\centering
\smallskip
$21$ \end{minipage}
&  \begin{minipage}{1 cm}
\centering
\smallskip
$13$ \end{minipage}
\\
\begin{minipage}{0.7 cm}
\centering
\smallskip
$ \mathcal{G}_3^2$  \end{minipage}
&   \begin{minipage}{10 cm}
\centering
\smallskip
$ 16T^3 + 14YT^2 + 14Y^2T + 14XT^2 + 16Y^3 + 12XYT + 14XY^2 + 14X^2T + 14X^2Y + 16X^3 $ \end{minipage}
&  \begin{minipage}{0.7 cm}
\centering
\smallskip
$21$ \end{minipage}
&  \begin{minipage}{1 cm}
\centering
\smallskip
$7$ \end{minipage}
\\
\hline
\end{tabular}
\end{minipage}\newline


$\; $ \phantom  \newline\newline\newline\newline\noindent
\begin{minipage}{\textwidth}
\centering
\textbf{Table 38} $ q = 19$, a. i. singular cubic curves giving incomplete $(n,3)$-arcs  \medskip \\
\medskip
\begin{tabular}
{|c|c|c|c|}
\hline
type & equation &$n$ &
\begin{minipage}{1.3 cm}
\smallskip
\centering
$\#$ of \\
residual\\
points
\end{minipage}\\
\hline
\begin{minipage}{0.7 cm}
\centering
\smallskip
$ \mathcal{N}_3^2$  \end{minipage}
&   \begin{minipage}{6 cm}
\centering
\smallskip
$ T^3 + Y^3 + XYT $ \end{minipage}
&  \begin{minipage}{0.7 cm}
\centering
\smallskip
$19$ \end{minipage}
&  \begin{minipage}{1.3 cm}
\centering
\smallskip
$3$ \end{minipage}
\\
\begin{minipage}{0.7 cm}
\centering
\smallskip
$ \mathcal{N}_1^1$  \end{minipage}
&   \begin{minipage}{6 cm}
\centering
\smallskip
$ T^3 + XY^2 $ \end{minipage}
&  \begin{minipage}{0.7 cm}
\centering
\smallskip
$20$ \end{minipage}
&  \begin{minipage}{1.3 cm}
\centering
\smallskip
$6$ \end{minipage}
\\
\begin{minipage}{0.7 cm}
\centering
\smallskip
$ \mathcal{N}_1^0$  \end{minipage}
&   \begin{minipage}{6 cm}
\centering
\smallskip
$ 6YT^2 + 17XT^2 + Y^3 + XY^2 $ \end{minipage}
&  \begin{minipage}{0.7 cm}
\centering
\smallskip
$21$ \end{minipage}
&  \begin{minipage}{1.3 cm}
\centering
\smallskip
$3$ \end{minipage}
\\
\begin{minipage}{0.7 cm}
\centering
\smallskip
$ \mathcal{N}_0^2$  \end{minipage}
&   \begin{minipage}{6 cm}
\centering
\smallskip
$ 2T^3 + Y^3 + XYT $ \end{minipage}
&  \begin{minipage}{0.7 cm}
\centering
\smallskip
$19$ \end{minipage}
&  \begin{minipage}{1.3 cm}
\centering
\smallskip
$4$ \end{minipage}
\\
\hline
\end{tabular}
\end{minipage}\newline\phantom  \newline\newline\newline\newline


\noindent\begin{minipage}{\textwidth}
\centering
\textbf{Table 39}  $ q = 19$, non-singular cubic curves giving incomplete $(n,3)$-arcs  \medskip \\
\medskip
\begin{tabular}
{|c|c|c|c|c|}
\hline
type & equation &$n$ &
\begin{minipage}{1.3 cm}
\smallskip
\centering
$\#$ of \\
residual\\
points
\end{minipage}   &
\begin{minipage}{1 cm}
\smallskip
\centering
$j-$ \\
inv. \\
\end{minipage}  \\
\hline
\begin{minipage}{0.7 cm}
\centering
\smallskip
$ \mathcal{G}_9^4$  \end{minipage}
&   \begin{minipage}{8 cm}
\centering
\smallskip
$ T^3 + Y^3 + 13XYT + X^3 $ \end{minipage}
&  \begin{minipage}{0.7 cm}
\centering
\smallskip
$18$ \end{minipage}
&  \begin{minipage}{1.3 cm}
\centering
\smallskip
$9$ \end{minipage}
&  \begin{minipage}{1 cm}
\centering
\smallskip
$1$ \end{minipage}
\\
\begin{minipage}{0.7 cm}
\centering
\smallskip
$ \mathcal{G}_3^1$  \end{minipage}
&   \begin{minipage}{8 cm}
\centering
\smallskip
$ T^3 + 3YT^2 + 3Y^2T + 3XT^2 + Y^3 + 7XYT + 3XY^2 + 3X^2T + 3X^2Y + X^3 $ \end{minipage}
&  \begin{minipage}{0.7 cm}
\centering
\smallskip
$18$ \end{minipage}
&  \begin{minipage}{1.3 cm}
\centering
\smallskip
$6$ \end{minipage}
&  \begin{minipage}{1 cm}
\centering
\smallskip
$14$ \end{minipage}
\\
\begin{minipage}{0.7 cm}
\centering
\smallskip
$ \mathcal{G}_3^1$  \end{minipage}
&   \begin{minipage}{8 cm}
\centering
\smallskip
$ 12T^3 + 17YT^2 + 17Y^2T + 17XT^2 + 12Y^3 + 16XYT + 17XY^2 + 17X^2T + 17X^2Y + 12X^3 $ \end{minipage}
&  \begin{minipage}{0.7 cm}
\centering
\smallskip
$18$ \end{minipage}
&  \begin{minipage}{1.3 cm}
\centering
\smallskip
$7$ \end{minipage}
&  \begin{minipage}{1 cm}
\centering
\smallskip
$13$ \end{minipage}
\\
\begin{minipage}{0.7 cm}
\centering
\smallskip
$ \mathcal{G}_3^1$  \end{minipage}
&   \begin{minipage}{8 cm}
\centering
\smallskip
$ 4T^3 + 12YT^2 + 12Y^2T + 12XT^2 + 4Y^3 + 6XYT + 12XY^2 + 12X^2T + 12X^2Y + 4X^3 $ \end{minipage}
&  \begin{minipage}{0.7 cm}
\centering
\smallskip
$15$ \end{minipage}
&  \begin{minipage}{1.3 cm}
\centering
\smallskip
$28$ \end{minipage}
&  \begin{minipage}{1 cm}
\centering
\smallskip
$10$ \end{minipage}
\\
\hline
\end{tabular}
\end{minipage}\newline\newline\newline\newline\noindent
\begin{minipage}{\textwidth}
\centering
\textbf{Table 39} Continue $1$  \medskip \\
\medskip
\begin{tabular}
{|c|c|c|c|c|}
\hline
type & equation &$n$ &
\begin{minipage}{1.3 cm}
\smallskip
\centering
$\#$ of \\
residual\\
points
\end{minipage}
& $j-$inv. \\
\hline
\begin{minipage}{0.7 cm}
\centering
\smallskip
$ \mathcal{G}_1$  \end{minipage}
&   \begin{minipage}{8 cm}
\centering
\smallskip
$ YT^2 + 7Y^3 + XY^2 + X^3 $ \end{minipage}
&  \begin{minipage}{0.7 cm}
\centering
\smallskip
$22$ \end{minipage}
&  \begin{minipage}{1.3 cm}
\centering
\smallskip
$1$ \end{minipage}
&  \begin{minipage}{1 cm}
\centering
\smallskip
$14$ \end{minipage}
\\
\begin{minipage}{0.7 cm}
\centering
\smallskip
$ \mathcal{G}_1$  \end{minipage}
&   \begin{minipage}{8 cm}
\centering
\smallskip
$ YT^2 + 6Y^3 + XY^2 + X^3 $ \end{minipage}
&  \begin{minipage}{0.7 cm}
\centering
\smallskip
$22$ \end{minipage}
&  \begin{minipage}{1.3 cm}
\centering
\smallskip
$1$ \end{minipage}
&  \begin{minipage}{1 cm}
\centering
\smallskip
$13$ \end{minipage}
\\
\begin{minipage}{0.7 cm}
\centering
\smallskip
$ \mathcal{G}_1$  \end{minipage}
&   \begin{minipage}{8 cm}
\centering
\smallskip
$ YT^2 + 3Y^3 + 2XY^2 + X^3 $ \end{minipage}
&  \begin{minipage}{0.7 cm}
\centering
\smallskip
$20$ \end{minipage}
&  \begin{minipage}{1.3 cm}
\centering
\smallskip
$3$ \end{minipage}
&  \begin{minipage}{1 cm}
\centering
\smallskip
$7$ \end{minipage}
\\
\begin{minipage}{0.7 cm}
\centering
\smallskip
$ \mathcal{G}_1$  \end{minipage}
&   \begin{minipage}{8 cm}
\centering
\smallskip
$ YT^2 + 16Y^3 + 2XY^2 + X^3 $ \end{minipage}
&  \begin{minipage}{0.7 cm}
\centering
\smallskip
$20$ \end{minipage}
&  \begin{minipage}{1.3 cm}
\centering
\smallskip
$2$ \end{minipage}
&  \begin{minipage}{1 cm}
\centering
\smallskip
$7$ \end{minipage}
\\
\hline
\end{tabular}
\end{minipage}\newline\phantom  \newline\newline


\noindent\begin{minipage}{\textwidth}
\centering
\textbf{Table 40}  $ q = 19$, non-singular cubic curves giving complete $(n,3)$-arcs  \medskip \\
\medskip
\begin{tabular}
{|c|c|c|c|}
\hline
type & equation &$n$ &
\begin{minipage}{1 cm}
\smallskip
\centering
$j-$ \\
inv. \\
\smallskip
\end{minipage}  \\
\hline
\begin{minipage}{0.7 cm}
\centering
\smallskip
$ \mathcal{E}_9^4$  \end{minipage}
&   \begin{minipage}{10 cm}
\centering
\smallskip
$ T^3 + Y^3 + X^3 $ \end{minipage}
&  \begin{minipage}{0.7 cm}
\centering
\smallskip
$27$ \end{minipage}
&  \begin{minipage}{1 cm}
\centering
\smallskip
$0$ \end{minipage}
\\
\begin{minipage}{0.7 cm}
\centering
\smallskip
$ \mathcal{G}_3^1$  \end{minipage}
&   \begin{minipage}{10 cm}
\centering
\smallskip
$ 11T^3 + 14YT^2 + 14Y^2T + 14XT^2 + 11Y^3 + 10XYT + 14XY^2 + 14X^2T + 14X^2Y + 11X^3 $ \end{minipage}
&  \begin{minipage}{0.7 cm}
\centering
\smallskip
$27$ \end{minipage}
&  \begin{minipage}{1 cm}
\centering
\smallskip
$3$ \end{minipage}
\\
\begin{minipage}{0.7 cm}
\centering
\smallskip
$ \mathcal{G}_3^1$  \end{minipage}
&   \begin{minipage}{10 cm}
\centering
\smallskip
$ 9T^3 + 8YT^2 + 8Y^2T + 8XT^2 + 9Y^3 + 17XYT + 8XY^2 + 8X^2T + 8X^2Y + 9X^3 $ \end{minipage}
&  \begin{minipage}{0.7 cm}
\centering
\smallskip
$24$ \end{minipage}
&  \begin{minipage}{1 cm}
\centering
\smallskip
$9$ \end{minipage}
\\
\begin{minipage}{0.7 cm}
\centering
\smallskip
$ \mathcal{G}_3^1$  \end{minipage}
&   \begin{minipage}{10 cm}
\centering
\smallskip
$ 17T^3 + 13YT^2 + 13Y^2T + 13XT^2 + 17Y^3 + 8XYT + 13XY^2 + 13X^2T + 13X^2Y + 17X^3 $ \end{minipage}
&  \begin{minipage}{0.7 cm}
\centering
\smallskip
$24$ \end{minipage}
&  \begin{minipage}{1 cm}
\centering
\smallskip
$17$ \end{minipage}
\\
\begin{minipage}{0.7 cm}
\centering
\smallskip
$ \mathcal{G}_3^1$  \end{minipage}
&   \begin{minipage}{10 cm}
\centering
\smallskip
$ 16T^3 + 10YT^2 + 10Y^2T + 10XT^2 + 16Y^3 + 2XYT + 10XY^2 + 10X^2T + 10X^2Y + 16X^3 $ \end{minipage}
&  \begin{minipage}{0.7 cm}
\centering
\smallskip
$24$ \end{minipage}
&  \begin{minipage}{1 cm}
\centering
\smallskip
$5$ \end{minipage}
\\
\hline
\end{tabular}
\end{minipage}\newline\newline\newline\newline\noindent
\begin{minipage}{\textwidth}
\centering
\textbf{Table 40} Continue $1$  \medskip \\
\medskip
\begin{tabular}
{|c|c|c|c|}
\hline
type & equation & n  &
\begin{minipage}{1 cm}
\smallskip
\centering
$j-$ \\
inv. \\
\smallskip
\end{minipage}  \\
\hline
\begin{minipage}{0.7 cm}
\centering
\smallskip
$ \mathcal{G}_0^1$  \end{minipage}
&   \begin{minipage}{10 cm}
\centering
\smallskip
$ 12T^3 + 4YT^2 + 3Y^3 + 10XYT + XY^2 + X^2T + 15X^3 $ \end{minipage}
&  \begin{minipage}{0.7 cm}
\centering
\smallskip
$27$ \end{minipage}
&  \begin{minipage}{1 cm}
\centering
\smallskip
$3$ \end{minipage}
\\
\begin{minipage}{0.7 cm}
\centering
\smallskip
$ \mathcal{G}_0^1$  \end{minipage}
&   \begin{minipage}{10 cm}
\centering
\smallskip
$ 6T^3 + 2YT^2 + 3Y^3 + 10XYT + XY^2 + X^2T + 11X^3 $ \end{minipage}
&  \begin{minipage}{0.7 cm}
\centering
\smallskip
$24$ \end{minipage}
&  \begin{minipage}{1 cm}
\centering
\smallskip
$9$ \end{minipage}
\\
\begin{minipage}{0.7 cm}
\centering
\smallskip
$ \mathcal{G}_0^1$  \end{minipage}
&   \begin{minipage}{10 cm}
\centering
\smallskip
$ 5T^3 + 4YT^2 + 6Y^3 + XYT + XY^2 + X^2T + 11X^3 $ \end{minipage}
&  \begin{minipage}{0.7 cm}
\centering
\smallskip
$24$ \end{minipage}
&  \begin{minipage}{1 cm}
\centering
\smallskip
$17$ \end{minipage}
\\
\begin{minipage}{0.7 cm}
\centering
\smallskip
$ \mathcal{G}_0^1$  \end{minipage}
&   \begin{minipage}{10 cm}
\centering
\smallskip
$ 3T^3 + 4YT^2 + 15Y^3 + 12XYT + XY^2 + X^2T + 18X^3 $ \end{minipage}
&  \begin{minipage}{0.7 cm}
\centering
\smallskip
$24$ \end{minipage}
&  \begin{minipage}{1 cm}
\centering
\smallskip
$5$ \end{minipage}
\\
\begin{minipage}{0.7 cm}
\centering
\smallskip
$ \mathcal{G}_0^1$  \end{minipage}
&   \begin{minipage}{10 cm}
\centering
\smallskip
$ 17T^3 + 2YT^2 + 18Y^3 + 3XYT + XY^2 + X^2T + 9X^3 $ \end{minipage}
&  \begin{minipage}{0.7 cm}
\centering
\smallskip
$24$ \end{minipage}
&  \begin{minipage}{1 cm}
\centering
\smallskip
$15$ \end{minipage}
\\
\hline
\end{tabular}
\end{minipage}\newline


$\; $ \phantom  \newline\newline\newline\newline\noindent
\begin{minipage}{\textwidth}
\centering
\textbf{Table 41} $ q = 23$, a. i. singular cubic curves giving incomplete $(n,3)$-arcs  \medskip \\
\medskip
\begin{tabular}
{|c|c|c|c|}
\hline
type & equation &$n$ &
\begin{minipage}{1.3 cm}
\smallskip
\centering
$\#$ of \\
residual\\
points
\end{minipage}\\
\hline
\begin{minipage}{0.7 cm}
\centering
\smallskip
$ \mathcal{N}_3^0$  \end{minipage}
&   \begin{minipage}{6 cm}
\centering
\smallskip
$ 14YT^2 + 3XT^2 + Y^3 + XY^2 $ \end{minipage}
&  \begin{minipage}{0.7 cm}
\centering
\smallskip
$25$ \end{minipage}
&  \begin{minipage}{1.3 cm}
\centering
\smallskip
$3$ \end{minipage}
\\
\begin{minipage}{0.7 cm}
\centering
\smallskip
$ \mathcal{N}_1^1$  \end{minipage}
&   \begin{minipage}{6 cm}
\centering
\smallskip
$ T^3 + XY^2 $ \end{minipage}
&  \begin{minipage}{0.7 cm}
\centering
\smallskip
$24$ \end{minipage}
&  \begin{minipage}{1.3 cm}
\centering
\smallskip
$1$ \end{minipage}
\\
\hline
\end{tabular}
\end{minipage}\newline


$\; $ \phantom  \newline\newline\newline\newline\noindent
\begin{minipage}{\textwidth}
\centering
\textbf{Table 42} $ q = 23$, a. i. singular cubic curves giving complete $(n,3)$-arcs  \medskip \\
\medskip
\begin{tabular}
{|c|c|c|}
\hline
type & equation &$n$ \\
\hline
\begin{minipage}{0.7 cm}
\centering
\smallskip
$ \mathcal{N}_1^2$  \end{minipage}
&   \begin{minipage}{8 cm}
\centering
\smallskip
$ T^3 + Y^3 + XYT $ \end{minipage}
&  \begin{minipage}{0.7 cm}
\centering
\smallskip
$23$ \end{minipage}
\\
\begin{minipage}{0.7 cm}
\centering
\smallskip
$ \mathcal{N}_0^0$  \end{minipage}
&   \begin{minipage}{8 cm}
\centering
\smallskip
$ 22T^3 + 14YT^2 + Y^2T + 3XT^2 + Y^3 + XY^2 $ \end{minipage}
&  \begin{minipage}{0.7 cm}
\centering
\smallskip
$25$ \end{minipage}
\\
\hline
\end{tabular}
\end{minipage}\newline\phantom  \newline\newline\newline\newline


\noindent\begin{minipage}{\textwidth}
\centering
\textbf{Table 43}  $ q = 23$, non-singular cubic curves giving incomplete $(n,3)$-arcs  \medskip \\
\medskip
\begin{tabular}
{|c|c|c|c|c|}
\hline
type & equation &$n$ &
\begin{minipage}{1.3 cm}
\smallskip
\centering
$\#$ of \\
residual\\
points
\end{minipage}   &
\begin{minipage}{1 cm}
\smallskip
\centering
$j-$ \\
inv. \\
\end{minipage}  \\
\hline
\begin{minipage}{0.7 cm}
\centering
\smallskip
$ \mathcal{G}_3^2$  \end{minipage}
&   \begin{minipage}{8 cm}
\centering
\smallskip
$ 9T^3 + 4YT^2 + 4Y^2T + 4XT^2 + 9Y^3 + 9XYT + 4XY^2 + 4X^2T + 4X^2Y + 9X^3 $ \end{minipage}
&  \begin{minipage}{0.7 cm}
\centering
\smallskip
$21$ \end{minipage}
&  \begin{minipage}{1.3 cm}
\centering
\smallskip
$4$ \end{minipage}
&  \begin{minipage}{1 cm}
\centering
\smallskip
$21$ \end{minipage}
\\
\begin{minipage}{0.7 cm}
\centering
\smallskip
$ \mathcal{G}_3^2$  \end{minipage}
&   \begin{minipage}{8 cm}
\centering
\smallskip
$ 7T^3 + 21YT^2 + 21Y^2T + 21XT^2 + 7Y^3 + 20XYT + 21XY^2 + 21X^2T + 21X^2Y + 7X^3 $ \end{minipage}
&  \begin{minipage}{0.7 cm}
\centering
\smallskip
$21$ \end{minipage}
&  \begin{minipage}{1.3 cm}
\centering
\smallskip
$3$ \end{minipage}
&  \begin{minipage}{1 cm}
\centering
\smallskip
$12$ \end{minipage}
\\
\begin{minipage}{0.7 cm}
\centering
\smallskip
$ \mathcal{G}_3^2$  \end{minipage}
&   \begin{minipage}{8 cm}
\centering
\smallskip
$ 4T^3 + 12YT^2 + 12Y^2T + 12XT^2 + 4Y^3 + 2XYT + 12XY^2 + 12X^2T + 12X^2Y + 4X^3 $ \end{minipage}
&  \begin{minipage}{0.7 cm}
\centering
\smallskip
$21$ \end{minipage}
&  \begin{minipage}{1.3 cm}
\centering
\smallskip
$7$ \end{minipage}
&  \begin{minipage}{1 cm}
\centering
\smallskip
$15$ \end{minipage}
\\
\begin{minipage}{0.7 cm}
\centering
\smallskip
$ \mathcal{G}_3^2$  \end{minipage}
&   \begin{minipage}{8 cm}
\centering
\smallskip
$ 5T^3 + 15YT^2 + 15Y^2T + 15XT^2 + 5Y^3 + 8XYT + 15XY^2 + 15X^2T + 15X^2Y + 5X^3 $ \end{minipage}
&  \begin{minipage}{0.7 cm}
\centering
\smallskip
$18$ \end{minipage}
&  \begin{minipage}{1.3 cm}
\centering
\smallskip
$25$ \end{minipage}
&  \begin{minipage}{1 cm}
\centering
\smallskip
$20$ \end{minipage}
\\
\hline
\end{tabular}
\end{minipage}\newline\newline\newline\newline\noindent
\begin{minipage}{\textwidth}
\centering
\textbf{Table 43} Continue $1$  \medskip \\
\medskip
\begin{tabular}
{|c|c|c|c|c|}
\hline
type & equation &$n$ &
\begin{minipage}{1.3 cm}
\smallskip
\centering
$\#$ of \\
residual\\
points
\end{minipage}
& $j-$inv. \\
\hline
\begin{minipage}{0.7 cm}
\centering
\smallskip
$ \mathcal{G}_0^2$  \end{minipage}
&   \begin{minipage}{8 cm}
\centering
\smallskip
$ T^3 + 14Y^2T + 16Y^3 + 17XYT + 3XY^2 + 14X^2T + 22X^3 $ \end{minipage}
&  \begin{minipage}{0.7 cm}
\centering
\smallskip
$21$ \end{minipage}
&  \begin{minipage}{1.3 cm}
\centering
\smallskip
$9$ \end{minipage}
&  \begin{minipage}{1 cm}
\centering
\smallskip
$12$ \end{minipage}
\\
\begin{minipage}{0.7 cm}
\centering
\smallskip
$ \mathcal{G}_0^2$  \end{minipage}
&   \begin{minipage}{8 cm}
\centering
\smallskip
$ T^3 + 8Y^2T + 15Y^3 + 5XYT + 3XY^2 + 8X^2T + 22X^3 $ \end{minipage}
&  \begin{minipage}{0.7 cm}
\centering
\smallskip
$18$ \end{minipage}
&  \begin{minipage}{1.3 cm}
\centering
\smallskip
$33$ \end{minipage}
&  \begin{minipage}{1 cm}
\centering
\smallskip
$18$ \end{minipage}
\\
\begin{minipage}{0.7 cm}
\centering
\smallskip
$ \mathcal{G}_0^2$  \end{minipage}
&   \begin{minipage}{8 cm}
\centering
\smallskip
$ T^3 + 20Y^2T + 20Y^3 + 9XYT + 3XY^2 + 20X^2T + 22X^3 $ \end{minipage}
&  \begin{minipage}{0.7 cm}
\centering
\smallskip
$18$ \end{minipage}
&  \begin{minipage}{1.3 cm}
\centering
\smallskip
$16$ \end{minipage}
&  \begin{minipage}{1 cm}
\centering
\smallskip
$14$ \end{minipage}
\\
\begin{minipage}{0.7 cm}
\centering
\smallskip
$ \mathcal{G}_0^2$  \end{minipage}
&   \begin{minipage}{8 cm}
\centering
\smallskip
$ T^3 + 18Y^2T + 20Y^3 + 15XYT + 3XY^2 + 18X^2T + 22X^3 $ \end{minipage}
&  \begin{minipage}{0.7 cm}
\centering
\smallskip
$18$ \end{minipage}
&  \begin{minipage}{1.3 cm}
\centering
\smallskip
$27$ \end{minipage}
&  \begin{minipage}{1 cm}
\centering
\smallskip
$22$ \end{minipage}
\\
\hline
\end{tabular}
\end{minipage}\newline\phantom  \newline\newline\newline\newline


\noindent\begin{minipage}{\textwidth}
\centering
\textbf{Table 44}  $ q = 23$, non-singular cubic curves giving complete $(n,3)$-arcs  \medskip \\
\medskip
\begin{tabular}
{|c|c|c|c|}
\hline
type & equation &$n$ &
\begin{minipage}{1 cm}
\smallskip
\centering
$j-$ \\
inv. \\
\smallskip
\end{minipage}  \\
\hline
\begin{minipage}{0.7 cm}
\centering
\smallskip
$ \mathcal{H}_3^2$  \end{minipage}
&   \begin{minipage}{10 cm}
\centering
\smallskip
$ T^3 + 3YT^2 + 3Y^2T + 3XT^2 + Y^3 + 7XYT + 3XY^2 + 3X^2T + 3X^2Y + X^3 $ \end{minipage}
&  \begin{minipage}{0.7 cm}
\centering
\smallskip
$24$ \end{minipage}
&  \begin{minipage}{1 cm}
\centering
\smallskip
$3$ \end{minipage}
\\
\begin{minipage}{0.7 cm}
\centering
\smallskip
$ \mathcal{H}_3^2$  \end{minipage}
&   \begin{minipage}{10 cm}
\centering
\smallskip
$ 18T^3 + 8YT^2 + 8Y^2T + 8XT^2 + 18Y^3 + 17XYT + 8XY^2 + 8X^2T + 8X^2Y + 18X^3 $ \end{minipage}
&  \begin{minipage}{0.7 cm}
\centering
\smallskip
$24$ \end{minipage}
&  \begin{minipage}{1 cm}
\centering
\smallskip
$3$ \end{minipage}
\\
\begin{minipage}{0.7 cm}
\centering
\smallskip
$ \mathcal{G}_3^2$  \end{minipage}
&   \begin{minipage}{10 cm}
\centering
\smallskip
$ 21T^3 + 17YT^2 + 17Y^2T + 17XT^2 + 21Y^3 + 12XYT + 17XY^2 + 17X^2T + 17X^2Y + 21X^3 $ \end{minipage}
&  \begin{minipage}{0.7 cm}
\centering
\smallskip
$33$ \end{minipage}
&  \begin{minipage}{1 cm}
\centering
\smallskip
$7$ \end{minipage}
\\
\begin{minipage}{0.7 cm}
\centering
\smallskip
$ \mathcal{G}_3^2$  \end{minipage}
&   \begin{minipage}{10 cm}
\centering
\smallskip
$ 6T^3 + 18YT^2 + 18Y^2T + 18XT^2 + 6Y^3 + 14XYT + 18XY^2 + 18X^2T + 18X^2Y + 6X^3 $ \end{minipage}
&  \begin{minipage}{0.7 cm}
\centering
\smallskip
$30$ \end{minipage}
&  \begin{minipage}{1 cm}
\centering
\smallskip
$18$ \end{minipage}
\\
\begin{minipage}{0.7 cm}
\centering
\smallskip
$ \mathcal{G}_3^2$  \end{minipage}
&   \begin{minipage}{10 cm}
\centering
\smallskip
$ 3T^3 + 9YT^2 + 9Y^2T + 9XT^2 + 3Y^3 + 19XYT + 9XY^2 + 9X^2T + 9X^2Y + 3X^3 $ \end{minipage}
&  \begin{minipage}{0.7 cm}
\centering
\smallskip
$30$ \end{minipage}
&  \begin{minipage}{1 cm}
\centering
\smallskip
$20$ \end{minipage}
\\
\hline
\end{tabular}
\end{minipage}\newline\newline\newline\newline\noindent
\begin{minipage}{\textwidth}
\centering
\textbf{Table 44} Continue $1$  \medskip \\
\medskip
\begin{tabular}
{|c|c|c|c|}
\hline
type & equation &$n$ &
\begin{minipage}{1 cm}
\smallskip
\centering
$j-$ \\
inv. \\
\smallskip
\end{minipage}  \\
\hline
\begin{minipage}{0.7 cm}
\centering
\smallskip
$ \mathcal{G}_1$  \end{minipage}
&   \begin{minipage}{10 cm}
\centering
\smallskip
$ YT^2 + 15Y^3 + 5XY^2 + X^3 $ \end{minipage}
&  \begin{minipage}{0.7 cm}
\centering
\smallskip
$32$ \end{minipage}
&  \begin{minipage}{1 cm}
\centering
\smallskip
$6$ \end{minipage}
\\
\begin{minipage}{0.7 cm}
\centering
\smallskip
$ \mathcal{G}_1$  \end{minipage}
&   \begin{minipage}{10 cm}
\centering
\smallskip
$ YT^2 + 13Y^3 + XY^2 + X^3 $ \end{minipage}
&  \begin{minipage}{0.7 cm}
\centering
\smallskip
$32$ \end{minipage}
&  \begin{minipage}{1 cm}
\centering
\smallskip
$8$ \end{minipage}
\\
\begin{minipage}{0.7 cm}
\centering
\smallskip
$ \mathcal{G}_1$  \end{minipage}
&   \begin{minipage}{10 cm}
\centering
\smallskip
$ YT^2 + 22Y^3 + 5XY^2 + X^3 $ \end{minipage}
&  \begin{minipage}{0.7 cm}
\centering
\smallskip
$31$ \end{minipage}
&  \begin{minipage}{1 cm}
\centering
\smallskip
$9$ \end{minipage}
\\
\begin{minipage}{0.7 cm}
\centering
\smallskip
$ \mathcal{G}_1$  \end{minipage}
&   \begin{minipage}{10 cm}
\centering
\smallskip
$ YT^2 + 19Y^3 + XY^2 + X^3 $ \end{minipage}
&  \begin{minipage}{0.7 cm}
\centering
\smallskip
$29$ \end{minipage}
&  \begin{minipage}{1 cm}
\centering
\smallskip
$11$ \end{minipage}
\\
\begin{minipage}{0.7 cm}
\centering
\smallskip
$ \mathcal{G}_1$  \end{minipage}
&   \begin{minipage}{10 cm}
\centering
\smallskip
$ YT^2 + 9Y^3 + XY^2 + X^3 $ \end{minipage}
&  \begin{minipage}{0.7 cm}
\centering
\smallskip
$28$ \end{minipage}
&  \begin{minipage}{1 cm}
\centering
\smallskip
$2$ \end{minipage}
\\
\hline
\end{tabular}%
\end{minipage}\newline


$\; $ \phantom  \newline\newline\newline\newline\noindent
\begin{minipage}{\textwidth}
\centering
\textbf{Table 46} $ q = 25$, a. i. singular cubic curves giving complete $(n,3)$-arcs  \medskip \\
\medskip
%
\end{minipage}\newline\phantom  \newline\newline


\noindent\begin{minipage}{\textwidth}
\centering
\textbf{Table 48}  $ q = 25$, non-singular cubic curves giving complete $(n,3)$-arcs  \medskip \\
\medskip
%
\end{minipage}\newline\phantom  \newline\newline

\noindent\begin{minipage}{\textwidth}
\centering
\textbf{Table 51}  $ q = 27$, non-singular cubic curves giving incomplete $(n,3)$-arcs  \medskip \\
\medskip
%
\end{minipage}\newline\phantom  \newline\newline


\noindent\begin{minipage}{\textwidth}
\centering
\textbf{Table 60}  $ q = 31$, non-singular cubic curves giving complete $(n,3)$-arcs  \medskip \\
\medskip
%
\end{minipage}\newline

\bigskip\bigskip\noindent\begin{minipage}{\textwidth}
\centering
\textbf{Table 61} $ q = 32$, a. i. singular cubic curves giving incomplete $(n,3)$-arcs  \medskip \\
\medskip
%
\end{minipage}\newline


$\; $ \phantom  \newline\newline\noindent\begin{minipage}{\textwidth}
\centering
\textbf{Table 62} $ q = 32$, a. i. singular cubic curves giving complete $(n,3)$-arcs  \medskip \\
\medskip
%
\end{minipage}\newline\phantom  \newline\newline


\noindent\begin{minipage}{\textwidth}
\centering
\textbf{Table 63}  $ q = 32$, non-singular cubic curves giving incomplete $(n,3)$-arcs  \medskip \\
\medskip
%
\end{minipage}\newline


$\; $ \phantom  \newline\newline\noindent\begin{minipage}{\textwidth}
\centering
\textbf{Table 66} $ q = 37$, a. i. singular cubic curves giving complete $(n,3)$-arcs  \medskip \\
\medskip
%
\end{minipage}\newline


$\; $ \phantom  \newline\newline\newline\noindent
\begin{minipage}{\textwidth}
\centering
\textbf{Table 78} $ q = 47$, a. i. singular cubic curves giving complete $(n,3)$-arcs  \medskip \\
\medskip
%
\end{minipage}\newline\phantom  \newline\newline\newline

\noindent\begin{minipage}{\textwidth}
\centering
\textbf{Table 79}  $ q = 47$, non-singular cubic curves giving incomplete $(n,3)$-arcs  \medskip \\
\medskip
%
\end{minipage}\newline


$\; $ \phantom  \newline\newline\noindent\begin{minipage}{\textwidth}
\centering
\textbf{Table 108} $ q = 73$, a. i. singular cubic curves giving complete $(n,3)$-arcs  \medskip \\
\medskip
%
\end{minipage}\newline\newline\newline\newline


\section{Appendix: Errata for Projective Geometries over Finite Fields,
Chapter 11}


Performing our computations, we found some typos in \cite[Chapter 11]{hirsh}.

\centering
\medskip%
%


\end{document}